\documentclass[]{amsart}
\usepackage[utf8]{inputenc}
\usepackage[document]{ragged2e}
\usepackage{amsmath,amssymb,amsthm,amsbsy}
\usepackage{epsfig}
\usepackage{graphicx}
\usepackage{relsize}
\usepackage{mdframed}
\usepackage{multirow}
\usepackage[width=1.0\textwidth, font=footnotesize]{caption}
\usepackage{mathtools}
\usepackage{adjustbox}
\usepackage{tcolorbox}
\usepackage{ragged2e}
\usepackage[linesnumbered,ruled,vlined]{algorithm2e}
\SetKwInput{KwInput}{Input}                    
\SetKwInput{KwOutput}{Output}              
\usepackage{hyperref}
\usepackage{subcaption}
\usepackage{bbm}
\usepackage{wrapfig}
\hypersetup{
	colorlinks=true,
	linkcolor=blue,
	citecolor=red,
	filecolor=red,      
	urlcolor=cyan,
}

\DeclareMathOperator*{\argmax}{argmax}

\newtheorem{remark}{Remark}

\newtheorem{example}{Example}

\newcommand{\abs}[1]{\left| #1 \right|}

\newcommand{\prt}[1]{\left( #1 \right)}
\newcommand{\cost}[1]{\mathrm{Cost}}

\newcommand{\bN}{\mathbb{N}}
\newcommand{\bP}{\mathbb{P}}
\newcommand{\bR}{\mathbb{R}}

\newcommand{\bT}{\mathbf{T}}

\newcommand{\cF}{\mathcal{F}}
\newcommand{\cO}{\mathcal{O}}

\newcommand{\cK}{\mathcal{K}}
\newcommand{\cL}{\mathcal{L}}

\newcommand{\cT}{\mathcal{T}}
\newcommand{\cU}{\mathcal{U}}

\newcommand{\cV}{\mathcal{V}}

\newcommand{\bu}{\boldsymbol{\upsilon}}
\newcommand{\bv}{\boldsymbol{\nu}}
\newcommand{\bw}{\boldsymbol{w}}

\newcommand{\bzeta}{\boldsymbol{\zeta}}

\newcommand{\vHat}{\hat{v}}

\newcommand{\E}{\mathbb{E}}
\newcommand{\V}{\mathbb{V}}

\newcommand{\Ind}[2]{\mathbbm{1}_{\{#2\}}}
\newcommand{\IndU}[1]{\mathbbm{1}_{\cU} \left(#1 \right)}

\newcommand{\red}[1]{{\color{red} #1}}

\title[Tracking rare events within EnKF]{Importance sampling for rare event tracking within the ensemble Kalman filtering framework}

\author[N. Ben Rached]{Nadhir Ben Rached} \address[{Nadhir Ben Rached}]{\newline
   School of Mathematics, Faculty of Engineering and Physical Sciences, University of Leeds, UK \newline (n.benrached@leeds.ac.uk)}

\author[E. von Schwerin]{Erik von Schwerin} \address[{Erik von Schwerin}]{\newline Applied Mathematics and Computational
	Sciences, KAUST, Thuwal, Saudi Arabia \newline
	(erik.vonschwerin@kaust.edu.sa)}

\author[G. Shaimerdenova]{Gaukhar
	Shaimerdenova$^*$}\thanks{$^*$Corresponding author: G.Shaimerdenova
	(gaukhar.shaimerdenova@kaust.edu.sa)} \address[{Gaukhar
	Shaimerdenova}]{\newline Applied Mathematics and Computational
	Sciences, KAUST, Thuwal, Saudi Arabia \newline
	(gaukhar.shaimerdenova@kaust.edu.sa)}

\author[R. Tempone]{Ra\'ul Tempone} \address[{Raul Tempone}]{\newline
	Chair of Mathematics for Uncertainty Quantification, RWTH Aachen
	University, Aachen, Germany \newline (tempone@uq.rwth-aachen.de)
	\newline
	\and
	\newline
	Applied Mathematics and Computational Sciences, KAUST, Thuwal, Saudi Arabia \newline (raul.tempone@kaust.edu.sa)}

\begin{document}

\begin{abstract}
\begin{justify}

In this work we employ importance sampling (IS) techniques to track a small over-threshold probability of a running maximum associated with the solution of a stochastic differential equation (SDE) within the framework of ensemble Kalman filtering (EnKF). Between two observation times of the EnKF, we propose to use IS with respect to the initial condition of the SDE, IS with respect to the Wiener process via a stochastic optimal control formulation, and combined IS with respect to both initial condition and Wiener process. Both IS strategies require the approximation of the solution of Kolmogorov Backward equation (KBE) with boundary conditions. In multidimensional settings, we employ a Markovian projection dimension reduction technique to obtain an approximation of the solution of the KBE by just solving a one dimensional PDE. The proposed ideas are tested on three illustrative examples: Double Well SDE, Langevin dynamics and noisy Charney-deVore model, and showcase a significant variance reduction compared to the standard Monte Carlo method and another sampling-based IS technique, namely, multilevel cross entropy. 
\end{justify}
  \bigskip
  \noindent
  
  
  
\end{abstract}

\maketitle
\begin{justify}
\section{Introduction}

Let $\{u_t\}_{t\geq 0}\in\bR^d$ be a Markov process which determines the state of the system arising from a stochastic dynamics model. In particlular, we assume the model dynamics is associated with a stochastic differental equation (SDE) driven by the Wiener process. Suppose also we are given partially observed observations $\{y_n\}_{n\in\bN}$ arriving sequentially at discrete time instances $\{t_n\}_{n\in\bN}$. We wish to approximate the expected value of a quantity of interest $\varphi: \bR^d \rightarrow \bR$ applied to $u_n:=u_{t_n}$, taking into account the accumulated observation data $Y_n:=(y_1, y_2, ..., y_n)$ up to time $t_n$, i.e., we seek to approximate $\E [\varphi (u_{n})|Y_n]$. 

Data assimilation is a powerful methodology to combine the observation data with the model forecast to produce more reliable estimates of true states of the system. There are various data assimilation techniques with their strengths and limitations. Among them, the ensemble Kalman filter (EnKF) has quickly become popular for performing inference online in many applications due to its simplicity in implementation and efficiency in high dimensional problems~\cite{evensen1994sequential, kalnay2003atmospheric, aanonsen2009ensemble, houtekamer2005atmospheric, reich2015probabilistic}. However, the EnKF is impractical when the system exhibits very unlikely, unusual fluctuations, i.e., in failure/rare event occurrences. Suppose we want to estimate the probability of a failure event of a dynamic system between two observation times $t_n$ and $t_{n+1}$ before obtaining the next measurement at time $t_{n+1}$. For example, we may want to detect a hazard for an aircraft caused by extreme turbulence given radar observations~\cite{gao2016distributed}, predict a financial crisis using stock market data~\cite{karmiani2019comparison}, or forecast a severe weather event provided different sources of weather measurements~\cite{jones2015simultaneous, zhang2015capturing, maejima2020impact}.

\textit{Objective:}  In this work, we are concerned with a particular type of rare event. More precisely, we want to track the probability of the running maximum of the projected one-dimensional stochastic process exceeding a critical threshold $\cK>0$ until the next observation time $t_{n+1}$ within the framework of the EnKF given the data $Y_n$:
\[
\alpha_n:=\bP{\Big(\max_{t_n\leq t\leq t_{n+1}} P_1 u_t \geq \cK \Big|Y_n\Big)}, 
\]
where $P_1\in \bR^{1\times d}$ is a pre-defined projection applied to $u_t\in \bR^d$. 
\begin{figure}[h!]
	\includegraphics[height=6cm, width=7cm]{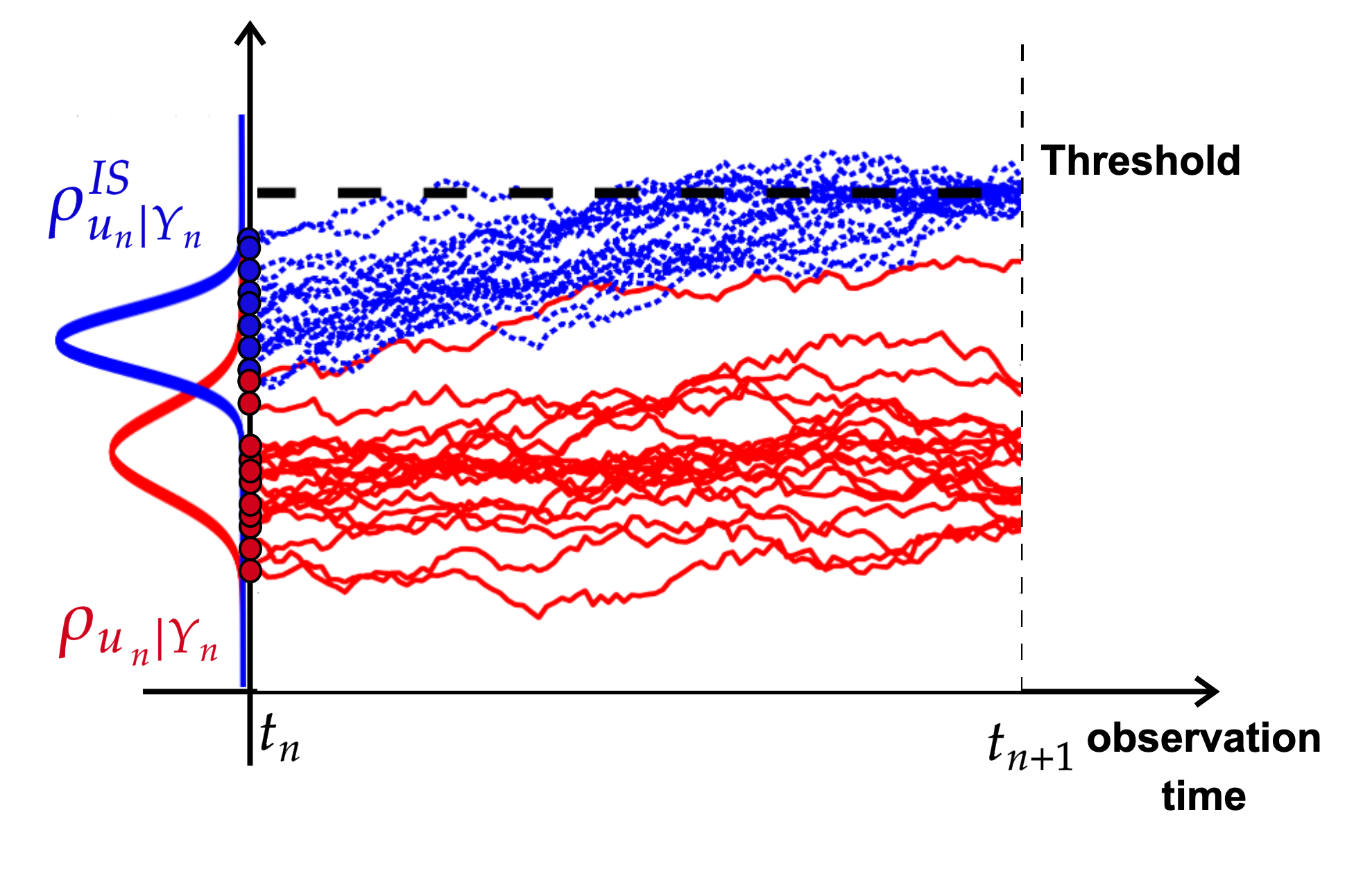}
	\caption{Illustration of the IS idea for the rare event tracking in the context of the EnKF.}
	\label{fig:ideaIS}
\end{figure}

Since the EnKF is a Monte Carlo (MC)-based filtering method, the estimation of a very small over-threshold probability by EnKF ensembles is computationally expensive for rare event regimes $\alpha_n \ll 1$. Various variance reduction techniques exist to improve the efficiency of the MC simulation. One famous technique to reduce the variance is importance sampling (IS). It is based on the idea of changing the original probability density to another biasing density such that the rare event is more likely to occur~\cite{kroese2013handbook, asmussen2007stochastic,benrachedIS}. To ensure the unbiasedness of the new IS estimator, one must reweight each sample by an appropriate likelihood ratio. However, a key challenge in IS is to select a good choice of biasing density, which results in a significantly smaller variance of the estimator compared to the one provided by the standard MC. 
 
 In this work, we incorporate IS techniques with the EnKF method to be able to track the rare event probability in real-time. We propose the following approaches that can be applied at each observation time $t_n$ within the framework of EnKF: IS with respect to (wrt) the initial condition of the SDE (i.e., a change of measure from the original density $\rho_{u_n|Y_n}$ to the biasing IS density $\rho_{u_n|Y_n}^{IS}$); IS wrt the Wiener process via a stochastic optimal control (SOC) problem~\cite{hartmann2019variational}; and the combination of IS wrt both initial condition and Wiener process. The idea is illustrated in Figure~\ref{fig:ideaIS}, where the red circles form the standard EnKF ensemble and the red trajectories are induced by the Wiener paths, while the blue circles are the auxiliary EnKF ensemble to estimate the rare event probability and the blue trajectories are driven by the controlled SDE obtained via change of measure wrt the Wiener process, thus, allowing to monitor the failure probability between observation times. Here, we note that the auxiliary EnKF ensemble size, which depends on the rareness of the estimated event, is not necessarily the same as the original EnKF ensemble size. The proposed IS schemes are based on solving the Kolmogorov backward equation (KBE), which is not computationally feasible in high dimensions. To overcome this issue, we employ the Markovian projection technique~\cite{gyongy1986mimicking,bayer2019implied} that allows us to approximate the multidimensional KBE by a lower-dimensional KBE. We compare our proposed approaches with the multilevel cross entropy (CE) method which is another widely-used technique to estimate rare event probabilities~\cite{kroese2013cross,kroese2013handbook}. 
 
\textit{Literature review.} An accurate estimation of rare event probabilities based on dynamic systems is a challenging problem that usually lacks analytical solutions and needs both sophisticated theory and efficient computational tools~\cite{bucklew2004introduction, beck2015rare, sapsis2021statistics}. A considerable amount of literature has been dedicated to analytical and semi-analytical approximation methods from the structural reliability analysis, for instance, the first and second order reliability methods (FORM and SORM) based on first and second order Taylor series expansions, respectively, to find
the most likely failure point~\cite{melchers2018structural, hasofer1974exact, der2005first}. However, these methods can be computationally costly in high-dimensional settings since they involve an optimization problem to solve and exhibit large approximation errors of estimates~\cite{valdebenito2010role}. An additional perspective to tackle the problem is based on the theory of large deviations which indicates that, in cases of low noise, trajectories may converge towards the most likely path leading to a rare event. The large deviation principle allows to estimate the rare event probability by solving a deterministic optimization problem based simply on the ordinary differential equations. However, this is possible only in the context of the small noise limit. The estimation of extreme events within the framework of large deviations is reviewed in~\cite{grafke2019numerical}. A combination of large deviations theory with the PDE-constrained optimization is proposed in~\cite{tong2021extreme}, with the importance sampling is in~\cite{ragone2018computation}, and with the cross entropy-based importance sampling is in~\cite{stadlerLDadaptive}. 

An alternative approach that has garnered significant attention and is being studied extensively is the use of sampling-based methodologies for rare events in dynamic models. Since a crude Monte Carlo method is unsuitable for the simulation of small probabilities~\cite{rubinstein2016simulation}, different advanced variance reduction techniques have been developed to speed up its efficiency while the accuracy is maintained (e.g., subset simulation/splitting algorithms~\cite{au2001estimation, botev2012efficient, cerou2019adaptive, brehier2019new, webber2019practical}, sequential MC and IS~\cite{heidelberger1995fast, cerou2012sequential, papaioannou2016sequential}, line sampling~\cite{de2015advanced,straub2015bayesian}, cross-entropy~\cite{kroese2013cross,de2005tutorial} methods). In general, sampling methods have demonstrated their strength and efficiency in estimating rare events accurately for high-dimensional settings, as evidenced by their successful implementation in various applications~\cite{au2003subset,beaurepaire2012reliability,zio2009estimation}. Another effective sampling method designed for rare events is referred to as genealogical particle analysis. This approach involves selecting and replicating the most favorable realizations of an ensemble of trajectories in an iterative form similar to the idea of MC acceptance-rejection techniques~\cite{del2005genealogical}. Recent studies have also explored the relationship between IS schemes and SOC theory when the underlying dynamics is governed by diffusion~\cite{hartmann2016model, hartmann2017variational, hartmann2018importance} with the further application for rare event simulation problems~\cite{hartmann2019variational}. The latter is based on a Gibbs variational principle to find the optimal change of measure, which can be formulated as a SOC problem and solved by the Feynman–Kac representation of the dynamic programming equations~\cite{hartmann2019variational}. The SOC theory is also applied for estimating rare event probabilities in the case of the McKean-Vlasov SDEs~\cite{ben2022double, rached2022multilevel} and stochastic reaction networks~\cite{ben2023learning} with further application of the Markovian projection technique in~\cite{hammouda2023automated}.

In the context of data assimilation, the approaches to tracking rare events proposed in this article align with the study presented in~\cite{vanden2013data}. The authors combine features of both sampling-based filters and variational methods, resulting in the development of hybrid data assimilation techniques, and test their idea on the Kuroshio model. The study is focused only on the low noise regime and with a small stochastic forcing term, which allows the application of large deviation analysis. The work \cite{wagner2022ensemble} employed the EnKF as a tool to estimate failure probabilities where the rare event problem is reformulated 
as a Bayesian inverse problem. In particular, the auxiliary EnKF ensembles are used to generate failure samples but not for tracking rare events in real-time, which is the focus of the current study.

\textit{Outline.} The subsequent sections of the work are structured in the following manner. Section~\ref{sec:problemsetup} describes a discrete-time filtering problem, followed by a concise summary of the EnKF framework. It also includes a discussion of how the problem relates to the concept of stopping time and the numerical approximation of the quantity of interest (QoI). In Section~\ref{sec:IStechniques}, we propose several IS techniques that can be combined with the EnKF framework. Finally, Section~\ref{sec:numerics} exhibits the numerical part where the proposed IS schemes are applied to three different examples: Double Well SDE, Langevin dynamics and noisy Charney-deVore model. 

\section{Problem setting}\label{sec:problemsetup}

We first describe the discrete-time filtering problem with observations arriving at times $\{t_n\}_{n\in \bN}$, and later return to the problem of estimating rare event probabilities between observations. 

Let the triple $\prt{\Omega, \cF, \bP}$ denote a complete
probability space on which we consider the filtering problem. The primary goal in the filtering problem is to track the underlying \textit{signal} $u_n$ governed by the stochastic dynamics 
\begin{equation}\label{eq:dynamics}
u_{n}(\omega)=\Psi (u_{n-1}, \omega), \quad n=1,2,...
\end{equation}
provided the accumulated data $Y_{n}$ up to time $t_n$ with \textit{noisy observations} given by 
\begin{equation}\label{eq:observations}
y_{n}(\omega) = Hu_{n}(\omega)+\eta_{n}(\omega),  \quad n=1,2,...
\end{equation}
where $\omega\in \Omega$, $H\in \bR^{m\times d}$, and $\{\eta_k\}_{k\in\bN}$ is an independent and identically distributed (i.i.d.) sequence being independent of $\{u_k\}_{k\in\bN \cup \{0\}}$ with $\eta_1 \sim N(0,\Gamma)$ for a positive definite $\Gamma \in \bR^{m\times m}$. For the sake of simplicity, we further drop the dependence on $\omega$ in the notation throughout the paper if there is no confusion. 

In other words, the filtering problem is concerned with an approximation of the conditional random variable $u_{n}|Y_{n}$ in real-time as soon as the data at time $t_n$ is acquired.  Let $\rho_{u_n|Y_n}$ denote the probability density function of $u_n$ given $Y_n$ and $\rho_{u_n|u_{n-1}}$ denote the transition density of $u_n$ given the previous value $u_{n-1}$ typically taken from the model equation in~\eqref{eq:dynamics}. Based on Bayes' theorem, the filtering problem can be expressed formally via two recursive steps of computing the probability density of $u_{n}|Y_{n}$ as follows
\[
\begin{split}
	\mbox{ \textbf{prediction step:}  }& \quad \rho_{u_{n}|Y_{n-1}}(u)
	\propto \int_{\bR^d}\rho_{u_{n}|u_{n-1}}(u) \rho_{u_{n-1}|Y_{n-1}}(v)dv,\\
	\mbox{ \textbf{update step:}  }& \quad \rho_{u_{n}|Y_{n}}(u)
	\propto  \exp\big(-\big| \Gamma^{-1/2}(y_{n}-Hu)\big|^2/2 \big) \; \rho_{u_{n}|Y_{n-1}}(u).
\end{split}
\]

If the initial density of $u_0|Y_0$ is assumed to be Gaussian and the dynamics $\Psi$ is linear with additive Gaussian noise, the filtering density becomes Gaussian, which can be entirely described by only two statistics (mean and covariance). In this case, the Kalman filter provides an exact algorithm to track the mean and covariance at each observation timestep. However, in a more general case, when $\Psi$ is nonlinear, an exact solution to the filtering problem is usually not attainable, and we mostly rely on approximation methods. In this work, we focus on an important class of such approximation methods in filtering, namely, the EnKF.

The EnKF is an extension of the Kalman filter to nonlinear settings, which produces an ensemble of particles whose empirical measure approximates the true filtering distribution of $u_{n}|Y_{n}$. Let $P$ denote the EnKF ensemble size and the pair  $\prt{v_{n,i}, \vHat_{n,i}}$ denote the prediction and updated
ensemble particles, respectively, corresponding to the sample $\omega_i$ at time $t_n$. Then, the EnKF algorithm with perturbed
observations consists of two iterative
steps for $n\in \bN$, namely
\begin{equation} \label{enkf:prediction}
	\centering \mbox{ \textbf{prediction step:}  }
	\left\{\begin{split}
		v_{n,i} &= \Psi_{n-1}^N(\vHat_{n-1,i}), \quad i=1,2,...,P,\\
		m_{n} &= \frac{1}{P} \sum_{i=1}^P v_{n,i},\\
		C_{n} &=\frac{1}{P-1}\sum_{i=1}^P \prt{v_{n,i}-m_{n}} \prt{v_{n,i}-m_{n}}^\bT,
	\end{split}\right.
\end{equation}
where $\Psi_{n-1}^N$ denotes the numerical discretization of the dynamics $\Psi_{n-1}$ using $N\geq1$ uniform timesteps, and
\begin{equation}\label{enkf:update}
	\centering \mbox{ \textbf{update step:}  }
	\left\{\begin{split}
		\tilde{y}_{n,i}&=y_{n}+\eta_{n,i}, \; i=1,2,...,P,\\
		\mathrm{K}_{n}&=C_{n}H^{\bT}(HC_{n}H^{\bT}+\Gamma)^{-1},\\
		\vHat_{n,i} &= (I-\mathrm{K}_{n}H)v_{n,i}+\mathrm{K}_{n}\tilde{y}_{n,i},\; i=1,2,...,P,\\
	\end{split}\right.
\end{equation}
where $\eta_{n,i}$ are i.i.d.~draws from $N(0, \Gamma)$, the matrix $\mathrm{K}_n$ is known as the Kalman gain and $H$ is the observation matrix in~\eqref{eq:observations} \cite{hoel2020multilevel,hoel2022multi}.


Let us consider discrete observation time points over the interval $[0,T]$, denoted by $0=t_0<t_1<...<t_O=T$. In the remainder of this work, we will consider low probability events between two consecutive observations, at times $t_n$ and $t_{n+1}$. As mentioned earlier, it will be assumed that the underlying signal process is considered as the evolution of the dynamics $\Psi$ associated with an SDE of the form
\begin{equation}\label{eq:gensde}
	\left\{\begin{split}
  &du_t=a(u_t)dt+b(u_t)dW_t^{\bP}, \quad t\in (t_n, t_{n+1}),\\
  &u_{n} \sim \rho_0,
	\end{split}\right.
\end{equation}
where $a: \bR^d \rightarrow \bR^d$ is called the drift coefficient, $b: \bR^d \rightarrow \bR^{d\times d_W}$ is called the diffusion coefficient, $W_t^{\bP}:\Omega \times [t_n, \infty) \rightarrow \bR^{d_W \times 1}$ is an $d_W$-dimensional Wiener process, independent of $u_k, \; k\leq n$, under the probability measure $\bP$, and the initial condition $\rho_0$ is defined by the updated EnKF density $\hat{\rho}_{u_n|Y_n}$ at time $t_n$. The drift and diffusion coefficients are generally expected to satisfy the regularity conditions to ensure the existence and uniqueness of the solution to~\eqref{eq:gensde} (see, e.g.~\cite[Section 4.5]{kloeden1992stochastic} for sufficient conditions). 

\begin{remark}[Tail assumptions] Note that this work does not consider more challenging cases such as SDEs resulting in heavy-tailed events. For that purpose, we assume that the drift coefficient $a$ is measurable and has a linear growth in space, and the diffusion coefficient $b$ is bounded and satisfies the usual uniform ellipticity condition. This ensures that the SDE~\eqref{eq:gensde} has a unique weak solution and admits the density whose tails decay as Gaussian (see, e.g.~\cite[Eq.~(1.2)-(1.3)]{menozzi2021density}). 
\end{remark}

According to the main objective of the work, we restrict our attention to the observation interval $[t_n, t_{n+1}]$ (although the proposed approaches in this work will be applicable for all observation time intervals). 
We define the first time when the running maximum process $\{\max_{t_n\leq t\leq t_{n+1}} P_1 u_t\}$ hits the threshold $\cK$ given $P_1u_n<\cK$ at the initial time $t_n$ by
\begin{equation*}
	\tau_{\cK} := \inf\{t>t_n: \; P_1u_t\geq \cK \, | \, P_1u_n<\cK \}. 
\end{equation*}
In the literature, the random variable $\tau_\cK$ is referred to as \textit{a first exit/passage time} and $\tau:=\min\{\tau_{\cK}, t_{n+1}\}$ is as \textit{a stopping time}. There is a connection between the events $\{\max_{t_n\leq t\leq t_{n+1}} P_1 u_t \geq \cK\}$ and $\{\tau_\cK \leq t_{n+1}\}$, which allows us to write 
\begin{equation}\label{eq:stoppingtime}
\alpha_n=\E_{\rho_0 \otimes \mathbb{P}}[\Ind_{\max_{t_n\leq t\leq t_{n+1}} P_1 u_t\geq \cK}|Y_n]=\E_{\rho_0 \otimes \mathbb{P}}[\Ind_{P_1u_{\tau}\geq \cK}|Y_n],
\end{equation}
where $\Ind{}{\cdot}$ is an indicator function. Furthermore, for convenience, we will omit the notation for dependence on the data $Y_n$.

A standard procedure to estimate the QoI in~\eqref{eq:stoppingtime} is to employ a time-discretization scheme such as Euler-Maruyama and use an MC simulation for computing the mean estimation of the exit events.  Let $t_n=t_{n,0}<t_{n,1}<...<t_{n,k}<...<t_{n,K}=t_{n+1}$ be a partition of the interval $[t_n, t_{n+1}]$ with $t_{n,k}:=t_n+k\Delta t$ and $\Delta t= \frac{t_{n+1}-t_n}{K}$. The discrete forward Euler-Maruyama (EM) approximation of the SDE~\eqref{eq:gensde} is 
\begin{equation}\label{eq:EMscheme}
	\begin{cases} 
	&\bar{u}_{n,k+1}=\bar{u}_{n,k}+a(\bar{u}_{n,k})\Delta t+b(\bar{u}_{n,k})\Delta W^{\bP}_{n,k},\\
	&\bar{u}_{n,0}=u_n,
	\end{cases}\,
\end{equation}
for $k=0,...,K-1$, with $\Delta W^{\bP}_{n,k}=W^{\bP}_{n,k+1}-W^{\bP}_{n,k} \sim N(0,\Delta t \mathbb{I}_{d_W})$ where $ \mathbb{I}_{d_W}$ is the $d_W$-dimensional identity matrix.

For sufficiently small $\Delta t$, we can approximate the QoI by a function of $K$ discrete EM points, i.e.,
\begin{equation}\label{eq:QoIapproxbyEM}
	\begin{split}
		\alpha_n \approx \E_{\rho _0\otimes \bP }[ \IndU{\bar{u}_{n,0:K}}],
	\end{split}
\end{equation}
where $\bar{u}_{n,0:K}:=\{\bar{u}_{n,k}\}_{k=0}^K$ and 
$$\cU:=\{x_k\in\bR^d, k=0,...,K: P_1x_k\geq \cK \mbox{  for some  } k\}.$$ Hereafter, we work with this EM approximation of the QoI to perform IS techniques. 

 \begin{figure}[h!]
	\includegraphics[height=6cm, width=7cm]{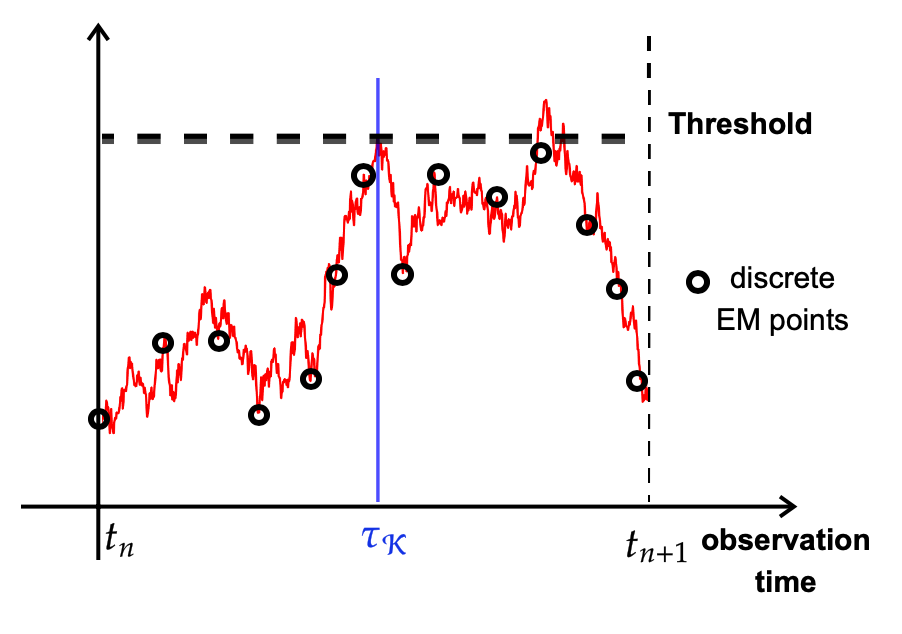}
	\caption{Illustration of the challenge with stopped diffusions in the boundary discussed in Remark~\ref{remark:brownianbridge}. The continuous path (red) hits the threshold at $\tau_\cK$ although the discrete EM solutions do not exit the boundary (or may exit after $\tau_\cK$).}
	\label{fig:stoppingtime}
\end{figure}
 
\begin{remark}[Stopped diffusions] \label{remark:brownianbridge}
One of the challenges in approximating the QoI involving stopped diffusions on the boundary is that a continuous path may leave the domain (i.e., hit the threshold) even if a discrete approximate solution does not exit the boundary (see Figure~\ref{fig:stoppingtime}). This leads to a slow convergence of the time-discretization error (for more discussions, see~\cite{gobet2000weak, gobet2009advanced, dzougoutov2005adaptive}). In particular, the EM scheme in~\eqref{eq:EMscheme} achieves a weak convergence, $\cO(\sqrt{\Delta t})$, of the relative error, as $\Delta t \rightarrow 0$, compared to a standard Euler-type convergence $\cO(\Delta t)$ under certain conditions~\cite{kloeden1992stochastic, bally1996law,higham2001algorithmic}. One way to recover the first-order convergence of the EM method is by performing a straightforward hitting test after each time step. This involves examining the distribution of the Brownian bridge between two discrete EM nodes to determine if the threshold was crossed during the given time step. This approach is referred to as \textit{the Brownian bridge technique}~\cite{gobet2000weak} and implemented in this work (see Algorithm~\ref{alg:adaptivetimestep} in Appendix~\ref{appx:AdaptiveScheme}). However, note that this is not the main challenge considered in this work.
\end{remark}

\section{Importance Sampling techniques}
\label{sec:IStechniques}
The primary focus of the work is to reduce the variance effectively via appropriate importance sampling techniques. For that purpose, observe that there are two sources of randomness in~\eqref{eq:gensde} coming from the initial condition and also from the Wiener process. Correspondingly, we can perform IS wrt $\rho_0$, IS wrt $W_t^\bP$ and IS wrt both $\rho_0$ and $W_t^\bP$. In this section, we describe the methodology of each proposal applied to the EM approximation of the QoI. 

\subsection{IS with respect to initial condition}
\label{ssec:ISrho0}
Let us consider the joint density $ \rho_{\bar{u}_{n,0:K}}$ of a random vector $(\bar{u}_{n,0:K})$ with $\rho_{\bar{u}_{n,0}}=\rho_0$. Now we introduce an alternative density function $\tilde{\rho}_{0}$, and using the shorthand $x_{1:K}$ for a $K$-tuple $(x_1,...,x_K)$, rewrite the approximation of the QoI as 

\begin{equation}\label{eq:QoIchangeinitial}
	\begin{split}
	\hat{\alpha}_n:=	\E_{\rho_0 \otimes \bP}[\IndU{\bar{u}_{n,0:K}}]&=\int  \IndU{x_{0:K}} \rho_{\bar{u}_{n,0:K}} (x_{0:K})dx_{0:K}\\
		& =\int \IndU{x_{0:K}} \rho_{0} (x_0) \rho_{ \bar{u}_{n,1:K}|\bar{u}_{n,0}} (x_{1:K}|x_0) dx_{0:K}\\
		&= \int \IndU{x_{0:K}} \frac{\rho_{0} (x_0)}{\tilde{\rho}_{0}(x_0)} \rho_{ \bar{u}_{n,1:K}|\bar{u}_{n,0}} ( x_{1:K}|x_0)\tilde{\rho}_{0}(x_0) dx_{0:K}\\
		&=\E_{\tilde{\rho}_0 \otimes \bP} \Big[ \IndU{\bar{u}_{n,0:K}} L_0( \bar{u}_{n,0})\Big],
	\end{split}
\end{equation}
where $L_0(x_0):=\frac{\rho_{0}(x_0)}{\tilde{\rho}_{0}(x_0)}$ is \textit{the likelihood ratio}. 

In order to select a proper density $\tilde{\rho}_0 $ that results in a smaller variance, let us consider the second moment of the QoI approximation,
\begin{equation}\label{eq:secondmomentestim}
	\begin{split}
	\E_{\tilde{\rho}_0 \otimes \bP} \Big[\IndU{\bar{u}_{n,0:K}} L^2(\bar{u}_{n,0})\Big]&=\int \IndU{x_{0:K}} \frac{\rho_{0}^2 (x_0)}{\tilde{\rho}_{0}^2(x_0)} \rho_{ \bar{u}_{n,1:K}|\bar{u}_{n,0}} ( x_{1:K}|x_0)\tilde{\rho}_{0}(x_0) dx_{0:K}\\
		&=\int\IndU{x_{0:K}} \frac{\rho_{0}^2 (x_0)}{\tilde{\rho}_{0}(x_0)} \rho_{ \bar{u}_{n,1:K}|\bar{u}_{n,0}} ( x_{1:K}|x_0) dx_{0:K}\\
		&=\int \E[\IndU{\bar{u}_{n,0:K}}| \bar{u}_{n,0}=x_0]  \frac{\rho_{0}^2 (x_0)}{\tilde{\rho}_{0}(x_0)}  dx_0
	\end{split}
\end{equation}
where the last equality follows from Fubini's theorem. Now applying the Cauchy-Schwartz inequality to the last term, we have 
\begin{equation}\label{eq:Cauchy-Schwartz}
	\begin{split}
		&\int\prt{\sqrt{\E[\IndU{\bar{u}_{n,0:K}}| \bar{u}_{n,0}=x_0]}  \frac{\rho_{0} (x_0)}{\sqrt{\tilde{\rho}_{0}(x_0)}}}^2  dx_0 \int \prt{\sqrt{\tilde{\rho}_{0}(x_0)}}^2 dx_0 \\
		&\geq \prt{\int \abs{\sqrt{\E[\IndU{\bar{u}_{n,0:K}}| \bar{u}_{n,0}=x_0]}\rho_{0} (x_0)}}^2,
	\end{split}
\end{equation}
where we can deduce that the optimal IS initial density is proportional to
\begin{equation}\label{eq:optimdensityPropto}
	\begin{split}
	\tilde{\rho}_{0}(x_0) &\propto \rho_{0}(x_0) \sqrt{ \E[ \IndU{\bar{u}_{n,0:K}}|  \bar{u}_{n,0}=x_0]}\\
	&= 	\rho_{0}(x_0) \sqrt{\bP\Big(\max_{0
			\leq k\leq K } P_1 \bar{u}_{n,k}\geq \cK \Big|\bar{u}_{n,0}=x_0\Big)}.
	\end{split}
\end{equation}
Here, $\rho_0$ is given and the exit probability needs to be approximated.

\subsubsection{\textbf{PDE-based method}}\label{ssec:PdeMethod} We propose a PDE-based approach to computing the probability in~\eqref{eq:optimdensityPropto}. Let us define the function
\[
\gamma(x,t):=\bP\Big(\max_{t\leq s\leq t_{n+1}} P_1 u_s \geq \cK \Big|  u_t =x\Big).
\]
This $\gamma(x,t)$ solves the Kolmogorov Backward Equation (KBE) associated with the SDE~\eqref{eq:dynamics}. We will first state this problem for $d=1$, where the standard discretization-based numerical PDE approximation methods are efficient, and then propose the use of Markovian projection~\cite{gyongy1986mimicking, bayer2019implied} as a way of approximating the solution of the $d$-dimensional PDE by a solution of a lower dimensional one.

\textit{One-dimensional setup.} For $d=1$, the process $P_1u_t=u_t$ is a one-dimensional Markov process, and thus, $\gamma(x,t)$ satisfies the KBE with a Dirichlet boundary condition on the domain $(-\infty, \cK] \times [t_n,t_{n+1}]$, 
\begin{equation}\label{eq:BKEpde}
\left\{
\begin{array}{ll}
	\frac{\partial \gamma}{\partial t}  = -a(x)	\frac{\partial \gamma}{\partial x}   -\frac{1}{2} b(x)b(x)^T 	\frac{\partial^2 \gamma}{\partial x},\\
	\gamma(x, t_{n+1}) =0, \quad \;\; x<\cK,\\
	\gamma(\cK, t) =1, \qquad \quad t\in [t_n, t_{n+1}],\\
	\lim\limits_{x\rightarrow -\infty} \gamma(x,t)=0, \;\, t\in [t_n, t_{n+1}].
\end{array}
\right.
\end{equation}
Note that there is a discontinuity in the data at $(\cK, t_{n+1})$. By numerically approximating~\eqref{eq:BKEpde} backward in time with the technique of freezing the coefficients in a neighborhood of the discontinuity (see Appendix~\ref{appx:PDEsolver} for more details about the numerical solver), we can obtain a sufficiently accurate approximation for the discrete-time exit probability in~\eqref{eq:optimdensityPropto} by a terminal time solution
\begin{equation*}
	   	\bP\Big(\max_{0
		\leq k\leq K } P_1 \bar{u}_{n,k}\geq \cK \Big|\bar{u}_{n,0}=x_0\Big) \approx 	\begin{cases} \gamma(x_0,t_n), \; &\mbox{ if } x_0< \cK,\\
		1, &\mbox{ if } x_0\geq\cK,
	\end{cases}
\end{equation*}
and compute for $x_0\in \bR$
\begin{equation}\label{eq:optimdensityPDE}
	\tilde{\rho}_{0}^{PDE,1}(x_0)\propto \rho_{0}(x_0) \sqrt{\gamma(x_0,t_n)}.
\end{equation}
For simplicity, due to the Gaussian assumption of EnKF, we fit a Gaussian density with parameters $(\mu_{0}^{\mathrm{fit}}, \sigma_{0}^{\mathrm{fit}})$ to the optimal IS initial density after normalizing the product $\rho_{0}(x_0) \sqrt{\gamma(x_0,t_n)}$ using a numerical quadrature rule for approximation of integrals.

Let us assume $\rho_0(x_0)\sim N(\mu_{0},\sigma_{0})$. Then, the likelihood $L_0$ is defined by the ratio of two Gaussian densities 
\[
L_0(x_0)=\frac{\rho_0(x_0)}{\tilde{\rho}_{0}^{PDE,1}(x_0)}=\frac{\sigma_{0}^{\mathrm{fit}}}{\sigma_{0}}\exp{\prt{-\frac{(x_0-\mu_{0})^2}{2\sigma_{0}^2}+\frac{\big(x_0-\mu_{0}^{\mathrm{fit}}\big)^2}{2\big(\sigma_{0}^{\mathrm{fit}}\big)^2}}}.
\]

\textit{High-dimensional setup.} For $d>1$,  the process $S_t:=P_1u_t$ with $u_t\in \bR^d$ is in general no longer Markovian, and thus, we apply a \textit{Markovian projection technique} to be able to employ the proposed PDE-based method above~\cite{gyongy1986mimicking, bayer2019implied}. We introduce a Markovian surrogate process $\check{S}_t\in \bR$ which follows the SDE
\begin{equation}\label{eq:MarkovSDE}
	\begin{cases} 
		&d\check{S}_t = \check{a}(\check{S}_t, t)dt+\check{b}(\check{S}_t, t)d\check{W}^{\bP}_t, \quad t\in (t_n, t_{n+1}),\\
		&\check{S}_{t_n}=P_1u_n, \quad u_n\sim \rho_0,
	\end{cases}\,
\end{equation}
where $\check{W}_t^{\bP}:\Omega \times [0, \infty) \rightarrow \bR$ is an one-dimensional Wiener process independent of the initial condition $\check{S}_{t_n}$, and $\check{a}, \check{b}: \bR\times [0, \infty) \rightarrow \bR$ are non-random coefficients defined by
\begin{equation}\label{eq:surrogateDriftAndDiff}
	\begin{cases} 
		\check{a}(y, t) = \E[P_1a(u_t)| P_1 u_t=y],\\
		\check{b}^2(y, t) = \E [(P_1 b b^T P_1^T)(u_t)| P_1 u_t=y].
	\end{cases}\,
\end{equation} 
It is proved that, in the case of a fixed initial condition, if $P_1a(x)$ and $P_1 b (x)$ are bounded measurable functions such that $(P_1 b b^T P_1^T)(x)$ is uniformly positive definite, then the SDE~\eqref{eq:MarkovSDE} admits a weak solution $\check{S}_t$ which has the same one-dimensional probability distribution as $S_t$ for every $t$~\cite{gyongy1986mimicking}. Since we have a random initial condition, it is assumed that the conditional expected values in~\eqref{eq:surrogateDriftAndDiff} exist and are well-defined too. One can use a discrete $L^2$ regression to approximate~\eqref{eq:surrogateDriftAndDiff}. We provide the details of the approximation method in Appendix~\ref{appx:L2regression}. 

For simplicity, we assume that the coordinate system is oriented so that one component is orthogonal to the hyperplane that defines the rare event threshold, e.g., $P_1u_t=u_t^i$ for a given $i=1,...,d$. For a general form of the projection $P_1$, we can use an orthogonal coordinate transformation such that the rare event is defined by a single component in the transformed coordinates with the condition that the first row of the rotation matrix and $P_1^T$ are collinear

Following the same idea as in one dimension, we solve the KBE corresponding to the surrogate process $\check{S}_t$ defined by
\begin{equation}\label{eq:BKEpdeSurrogate}
	\left\{
	\begin{array}{ll}
		\frac{\partial \check{\gamma}}{\partial t}  = -\check{a}(y, t)	\frac{\partial \check{\gamma}}{\partial y}   -\frac{1}{2} \check{b}(y,t) \check{b}(y, t)^T 	\frac{\partial^2 \check{\gamma}}{\partial y},\\
		\check{\gamma}(y, t_{n+1}) =0, \quad \;\; y<\cK,\\
	\check{\gamma}(\cK, t) =1, \qquad \quad t\in [t_n, t_{n+1}],\\
	\lim\limits_{s\rightarrow -\infty} \check{\gamma}(y,t)=0, \;\, t\in [t_n, t_{n+1}].
	\end{array}
	\right.
\end{equation}
and approximate the high-dimensional optimal IS initial density for $x_0 \in \bR^d$ by 
\begin{equation}\label{eq:optimindenhighdim}
	\tilde{\rho}_{0}^{PDE,1}(x_0)\propto \rho_{0}(x_0) \sqrt{\gamma(x_0,t_n)}  \approx \rho_{0}(x_0) \sqrt{\check{\gamma}(P_1x_0,t_n)}.
\end{equation}
Note that $\rho_0(x_0)$ is here a multivariate Gaussian density. For simplicity, let us provide an example of how we fit a Gaussian density to the product $\rho_{0}(x_0) \sqrt{\check{\gamma}(P_1x_0,t_n)}$ in the two-dimensional case.

\begin{example} 
Let  $x_0 = (x_{0}^{1}, x_{0}^{2}) \in \bR^2$ and the projection $P_1x_0=x_0^1$. Let $\rho_0(x_0)$ be a bivariate Gaussian density with mean 
$\begin{bmatrix}
	\mu_{0}^1\\
	\mu_{0}^2
\end{bmatrix}$ and covariance 
$\begin{bmatrix}
	(\sigma_{0}^1)^2 & \varrho\sigma_{0}^1\sigma_{0}^2 \\
	\varrho\sigma_{0}^1\sigma_{0}^2 & (\sigma_{0}^2)^2
\end{bmatrix}$ where $\varrho$ is a correlation coefficient between $x_0^1$ and $x_0^2$. Then, the optimal IS initial density $\tilde{\rho}_{0}(x_{0}^{1}, x_{0}^{2})\propto \rho_0(x_0^2|x_0^1)\rho_0(x_0^1)\sqrt{\bar{\gamma}(x_{0}^{1}, x_{0}^{2};t_n)}$ can be approximated by the probability density
 \begin{equation*}
 	\begin{split}
 	\check{\tilde{\rho}}_{0}(x_{0}^{1}, x_{0}^{2}) \propto \rho_0(x_0^2|x_0^1)\rho_0(x_0^1)\sqrt{\check{\gamma}(x_0^1,t_n)}.
 \end{split}
 \end{equation*}
We can fit the univariate Gaussian density to the product $\rho_0(x_0^1)\sqrt{\check{\gamma}(x_0^1,t_n)}$ with parameters $\mu_{0}^{\mathrm{fit}}$ and $\sigma_{0}^{\mathrm{fit}}$. Note that $\check{\tilde{\rho}}_0(x_0^2|x_0^1) = \rho_0(x_0^2|x_0^1)$ and
\[
\check{\tilde{\rho}}_0(x_0^1) \propto \int \rho_0(x_0^2|x_0^1)\rho_0(x_0^1)\sqrt{\check{\gamma}(x_0^1,t_n)} dx_0^2 = \rho_0(x_0^1)\sqrt{\check{\gamma}(x_0^1,t_n)}. 
\]
Applying a definition of the conditional bivariate Gaussian distribution, we have
\[
\check{\tilde{\rho}}_0(x_0^2|x_0^1) \sim N\Big(\mu_{0}^2+\varrho\frac{\sigma_{0}^2}{\sigma_{0}^1}(x_0^1-\mu_{0}^1), (1-\varrho^2)(\sigma_{0}^2)^2\Big),
\]
where $x_0^1\sim N(\mu_{0}^{\mathrm{fit}}, \sigma_{0}^{\mathrm{fit}})$. 
The likelihood ratio, in this case, is defined by
\[
L_0(x_0)= \frac{\rho_0(x_0^1)}{\check{\tilde{\rho}}_0(x_0^1)}=\frac{\sigma_{0}^{\mathrm{fit}}}{\sigma_{0}^1}\exp{\prt{-\frac{(x_0^1-\mu_{0}^1)^2}{2(\sigma_{0}^1)^2}+\frac{(x_0^1-\mu_{0}^{\mathrm{fit}})^2}{2(\sigma_{0}^{\mathrm{fit}})^2}}}
\]
A similar approach is used in the numerical section~\ref{ssec:Langevin}.
\end{example}

\begin{remark}[Alternative to Gaussian fitting]
	Instead of Gaussian fitting, another alternative approach, which is not studied in this work, is to consider non-parametric methods such as the kernel density estimations (KDE) based on the ensemble of particles. In this case, one should explore how to effectively sample from the KDE and compute the corresponding likelihood. 
\end{remark}

\subsubsection{\textbf{Multilevel Cross-Entropy method}}\label{ssec:CEMethod} For the purpose of comparison, we consider a simple, efficient iterative procedure known as the multilevel cross-entropy (CE) method to estimate the optimal IS density~\eqref{eq:optimdensityPropto}. The core idea of the CE is to assume $\tilde{\rho}_{0}(x_{0})$ to be from a particular family of densities $\{\rho_{0} (x_0; \bv) \}$ (in the framework of EnKF - Gaussian family) and select the parameter vector $\bv$ such that the distance between $\rho_{0}(x_{0}; \bv)$ and $\tilde{\rho}_{0} (x_0)$ in~\eqref{eq:optimdensityPropto} is minimal~\cite{de2005tutorial, kroese2013cross}. A methodology and basic algorithm of the multilevel CE adjusted to the EnKF setting are briefly provided in Appendix~\ref{appx:CEMethod}. We denote the optimal IS density obtained via the multilevel CE method by  $\tilde{\rho}_{0}^{CE}$. 
\subsection{IS with respect to Wiener processes using SOC}
\label{ssec:ISdW}

Here, we base the analysis on the continuous time problem~\eqref{eq:gensde}. We apply IS wrt Wiener processes obtained by stochastic optimal control for the  SDEs of the form~\eqref{eq:gensde} with a fixed initial condition $x_0$. 

\textit{One-dimensional setup.} Let us introduce a process $$W_t^{\mathbb{Q}}=W_t^{\mathbb{P}}-\int_{t_n}^{t} \xi(u_s, s)ds$$ 
where $\xi: \bR \times [t_n, t_{n+1}] \rightarrow \bR^{d_W\times 1}$ and substitute it in~\eqref{eq:gensde} for $t\in (t_n,t_{n+1})$ to get
\begin{equation}
	\begin{cases}
		du_t=  (a(u_t)+b(u_t)\xi(u_t, t))dt+b(u_t)dW_t^{\mathbb{Q}}, \\
		u_n =x_0.
	\end{cases}\,
\end{equation}

Girsanov's theorem~\cite{karatzas2012brownian} states under certain conditions that $(W^\mathbb{Q}(t))_{t\geq0}$  is a Wiener process under the probabilty measure $\mathbb{Q}$ with the likelihood ratio
\begin{equation}\label{eq:lklhW}
\frac{d\mathbb{P}}{d\mathbb{Q}}|_{[t_n, t]}=\exp\prt{{-\int_{t_n}^{t} \xi(u_s, s)^TdW_s^{\mathbb{Q}}-\frac{1}{2}\int_{t_n}^{t}\xi(u_s, s)^T\xi(u_s,s)ds}}.
\end{equation}
Then, a change of measure wrt the Wiener process is as follows
\begin{equation*}
\E_{\rho_0 \otimes \mathbb{P}}[\Ind_{u_{\tau}\geq\cK}]=\E_{\rho_0 \otimes \mathbb{Q}}\Big[\Ind_{u_{\tau}^{\xi}\geq\cK}\frac{d\mathbb{P}}{d \mathbb{Q}}|_{[t_n, \tau]}\Big],
\end{equation*}
where $u_{t}^{\xi}$ follows the controlled SDE
\begin{equation*}\label{eq:controlledsde1d}
	\begin{cases}
		du_t^{\xi}=  (a(u_t^{\xi})+b(u_t^{\xi})\xi(u_t^{\xi}, t))dt+b(u_t^{\xi})dW_t^{\mathbb{Q}}, \\
		u_n^{\xi}=x_0.
	\end{cases}\,
\end{equation*}
Now, the aim is to find a control $\xi$ that minimizes the variance of the QoI. Note that it is sufficient to minimize the second moment of the estimator that is
\begin{equation*}
	\begin{split}
\min_{\xi} \mathbb{E}_{\rho_0 \otimes \mathbb{Q}} \Big[\Ind_{u_{\tau}^{\xi}\geq\cK}  e^{-2\int_{t_n}^{\tau} \xi(u_s^{\xi}, s)^TdW_s^{\mathbb{Q}}-\int_{t_n}^{\tau}\xi(u_s^\xi, s)^T\xi(u_s^\xi,s)ds} \Big | u_n^{\xi}=x_0\Big].
\end{split}
\end{equation*}
Let us define the value function that minimizes the second moment of the MC estimator:
\begin{equation*}
	\begin{split}
	\theta(x,t)=
		\min_{\xi} \mathbb{E}_{\rho_0 \times \mathbb{Q}} \Big[\Ind_{u_{\tau}^{\xi}\geq\cK}  e^{-2\int_{t}^{\tau} \xi(u_s^{\xi}, s)^TdW_s^{\mathbb{Q}}-\int_{t}^{\tau}\xi(u_s^\xi, s)^T\xi(u_s^\xi,s)ds} \Big | u_t^{\xi}=x\Big].
	\end{split}
\end{equation*}
Then, $\theta(x,t)$ satisfies the following initial-boundary problem
\begin{equation}\label{eq:SOCpde}
	\left\{
	\begin{array}{ll}
		\frac{\partial \theta }{\partial t} = -a(x) \frac{\partial \theta }{\partial x}  -\frac{1}{2} b(x) b(x)^T \frac{\partial^2 \theta }{\partial x^2} +\frac{1}{2 \theta}b(x)b(x)^T 	\prt{\frac{\partial \theta }{\partial x}}^2,\\
		\theta(x, t_{n+1}) =0, \quad  \; \,  x<\cK,\\
		\theta(\cK, t) =1, \quad \quad \; \; \; t\in[t_n,t_{n+1}]\\
		\lim\limits_{x\rightarrow -\infty}	\theta(x, t)=0,\; t\in[t_n,t_{n+1}],
	\end{array}
	\right.
\end{equation}
with the optimal control
\[
\xi^{\star}(x,t)=\frac{1}{2}b(x)^T\frac{\partial \log \theta(x,t)}{\partial x}.
\]


Alternatively, we can consider $\theta(x,t)=\gamma^2(x,t)$ which leads to the KBE defined in~\eqref{eq:BKEpde} with the optimal control
\begin{equation}\label{eq:optcontrol1d}
\xi^*(x,t) = b(x)^{T} \frac{\partial \log \gamma (x,t)}{\partial x}.
\end{equation}

\textit{High-dimensional setup}. In more than one dimension, similarly to Section~\ref{ssec:PdeMethod}, we apply the Markovian projection technique and apply a change of measure wrt the Wiener process to the one-dimensional surrogate process $\check{S}_t$. Then, the value function satisfies the above-mentioned PDEs since the Markov property for $\check{S}_t$ holds. The high-dimensional optimal control with $x \in \bR^d$ can be approximated via a solution of the one-dimensional PDE by
\begin{equation}\label{eq:optcontrolmultid}
\xi^*(x,t) \approx b(x)^TP_1^T\frac{\partial \log \check{\gamma}  (P_1x,t)}{\partial y},
\end{equation}
where the partial derivative is taken wrt the component on which we project.
 
Correspondingly, the high-dimensional controlled SDE takes the form 
\begin{equation*}
	\begin{cases}
		du_t^{\xi}=  (a(u_t^{\xi})+b(u_t^{\xi})\xi(u_t^{\xi}, t))dt+b(u_t^{\xi})dW_t^{\mathbb{Q}}, \\
		u_n^{\xi}=x_0,
	\end{cases}\,
\end{equation*}
where $\xi: \bR^d \times [t_n, t_{n+1}] \rightarrow \bR^{d_W \times 1}$ and the likelihood ratio is 
\begin{equation}\label{eq:lklhWmultiD}
	\frac{d\mathbb{P}}{d\mathbb{Q}}|_{[t_n, t]}=\exp\prt{{-\int_{t_n}^{t} \xi(u_s, s)^T dW_s^{\mathbb{Q}}-\frac{1}{2}\int_{t_n}^{t}\xi(u_s, s)^T\xi(u_s,s)ds}}.
\end{equation}

\subsection{IS with respect to both initial condition and Wiener processes }
\label{ssec:ISboth}
Recall that $ \rho_{\bar{u}_{n,0:K}}$ denotes the joint density of a random vector $(\bar{u}_{n,0:K})$ with $\rho_{\bar{u}_{n,0}}=\rho_0$. Now we introduce an alternative joint density function $\tilde{\rho}_{\bar{u}_{n,0:K}^{\xi}}$ and consider
\begin{equation}\label{eq:QoIchangeboth}
	\begin{split}
		\E_{\rho_0 \otimes \bP}[\IndU{ \bar{u}_{n,0:K}}]&=\int  \IndU{x_{0:K}} \rho_{\bar{u}_{n,0:K}} (x_{0:K})dx_{0:K}\\
		&= \int \IndU{x_{0:K}} \frac{\rho_{\bar{u}_{n,0:K}} (x_{0:K})}{\tilde{\rho}_{\bar{u}_{n,0:K}^{\xi}} (x_{0:K})}  \tilde{\rho}_{\bar{u}_{n,0:K}^{\xi}} (x_{0:K})dx_{0:K}\\
		&=\E_{\tilde{\rho}_0 \otimes \mathbb{Q}} \Big[ \IndU{ \bar{u}_{n,0:K}^{\xi}} \frac{\rho_{\bar{u}_{n,0:K}} (\bar{u}_{n,0:K}^{\xi})}{\tilde{\rho}_{\bar{u}_{n,0:K}^{\xi}} (\bar{u}_{n,0:K}^{\xi})}\Big]\\
		&=\E_{\tilde{\rho}_0 \otimes \mathbb{Q}} \Big[ \IndU{ \bar{u}_{n,0:K}^{\xi}} L_0(\bar{u}_{n,0}^{\xi})L_W(\bar{u}_{n,1:K}^{\xi})\Big],
	\end{split}
\end{equation}
where $L_0(\cdot)$ is a likelihood ratio of the IS wrt the initial condition and $L_W(\cdot)$ is a likelihood ratio of the IS wrt the Wiener processes.
 
Using the Cauchy-Schwartz inequality, similarly to~\eqref{eq:Cauchy-Schwartz}, the optimal importance joint density that reduces the variance of the estimator is proportional to
\[
\tilde{\rho}_{\bar{u}_{n,0:K}^{\xi}} (x_{0:K}) \propto  \IndU{x_{0:K}}\rho_{\bar{u}_{n,0:K}}(x_{0:K}).
\]
Furthermore, we can deduce that
\begin{equation}\label{eq:optiminitdensityISboth}
	\begin{split}
		&\tilde{\rho}_{\bar{u}_{n,0}^{\xi}} (x_{0}) = \int \tilde{\rho}_{\bar{u}_{n,0:K}^{\xi}} (x_{0:K})dx_{1:K}  \propto \int   \IndU{x_{0:K}}\rho_{\bar{u}_{n,0:K}}(x_{0:K}) dx_{1:K}\\
		& = \int  \IndU{x_{0:K}} \rho_{0}(x_{0}) \rho_{ \bar{u}_{n,1:K}|\bar{u}_{n,0}} ( x_{1:K}|x_0)  dx_{1:K}\\
		&=\rho_{0}(x_{0})  \int \big( \Ind_{P_1x_0\geq \cK}(x_0)+\Ind_{P_1x_0<\cK}(x_0)\IndU{x_{1:K}}\big)\rho_{ \bar{u}_{n,1:K}|\bar{u}_{n,0}} ( x_{1:K}|x_0)  dx_{1:K}\\
		&= \rho_{0}(x_{0}) \begin{cases}
			 \bP(\max_{1\leq k\leq K } P_1 \bar{u}_{n,k} \geq \cK|\bar{u}_{n,0}=x_0), \quad P_1x_0<\cK,\\
			 	1, \quad P_1x_0\geq \cK.
		\end{cases}
	\end{split}
\end{equation}

\textit{One-dimensional setup.} Similarly to Section~\ref{ssec:ISrho0}, for $d=1$, the optimal importance initial density in the IS wrt both initial condition and Wiener process is given by
\[
	\tilde{\rho}_{0}^{PDE,2}(x_0)\propto \rho_{0}(x_0)\gamma(x_0,t_n), \quad x_0<\cK,
\]
where $\gamma(x_0,t_n)$ is an approximation for the terminal time solution of the KBE~\eqref{eq:BKEpde}. Notice that~\eqref{eq:optimdensityPropto} differs from~\eqref{eq:optiminitdensityISboth} by having a square root in the exit probability. 

The optimal control is defined similar to~\eqref{eq:optcontrol1d} and the likelihood ratios are correspondingly given by
\begin{equation}\label{eq:discreteliklBoth1d}
	\begin{split}
	L_0(\bar{u}_{n,0}^{\xi})&=\frac{\rho_{0} (\bar{u}_{n,0}^{\xi})}{\tilde{\rho}_{0}^{PDE,2} (\bar{u}_{n,0}^{\xi})}=\frac{\sigma_{0}^{\mathrm{fit}}}{\sigma_{0}}\exp{\prt{-\frac{(x_0-\mu_{0})^2}{2\sigma_{0}^2}+\frac{\big(x_0-\mu_{0}^{\mathrm{fit}}\big)^2}{2\big(\sigma_{0}^{\mathrm{fit}}\big)^2}}},\\
	L_W(\bar{u}_{n,1:K}^{\xi})&=\prod_{k=1}^{K-1}\exp{\prt{-\xi(t_{n,k}, \bar{u}_{n,k}^\xi)^T\Delta W_{n,k}^{\mathbb{Q}}-\frac{1}{2}\Delta t_{n,k} \xi(t_{n,k}, \bar{u}_{n,k}^\xi)^T\xi(t_{n,k}, \bar{u}_{n,k}^\xi)}},
\end{split}
\end{equation}
where the control $\xi$ is the zero vector after a stopping time. Note that $L_W$ here is the discrete form of the likelihood defined in~\eqref{eq:lklhW}.

\textit{High-dimensional setup.} For $d>1$, similarly to the previous sections, as a result of the Markovian projection technique, the optimal control is given by~\eqref{eq:optcontrolmultid}, the optimal initial density is approximated by 
\[
\tilde{\rho}_{0}^{PDE,2}(x_0)\propto \rho_{0}(x_0)\gamma(x_0,t_n)\approx \rho_{0}(x_0) \check{\gamma}(P_1x_0,t_n), \quad P_1x_0<\cK,
\]
 and the likelihood ratios are respectively given by
\begin{equation}\label{eq:discreteliklbothmultid}
	\begin{split}
		L_0(\bar{u}_{n,0}^{\xi})&=\frac{\rho_{0} (\bar{u}_{n,0}^{\xi})}{\tilde{\rho}_{0}^{PDE,2} (\bar{u}_{n,0}^{\xi})},\\
		L_W(\bar{u}_{n,1:K}^{\xi})&=\prod_{k=1}^{K-1}\exp{\prt{-\xi(t_{n,k}, \bar{u}_{n,k}^\xi)^T \Delta W_{n,k}^{\mathbb{Q}}-\frac{1}{2}\xi(t_{n,k}, \bar{u}_{n,k}^\xi)^T\xi(t_{n,k}, \bar{u}_{n,k}^\xi)\Delta t_{n,k}}},
	\end{split}
\end{equation}
where the control $\xi$ is the zero vector after a stopping time. Note that $L_W$ here is the discrete form of the likelihood defined in~\eqref{eq:lklhWmultiD}.


\section{Numerical examples}
\label{sec:numerics}
In this section, we provide numerical examples that showcase the application of the 
IS techniques presented in Section~\ref{sec:IStechniques}. We test the proposed IS approaches in a range of problems: 
a Double Well SDE, Langevin dynamics and a noisy Charney-deVore model~\cite{hoel2020multilevel, hoel2022multi, grafke2019numerical}. Experimental support for the variance reduction of the proposed methods is provided in comparison with the standard MC and the multilevel CE methods. 

The standard MC estimator to approximate the probability $\hat{\alpha}_n$ in~\eqref{eq:QoIchangeinitial} of the discretized dynamics  is
\[
\hat{\alpha}_n^{MC}:=\frac{1}{J} \sum_{i=1}^{J} \IndU{\bar{u}_{n,0:K}^{[i]}}
\]
where $\Big\{\IndU{\bar{u}_{n,0:K}^{[i]}}\Big\}_{i=1}^J$ are i.i.d realizations of $\big\{\IndU{\bar{u}_{n,0:K}}\big\}$. The estimator is unbiased, which means $\E[\hat{\alpha}_n^{MC}]=\hat{\alpha}_n$, and the variance is
\[
\V[\hat{\alpha}_n^{MC}]=\frac{\hat{\alpha}_n-\hat{\alpha}_n^2}{J}.
\]    
The asymptotic result in the CLT motivates the confidence interval (CI) for $\hat{\alpha}_n$ as
\[
\hat{\alpha}_n^{MC} \pm z_c \sqrt{\V[\hat{\alpha}_n^{MC}]},
\]
where $z_c$ is a constant corresponding to a given level of confidence. We consider $95\%$ CI in our numerical simulations by setting $z_c=1.96$. 

In the regime of rare events, that is, $\hat{\alpha}_n$ tends to be very small, it is appropriate to measure the accuracy of the estimation as a relative statistical error, which we estimate by the ratio of the standard deviation and the mean
 $$\epsilon_{st}^{MC}=z_c \frac{\sqrt{\V[\hat{\alpha}_n^{MC}]}}{\hat{\alpha}_n^{MC} }.$$ 

The method of numerical approximation to~\eqref{eq:BKEpde} is desribed in Appendix~\ref{appx:PDEsolver}.  A numerical approximation of the corresponding optimal control is obtained by a central difference with linear interpolation throughout the domain (see~\cite[Appendix E]{ben2022double} for more details).

\subsection{Double-Well problem}\label{ssec:dw}
We consider nonlinear dynamics with a constant diffusion and a drift arising for a double well (DW) potential
\[
du_t=-V'(u_t)dt+b dW_t,
\]
where $V(u)=\frac{1}{2+4u^2}+\frac{u^2}{4}$. 
Figure~\ref{fig:DWtrajAndPotential} illustrates the metastable behaviour of the DW dynamics between two wells over $T=100$ observation times. Notice that the trajectory tends to remain in one well for a considerably longer time and transit to another well.

\begin{figure}[h!]
	\includegraphics[height=6cm, width=8cm]{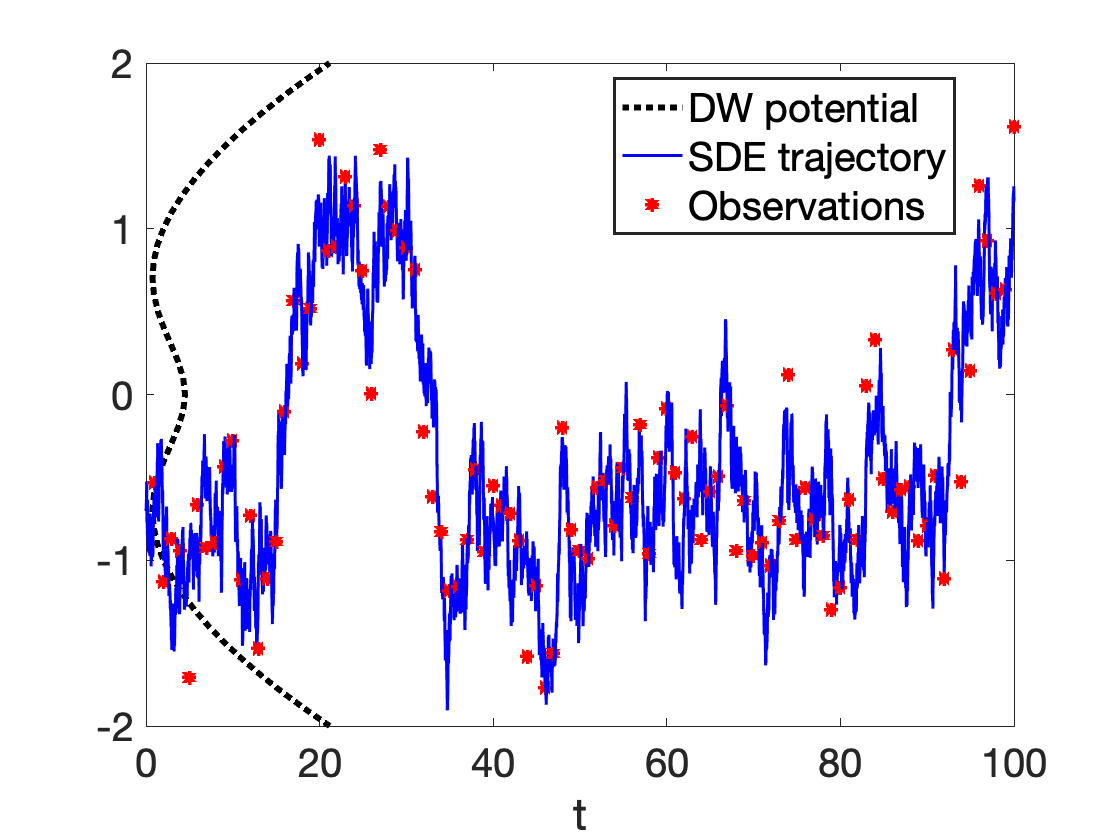}
	\caption{A trajectory of the DW SDE from Section~\ref{ssec:dw} over $T=100$ observation times with given observations and the double-well potential. The model parameters:  $b = 0.5$, $u_0 \sim N(-0.7,0.1)$. The EnKF parameters: $H=1$, $\Gamma=0.1$.}
	\label{fig:DWtrajAndPotential}
\end{figure}

\begin{table}[h!]
	\begin{adjustbox}{width=\columnwidth,center}
		\begin{tabular}{|c|c|c|c|c|l|}
			\hline
			\multirow{2}{*}{$\cK$} & \multicolumn{5}{c|}{\textbf{Estimator with 95\% CI} } \\ \cline{2-6}  
			&  \multicolumn{1}{c|}{\multirow{2}{*}{$\hat{\alpha}^{MC}$}} & \multicolumn{1}{c|}{\multirow{2}{*}{$\hat{\alpha}^{CE,\tilde{\rho}_{0}}$} }& \multicolumn{1}{c|}{\multirow{2}{*}{$\hat{\alpha}^{PDE,\tilde{\rho}_{0}}$}} & \multicolumn{1}{c|}{\multirow{2}{*}{$\hat{\alpha}^{PDE,W_t}$} }&\multicolumn{1}{c|}{\multirow{2}{*}{$\hat{\alpha}^{PDE,\text{both}}$} }\\  [2ex]\hline
			\multirow{2}{*}{$\mathbf{0}$}    & \multicolumn{1}{c|}{$\mathbf{4.65 \times 10^{-2}}$}& \multicolumn{1}{c|}{$\mathbf{4.66 \times 10^{-2}}$}     &  \multicolumn{1}{c|}{$\mathbf{4.65 \times 10^{-2}}$} & \multicolumn{1}{c|}{$\mathbf{4.63 \times 10^{-2}}$} & \multicolumn{1}{c|}{$\mathbf{4.63 \times 10^{-2}}$} \\
			& $[0.0461, 0.0469]$& $[0.0462, 0.0470]$&$[0.0461, 0.0469]$ &$[0.0462, 0.0464]$&$[0.04626, 0.04632]$\\  \hline
			
			\multirow{2}{*}{$\mathbf{0.5}$}    & \multicolumn{1}{c|}{$\mathbf{4.83 \times 10^{-3}}$}& \multicolumn{1}{c|}{$\mathbf{4.76 \times 10^{-3}}$}     &  \multicolumn{1}{c|}{$\mathbf{4.84 \times 10^{-3}}$} & \multicolumn{1}{c|}{$\mathbf{4.81 \times 10^{-3}}$} & \multicolumn{1}{c|}{$\mathbf{4.82 \times 10^{-3}}$} \\
			& $[0.0047, 0.0050]$& $[0.0046, 0.0049]$&$[0.0047, 0.0050]$ &$[0.00480, 0.00483]$&$[0.004815, 0.004822]$\\  \hline
			\multirow{2}{*}{$\mathbf{1}$}    & \multicolumn{1}{c|}{$\mathbf{1.81 \times 10^{-4}}$}& \multicolumn{1}{c|}{-}     &  \multicolumn{1}{c|}{$\mathbf{1.90 \times 10^{-4}}$} & \multicolumn{1}{c|}{$\mathbf{1.89 \times 10^{-4}}$} & \multicolumn{1}{c|}{$\mathbf{1.89 \times 10^{-4}}$} \\
			& $[1.6\text{e-}04, 2.1\text{e-}04]$& &$[1.7\text{e-}04, 2.1\text{e-}04]$ &$[1.88\text{e-}04, 1.90\text{e-}04]$&$[1.891\text{e-}04, 1.895\text{e-}04]$\\  \hline
			\multirow{2}{*}{$\mathbf{1.2}$}    & \multicolumn{1}{c|}{$\mathbf{3.70 \times 10^{-5}}$}& \multicolumn{1}{c|}{-}     &  \multicolumn{1}{c|}{$\mathbf{3.18 \times 10^{-5}}$} & \multicolumn{1}{c|}{$\mathbf{3.12 \times 10^{-5}}$} & \multicolumn{1}{c|}{$\mathbf{3.13 \times 10^{-5}}$}\\
			& $[2.5\text{e-}05, 4.9\text{e-}05]$& &$[2.4\text{e-}05, 4.0\text{e-}05]$ &$[3.1\text{e-}05, 3.2\text{e-}05]$&$[3.125\text{e-}05, 3.132\text{e-}05]$\\  \hline
		\end{tabular}
	\end{adjustbox}
\newline
\vspace*{0.2 cm}
	\begin{adjustbox}{width=\columnwidth,center}
		\begin{tabular}{|c|ccccl|l|l|l|l|}
			\hline
			\multirow{2}{*}{$\cK$}     & \multicolumn{5}{c|}{\textbf{Relative statistical error}}   &\multicolumn{4}{c|}{\textbf{Variance reduction}}\\\cline{2-10}  
			& \multicolumn{1}{c|}{\multirow{2}{*}{$\epsilon_{st}^{MC}$} }& \multicolumn{1}{c|}{\multirow{2}{*}{$\epsilon_{st}^{CE,\tilde{\rho}_{0}}$}} & \multicolumn{1}{c|}{\multirow{2}{*}{$\epsilon_{st}^{PDE,\tilde{\rho}_{0}}$} }&\multicolumn{1}{c|}{\multirow{2}{*}{$\epsilon_{st}^{PDE,W_t}$} }& \multicolumn{1}{c|}{\multirow{2}{*}{$\epsilon_{st}^{PDE,\text{both}}$} }& \multirow{2}{*}{$\frac{\mathbb{V}_{MC}}{\mathbb{V}_{CE,\tilde{\rho}_{0}}}$}& \multirow{2}{*}{$\frac{\mathbb{V}_{MC}}{\mathbb{V}_{PDE,\tilde{\rho}_{0}}}$} & \multirow{2}{*}{$\frac{\mathbb{V}_{MC}}{\mathbb{V}_{PDE,W_t}}$}&\multirow{2}{*}{$\frac{\mathbb{V}_{MC}}{\mathbb{V}_{PDE, \text{both}}}$}\\  [2ex]\hline
			$\mathbf{0}$          & \multicolumn{1}{c|}{$\mathbf{0.9\%}$}     &  \multicolumn{1}{c|}{$\mathbf{0.8\%}$} & \multicolumn{1}{c|}{$\mathbf{0.9\%}$} & \multicolumn{1}{c|}{$\mathbf{0.2\%}$}&\multicolumn{1}{c|}{$\mathbf{0.06\%}$} &\multicolumn{1}{c|}{$\mathbf{1.00} $} & \multicolumn{1}{c|}{$\mathbf{1.14} $} &\multicolumn{1}{c|}{$\mathbf{28} $} &\multicolumn{1}{c|}{$\mathbf{208} $} \\ \hline
			$\mathbf{0.5}$          & \multicolumn{1}{c|}{$\mathbf{2.8\%}$}     &  \multicolumn{1}{c|}{$\mathbf{3.0\%}$} & \multicolumn{1}{c|}{$\mathbf{2.5\%}$} & \multicolumn{1}{c|}{$\mathbf{0.3\%}$}&\multicolumn{1}{c|}{$\mathbf{0.07\%}$} &\multicolumn{1}{c|}{$\mathbf{0.90} $} & \multicolumn{1}{c|}{$\mathbf{1.29} $} &\multicolumn{1}{c|}{$\mathbf{96} $} &\multicolumn{1}{c|}{$\mathbf{1534} $} \\ \hline
			$\mathbf{1}$          & \multicolumn{1}{c|}{$\mathbf{14.6\%}$}     &  \multicolumn{1}{c|}{-} & \multicolumn{1}{c|}{$\mathbf{10.4\%}$} & \multicolumn{1}{c|}{$\mathbf{0.5\%}$}&\multicolumn{1}{c|}{$\mathbf{0.09\%}$} &\multicolumn{1}{c|}{-} & \multicolumn{1}{c|}{$\mathbf{1.78} $} &\multicolumn{1}{c|}{$\mathbf{691} $} &\multicolumn{1}{c|}{$\mathbf{20992} $} \\ \hline
			$\mathbf{1.2}$          & \multicolumn{1}{c|}{$\mathbf{32.2\%}$}     &  \multicolumn{1}{c|}{-} & \multicolumn{1}{c|}{$\mathbf{26.2\%}$} & \multicolumn{1}{c|}{$\mathbf{0.7\%}$}&\multicolumn{1}{c|}{$\mathbf{0.11\%}$} &\multicolumn{1}{c|}{-} & \multicolumn{1}{c|}{$\mathbf{2.06} $} &\multicolumn{1}{c|}{$\mathbf{2975} $} &\multicolumn{1}{c|}{$\mathbf{124939} $}  \\ \hline
		\end{tabular}
	\end{adjustbox}
	\caption{\textbf{ Double Well example.} Model parameter: $b=0.5$. Simulation parameters: $T=1$, $\Delta t= 0.01$,  $J = 10^6$, $u_0 \sim N(\mu_0,\sigma_0)$ with $\mu_0=-1$, $\sigma_0=0.2$. 
		Numerical results for the last two lines of the CE method are missing due to the impracticability of the algorithm for the given parameter setting and a limit of computer capacity to accomplish the simulation.}
	\label{table:DWsmallsigma}
\end{table}

\begin{figure}[h!]
	\includegraphics[height=6.5cm, width=6.5cm]{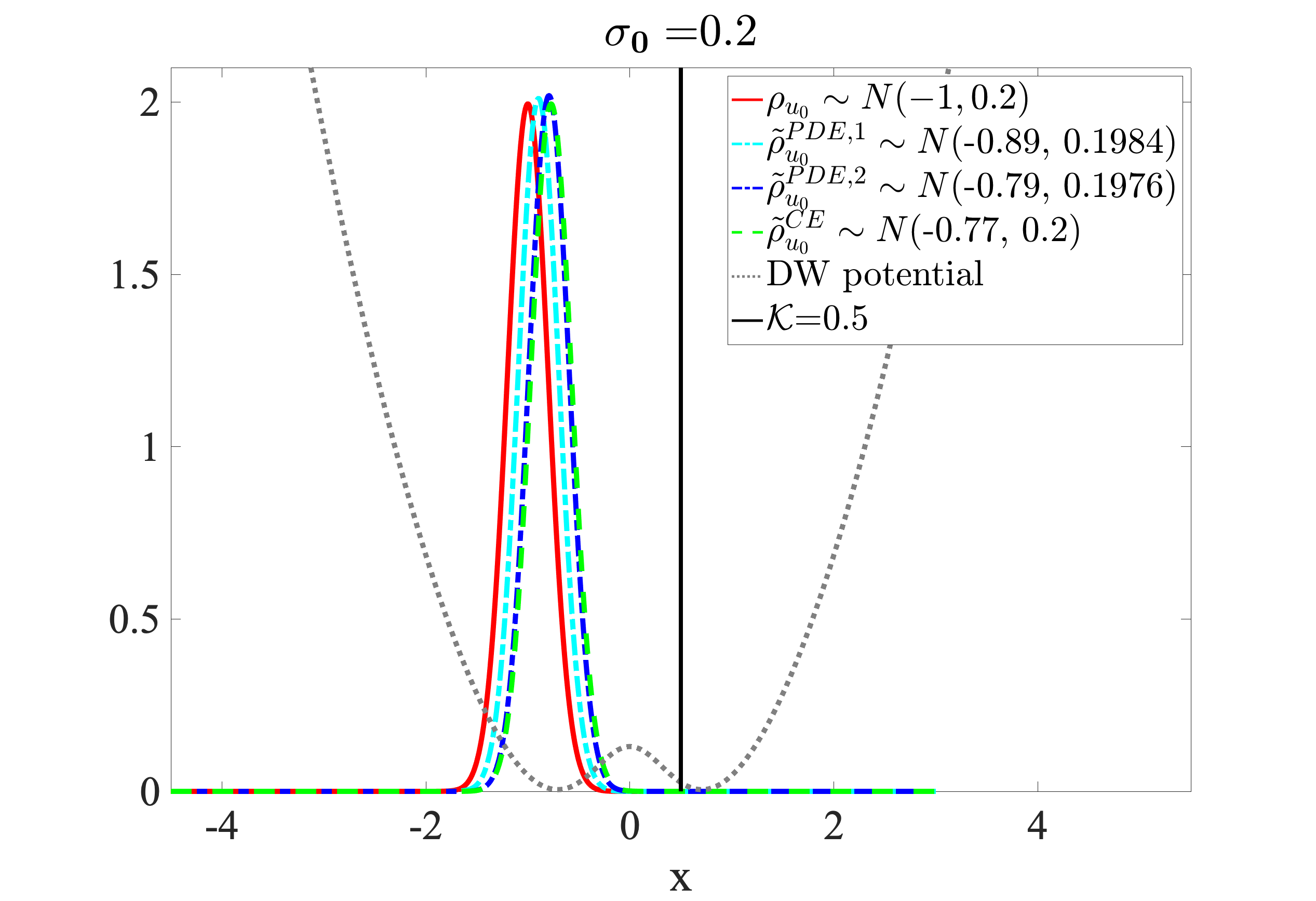}
	\hspace*{-0.63cm}
    \includegraphics[height=6.5cm, width=6.5cm]{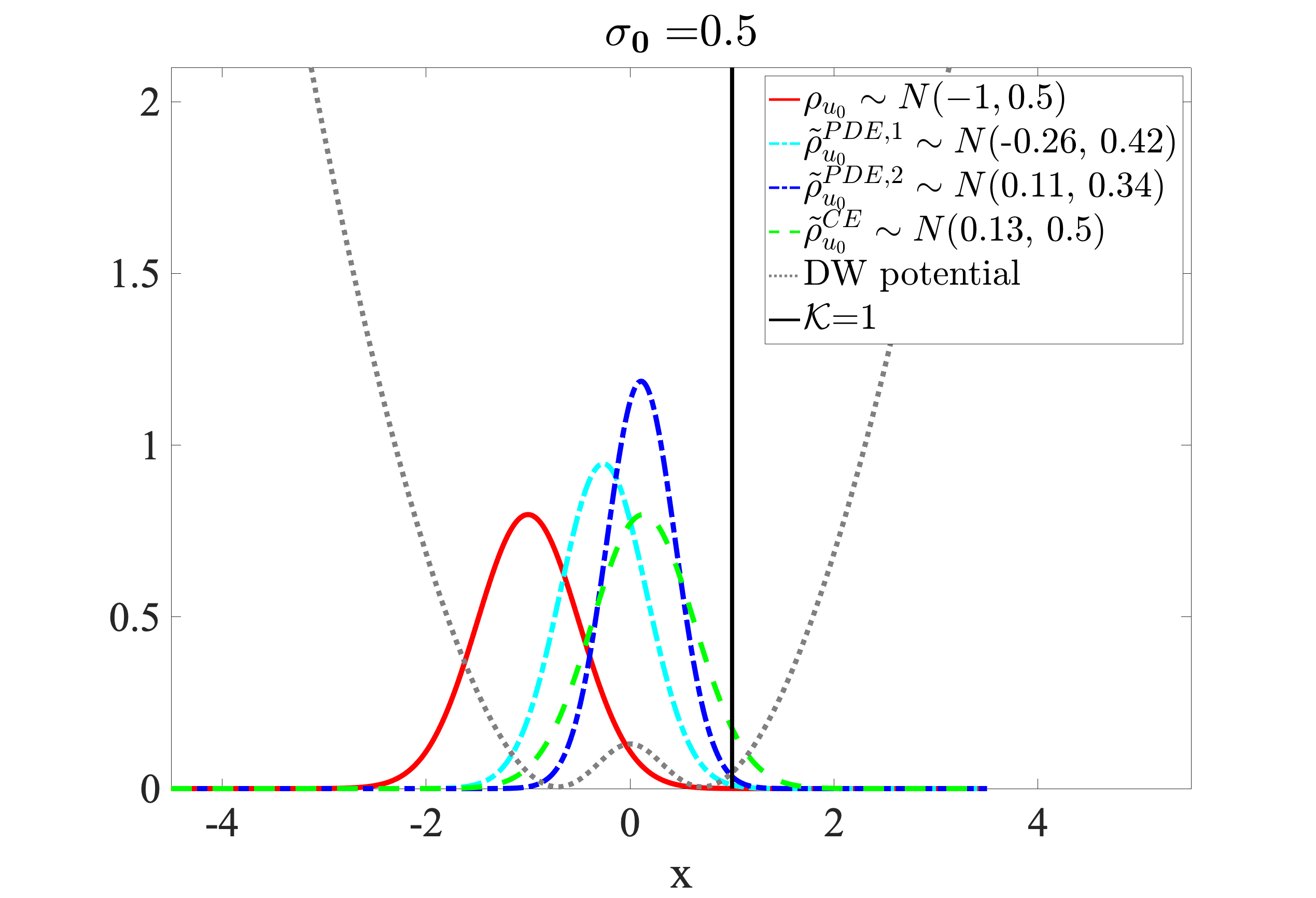}
    	\hspace*{-0.63cm}
	\includegraphics[height=6.5cm, width=6.5cm]{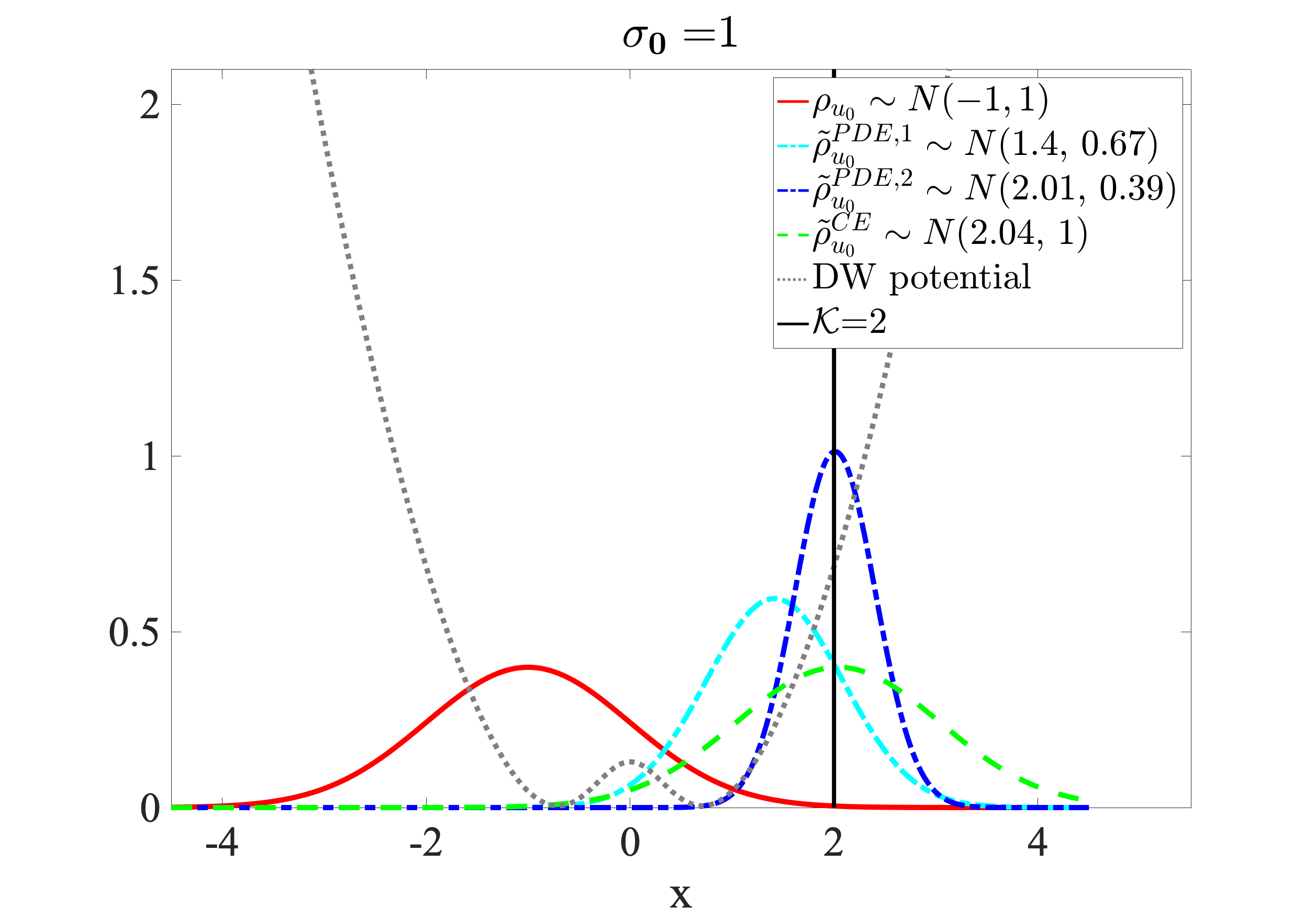}
	\caption{\textbf{Double Well example.} The original initial density ($\rho_{u_0}\sim N(-1, \sigma_0)$ in red solid line) is centered in one well with the threshold $\cK$ in the other well. For all three values of $\sigma_0$, $\cK$ is such that the rare event probability is around $10^{-3}$. A comparision of the approximated optimal initial densities based on: IS wrt $\rho_0$ ($\tilde\rho_{u_0}^{PDE,1}$ in cyan dash-dotted line), IS wrt both $\rho_0$ and $W(t)$ ($\tilde\rho_{u_0}^{PDE,2}$ in blue dash-dotted line), and CE-based IS ($\tilde\rho_{u_0}^{CE}$ in green dashed line). }
	\label{fig:DWchangeofrho0}
\end{figure}

 In Figure~\ref{fig:DWchangeofrho0}, we illustrate a comparison of the approximations to the optimal IS initial density $\tilde{\rho}_{u_0}$ obtained via the PDE-based and multilevel CE-based methods when the original initial density $\rho_{u_0}$ is given in one well and the threshold $\cK$ is set in the other well. Three subplots demonstrate how the IS densities change from the original one in the cases of different initial standard deviation $\sigma_0$. With a relatively small $\sigma_0$, all IS densities shift towards the threshold slightly, whereas with larger $\sigma_0$ values, they shift significantly to the rare event region. In particular, the CE-based IS density exhibits a considerable shift towards the threshold, however, it has the same scale as the original density by design. The PDE-based IS density $\tilde{\rho}_{u_0}^{PDE,2}$ wrt both initial condition and Wiener process has a larger shift and a higher concentration in the rare event region compared to the initial density $\tilde{\rho}_{u_0}^{PDE,1}$ obtained by IS wrt the initial condition. This is because of the difference in the optimal IS initial density formulas of each approach presented in~\eqref{eq:optimdensityPropto} and~\eqref{eq:optiminitdensityISboth}.
 The IS approach wrt the initial condition gives a substantial variance reduction only for comparatively large standard deviation values (see Table~\ref{table:DWlargesigma} for the example $\sigma_0=1$). In the regime of small $\sigma_0$, the IS approach wrt Wiener processes via the SOC theory proposed in Section~\ref{ssec:ISdW} works well (see Table~\ref{table:DWsmallsigma} for the example $\sigma_0=0.2$. The combination of both approaches (i.e., IS wrt both initial condition and Wiener processes) provides considerably better variance reduction over a crude MC method for all three cases of $\sigma_0$ (see Table~\ref{table:DWsmallsigma} and~\ref{table:DWlargesigma}). Figure~\ref{fig:DWci_Sig0_05} illustrates the consistency of the rare event probability estimates produced by different methods for the case of $\sigma_0=0.5$, where the $95\%$ confidence intervals (CI) overlap and converge as the MC sample size increases, and even with relatively small sample sizes, the IS wrt both $\rho_0$ and $W_t$ demonstrates superior improvement in reducing variance. A robustness of the numerical results is verified via bootstrapping. In particular, Figure~\ref{fig:DWrobustnessPDE} shows $95\%$ CI of two independent runs, which are overlapping, for the sample standard deviation (std) produced by the CE- and the PDE-based methods via bootstrap sampling with a bootstrap number $B=10^4$ over $J=10^6$ i.i.d. samples for the parameter setting of Table~\ref{table:DWlargesigma}. 

\begin{table}[h!]
	\begin{adjustbox}{width=\columnwidth,center}
		\begin{tabular}{|c|c|c|c|c|l|}
			\hline
			\multirow{2}{*}{$\cK$} & \multicolumn{5}{c|}{\textbf{Estimator with 95\% CI} } \\ \cline{2-6}  
			&  \multicolumn{1}{c|}{\multirow{2}{*}{$\hat{\alpha}^{MC}$}} & \multicolumn{1}{c|}{\multirow{2}{*}{$\hat{\alpha}^{CE,\tilde{\rho}_{0}}$} }& \multicolumn{1}{c|}{\multirow{2}{*}{$\hat{\alpha}^{PDE,\tilde{\rho}_{0}}$}} & \multicolumn{1}{c|}{\multirow{2}{*}{$\hat{\alpha}^{PDE,W_t}$} }&\multicolumn{1}{c|}{\multirow{2}{*}{$\hat{\alpha}^{PDE,\text{both}}$} }\\  [2ex]\hline
			\multirow{2}{*}{$\mathbf{1.5}$}    & \multicolumn{1}{c|}{$\mathbf{1.48 \times 10^{-2}}$}& \multicolumn{1}{c|}{$\mathbf{1.47 \times 10^{-2}}$}     &  \multicolumn{1}{c|}{$\mathbf{1.48 \times 10^{-2}}$} & \multicolumn{1}{c|}{$\mathbf{1.47 \times 10^{-2}}$} & \multicolumn{1}{c|}{$\mathbf{1.48 \times 10^{-2}}$} \\
			& $[0.0145, 0.0150]$& $[0.0147, 0.0148]$&$[0.0147, 0.0149]$ &$[0.0145, 0.0149]$&$[0.01479, 0.01480]$\\  \hline
			
			\multirow{2}{*}{$\mathbf{2}$}    & \multicolumn{1}{c|}{$\mathbf{2.5\times 10^{-3}}$}& \multicolumn{1}{c|}{$\mathbf{2.58 \times 10^{-3}}$}     &  \multicolumn{1}{c|}{$\mathbf{2.59 \times 10^{-3}}$} & \multicolumn{1}{c|}{$\mathbf{2.53 \times 10^{-3}}$} & \multicolumn{1}{c|}{$\mathbf{2.58 \times 10^{-3}}$} \\
			& $[0.0024, 0.0026]$& $[0.00256, 0.00260]$&$[0.00257, 0.00260]$ &$[0.0025, 0.0026]$&$[0.002578, 0.002584]$\\  \hline
			\multirow{2}{*}{$\mathbf{2.5}$}    & \multicolumn{1}{c|}{$\mathbf{3.8 \times 10^{-4}}$}& \multicolumn{1}{c|}{$\mathbf{3.98 \times 10^{-4}}$}     &  \multicolumn{1}{c|}{$\mathbf{3.98 \times 10^{-4}}$} & \multicolumn{1}{c|}{$\mathbf{3.81 \times 10^{-4}}$} & \multicolumn{1}{c|}{$\mathbf{3.97 \times 10^{-4}}$} \\
			& $[3.4\text{e-}04, 4.2\text{e-}04]$& $[3.96\text{e-}04, 4.00\text{e-}04]$&$[3.96\text{e-}04, 3.99\text{e-}04]$ &$[3.5\text{e-}04, 4.1\text{e-}04]$&$[3.96\text{e-}04, 3.97\text{e-}04]$\\  \hline
			\multirow{2}{*}{$\mathbf{3}$}    & \multicolumn{1}{c|}{$\mathbf{5.8 \times 10^{-5}}$}& \multicolumn{1}{c|}{$\mathbf{5.06 \times 10^{-5}}$}     &  \multicolumn{1}{c|}{$\mathbf{5.06 \times 10^{-5}}$} & \multicolumn{1}{c|}{$\mathbf{5.98 \times 10^{-5}}$} & \multicolumn{1}{c|}{$\mathbf{5.07 \times 10^{-5}}$}\\
			& $[4.3\text{e-}05, 7.3\text{e-}05]$& $[5.03\text{e-}05, 5.09\text{e-}05]$&$[5.04\text{e-}05, 5.09\text{e-}05]$ &$[4.6\text{e-}05, 7.3\text{e-}05]$&$[5.06\text{e-}05, 5.07\text{e-}05]$\\  \hline
		\end{tabular}
	\end{adjustbox}
\newline
\vspace*{0.2 cm}
	\begin{adjustbox}{width=\columnwidth,center}
		\begin{tabular}{|c|ccccl|l|l|l|l|}
			\hline
			\multirow{2}{*}{$\cK$}     & \multicolumn{5}{c|}{\textbf{Relative statistical error}}   &\multicolumn{4}{c|}{\textbf{Variance reduction}}\\\cline{2-10}  
			& \multicolumn{1}{c|}{\multirow{2}{*}{$\epsilon_{st}^{MC}$} }& \multicolumn{1}{c|}{\multirow{2}{*}{$\epsilon_{st}^{CE,\tilde{\rho}_{0}}$}} & \multicolumn{1}{c|}{\multirow{2}{*}{$\epsilon_{st}^{PDE,\tilde{\rho}_{0}}$} }&\multicolumn{1}{c|}{\multirow{2}{*}{$\epsilon_{st}^{PDE,W_t}$} }& \multicolumn{1}{c|}{\multirow{2}{*}{$\epsilon_{st}^{PDE,\text{both}}$} }& \multirow{2}{*}{$\frac{\mathbb{V}_{MC}}{\mathbb{V}_{CE,\tilde{\rho}_{0}}}$}& \multirow{2}{*}{$\frac{\mathbb{V}_{MC}}{\mathbb{V}_{PDE,\tilde{\rho}_{0}}}$} & \multirow{2}{*}{$\frac{\mathbb{V}_{MC}}{\mathbb{V}_{PDE,W_t}}$}&\multirow{2}{*}{$\frac{\mathbb{V}_{MC}}{\mathbb{V}_{PDE, \text{both}}}$}\\  [2ex]\hline
			$\mathbf{1.5}$          & \multicolumn{1}{c|}{$\mathbf{1.6\%}$}     &  \multicolumn{1}{c|}{$\mathbf{0.5\%}$} & \multicolumn{1}{c|}{$\mathbf{0.4\%}$} & \multicolumn{1}{c|}{$\mathbf{1.3\%}$}&\multicolumn{1}{c|}{$\mathbf{0.06\%}$} &\multicolumn{1}{c|}{$\mathbf{9} $} & \multicolumn{1}{c|}{$\mathbf{15} $} &\multicolumn{1}{c|}{$\mathbf{1.6} $} &\multicolumn{1}{c|}{$\mathbf{677} $} \\ \hline
			$\mathbf{2}$          & \multicolumn{1}{c|}{$\mathbf{3.9\%}$}     &  \multicolumn{1}{c|}{$\mathbf{0.6\%}$} & \multicolumn{1}{c|}{$\mathbf{0.5\%}$} & \multicolumn{1}{c|}{$\mathbf{3.2\%}$}&\multicolumn{1}{c|}{$\mathbf{0.10\%}$} &\multicolumn{1}{c|}{$\mathbf{40} $} & \multicolumn{1}{c|}{$\mathbf{61} $} &\multicolumn{1}{c|}{$\mathbf{1.4} $} &\multicolumn{1}{c|}{$\mathbf{1327} $} \\ \hline
			$\mathbf{2.5}$          & \multicolumn{1}{c|}{$\mathbf{10.1\%}$}     &  \multicolumn{1}{c|}{$\mathbf{0.6\%}$} & \multicolumn{1}{c|}{$\mathbf{0.5\%}$} & \multicolumn{1}{c|}{$\mathbf{8.7\%}$}&\multicolumn{1}{c|}{$\mathbf{0.10\%}$} &\multicolumn{1}{c|}{$\mathbf{261} $} & \multicolumn{1}{c|}{$\mathbf{360} $} &\multicolumn{1}{c|}{$\mathbf{1.3} $} &\multicolumn{1}{c|}{$\mathbf{6150} $} \\ \hline
			$\mathbf{3}$          & \multicolumn{1}{c|}{$\mathbf{25.7\%}$}     &  \multicolumn{1}{c|}{$\mathbf{0.6\%}$} & \multicolumn{1}{c|}{$\mathbf{0.4\%}$} & \multicolumn{1}{c|}{$\mathbf{22.4\%}$}&\multicolumn{1}{c|}{$\mathbf{0.10\%}$} &\multicolumn{1}{c|}{$\mathbf{2641} $} & \multicolumn{1}{c|}{$\mathbf{4580} $} &\multicolumn{1}{c|}{$\mathbf{1.2} $} &\multicolumn{1}{c|}{$\mathbf{46372} $}  \\ \hline
		\end{tabular}
	\end{adjustbox}
	\caption{\textbf{ Double Well example.} Model parameter: $b=0.5$. Simulation parameters: $T=1$, $\Delta t= 0.01$,  $J = 10^6$, $u_0 \sim N(\mu_0,\sigma_0)$ with $\mu_0=-1$, $\sigma_0=1$. 
	}
	\label{table:DWlargesigma}
\end{table}

\begin{figure}[h!]
	\includegraphics[height=6cm, width=8cm]{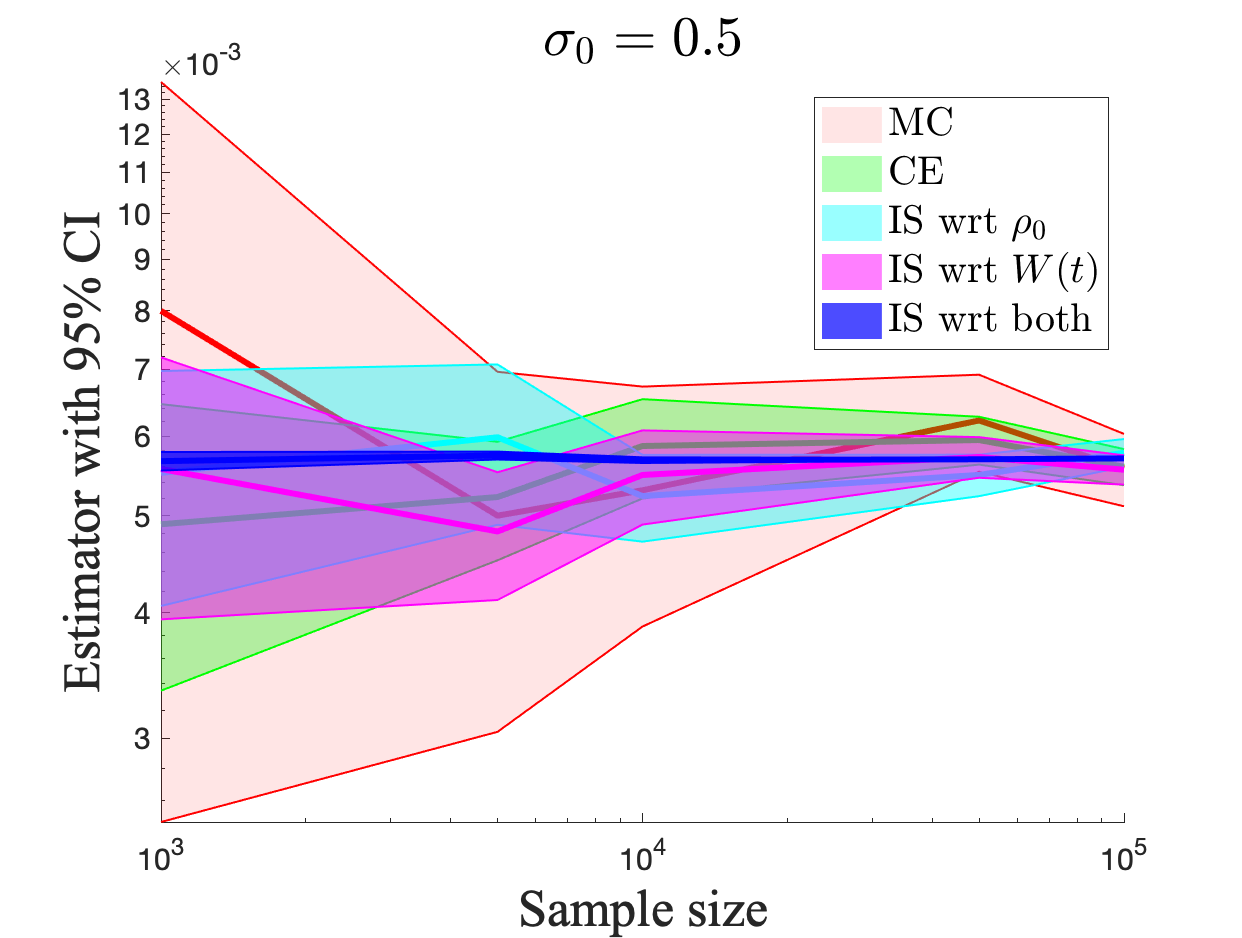}
	\caption{\textbf{Double Well example} Model parameter: $b=0.5$. 95\% CI of the estimator with the simulation parameters: $T=1$, $\Delta t= 0.01$,  $\cK=1$, $u_0 \sim N(\mu_0,\sigma_0)$ with $\mu_0=-1$, $\sigma_0=0.5$.}
	\label{fig:DWci_Sig0_05}
\end{figure}

\begin{figure}[h!]
	\centering
	\hspace*{-0.6cm}
	\includegraphics[height=6.3cm, width=6.3cm]{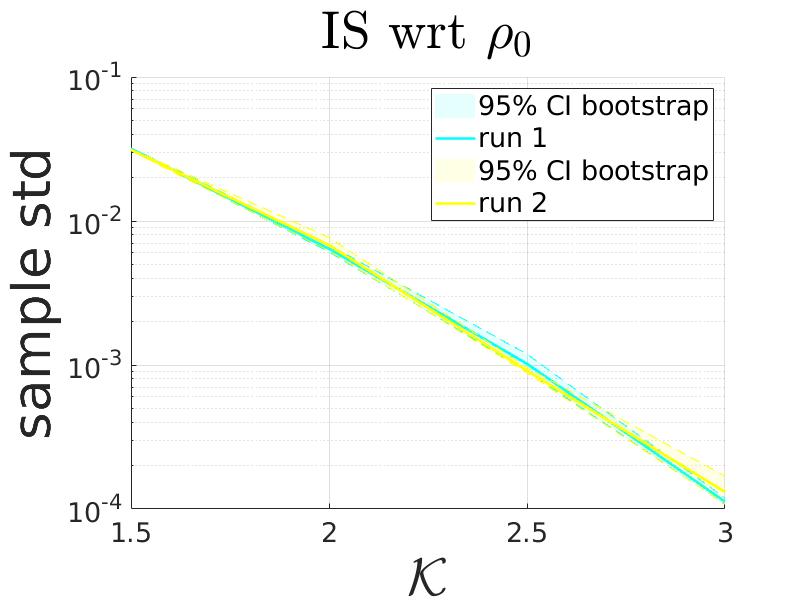}
	\hspace*{-0.1cm}
	\includegraphics[height=6.3cm, width=6.3cm]{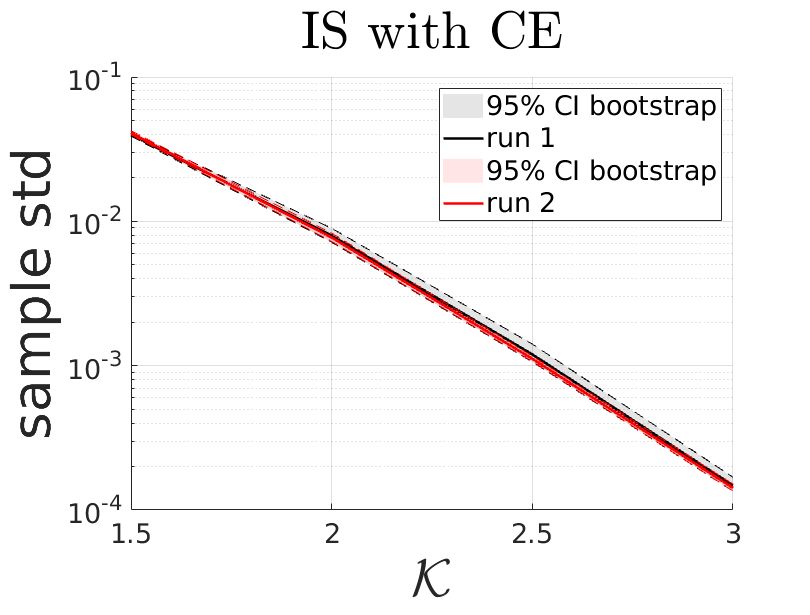}\\
    \hspace*{-0.6cm}
	\includegraphics[height=6.3cm, width=6.3cm]{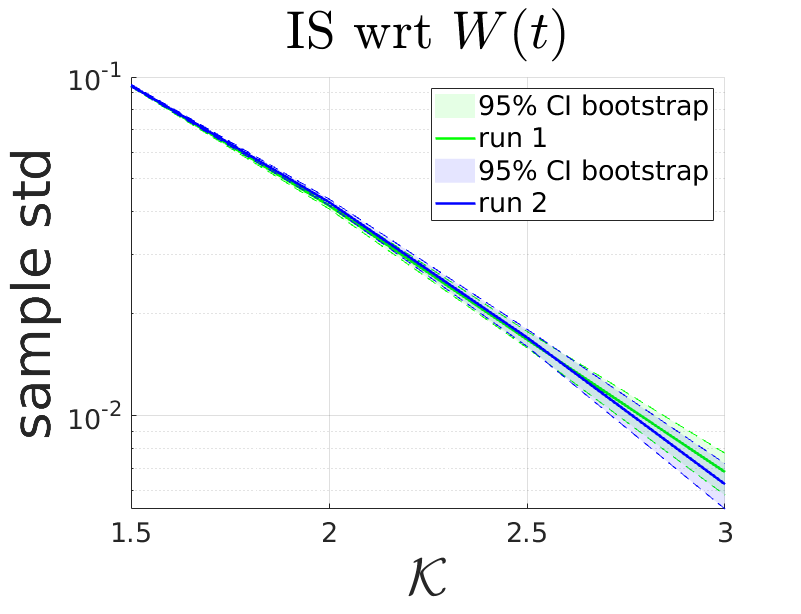}
	\hspace*{-0.1cm}
	\includegraphics[height=6.3cm, width=6.3cm]{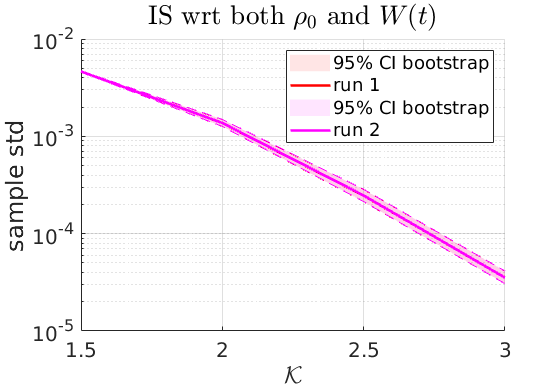}
	\caption{\textbf{Double Well example.} Two independent simulation results of the bootstrapping for a sample standard deviation produced by the CE- and PDE-based methods with parameters of the  Table~\ref{table:DWlargesigma} and with the bootstrap number $B=10^4$.}
	\label{fig:DWrobustnessPDE}
\end{figure}

\subsection{Langevin dynamics} \label{ssec:Langevin}

\begin{figure}[h!]
	\hspace*{-0.6cm}
	\includegraphics[height=4.5cm, width=4.5cm]{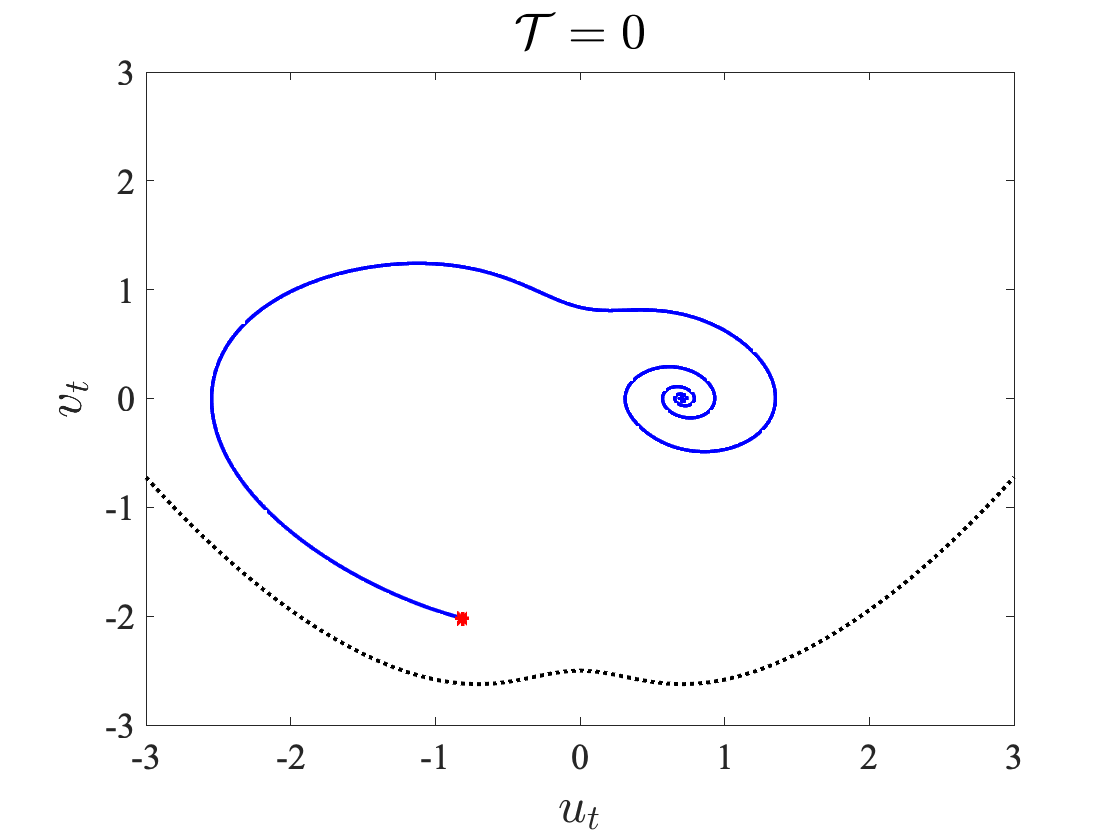}
	\hspace*{-0.6cm}
	\includegraphics[height=4.5cm, width=4.5cm]{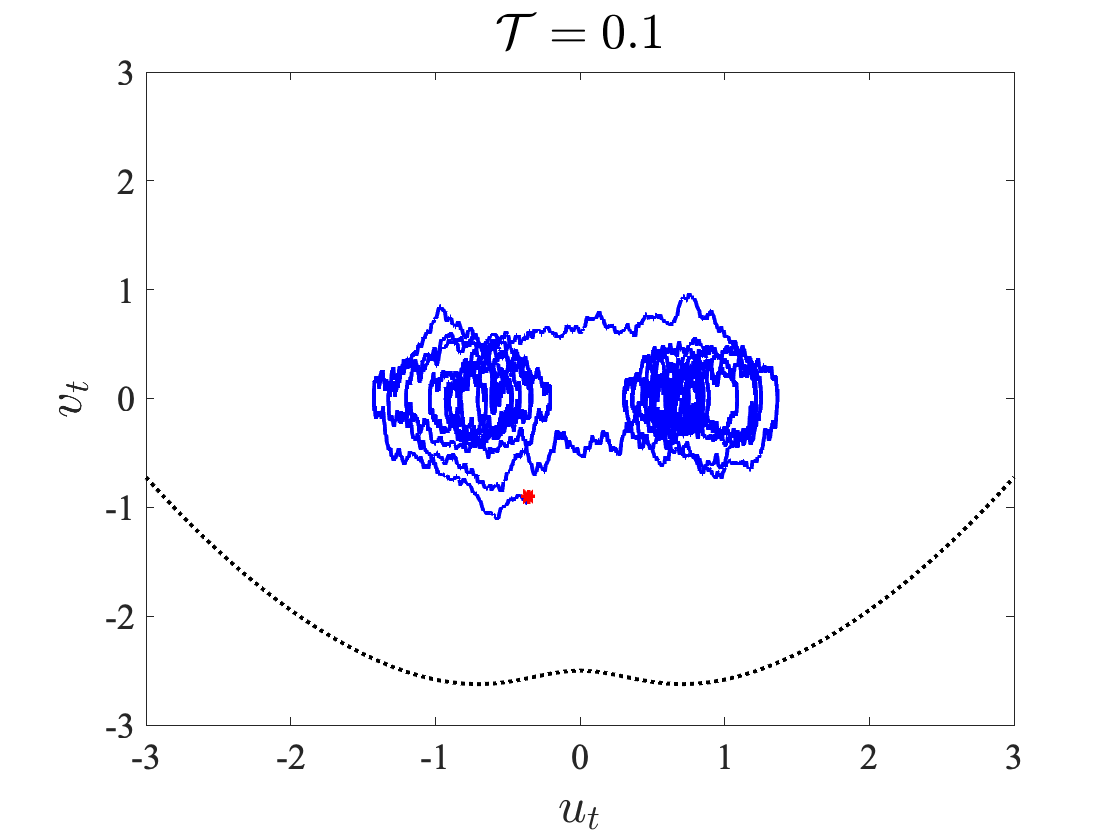}
	\hspace*{-0.6cm}
	\includegraphics[height=4.5cm, width=4.5cm]{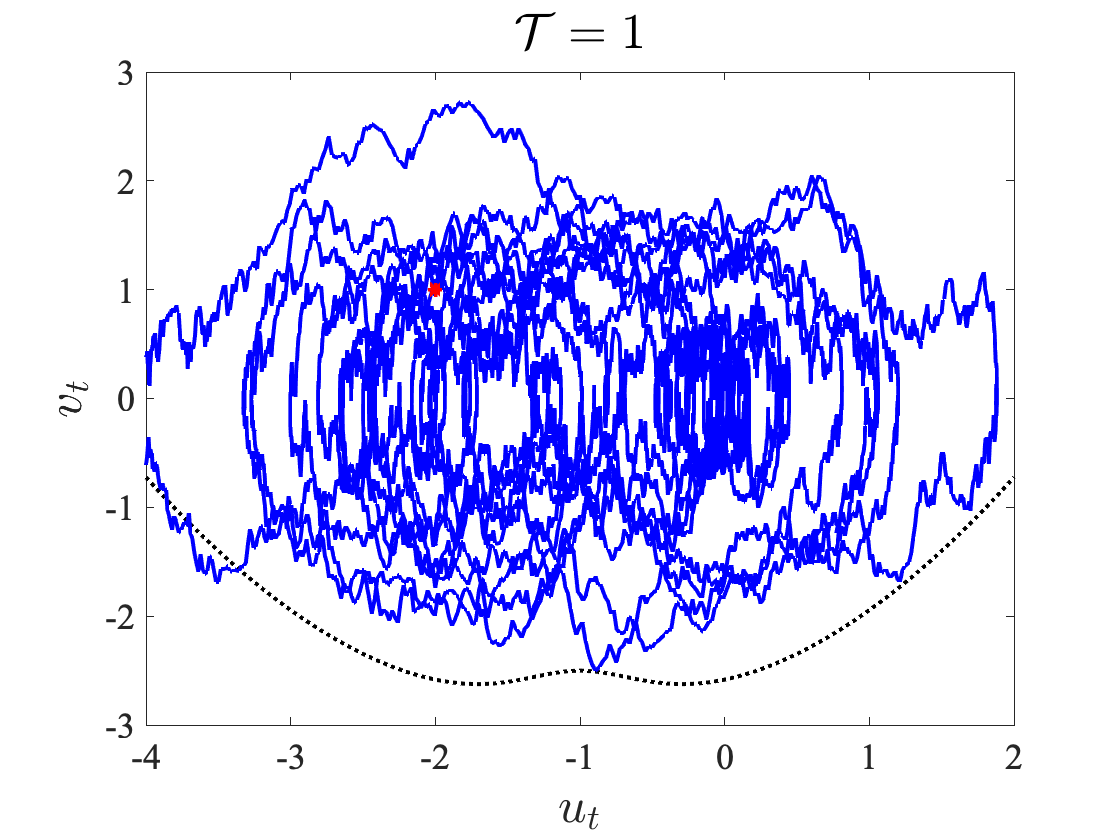}
	\caption{Example: trajectories of the Langevin dynamics with different
temperature $\cT$ values up to time $T=100$. The red
star represents the initial value and the dashed black line represents the double well potential. }
\label{fig:LangevinTraj}
\end{figure}


We consider	two-dimensional Langevin dynamics 
\begin{equation}\label{Langevin}
	\begin{cases} 
		du_t= v_t dt, \\
		dv_t = -V'(u_t)dt-\kappa v_t dt + (2\kappa \cT)^{1/2}dW_t, \;
	\end{cases}\,
\end{equation}
where $V$ is the same double-well potential as in Section~\ref{ssec:dw}. Figure~\ref{fig:LangevinTraj} illustrates trajectories of the Langevin dynamics in different temperature $\cT$ settings. The rare event threshold is defined by the velocity, so that the projection is $P_1=[0 \; 1]$.  In Figure~\ref{fig:Langevinrho0}, a change of the initial density towards the importance region with $\cK=3$ is shown for, respectively, PDE- and CE-based IS techniques wrt $\rho_0$. 

\begin{figure}[h!]
	\hspace*{-0.65cm}
	\includegraphics[height=6.5cm, width=6.5cm]{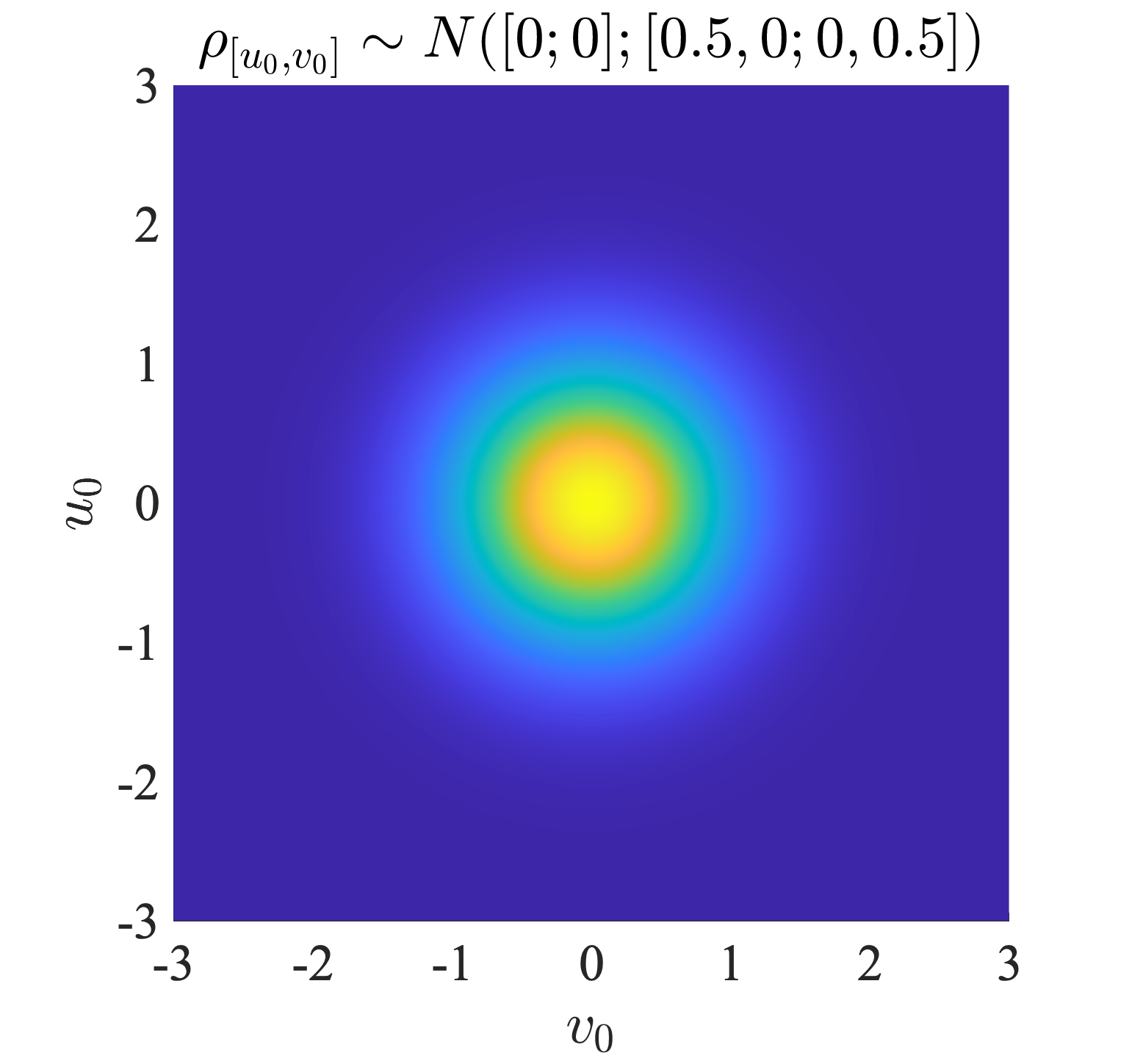}
	\hspace*{-0.6cm}
	\includegraphics[height=6.5cm, width=6.5cm]{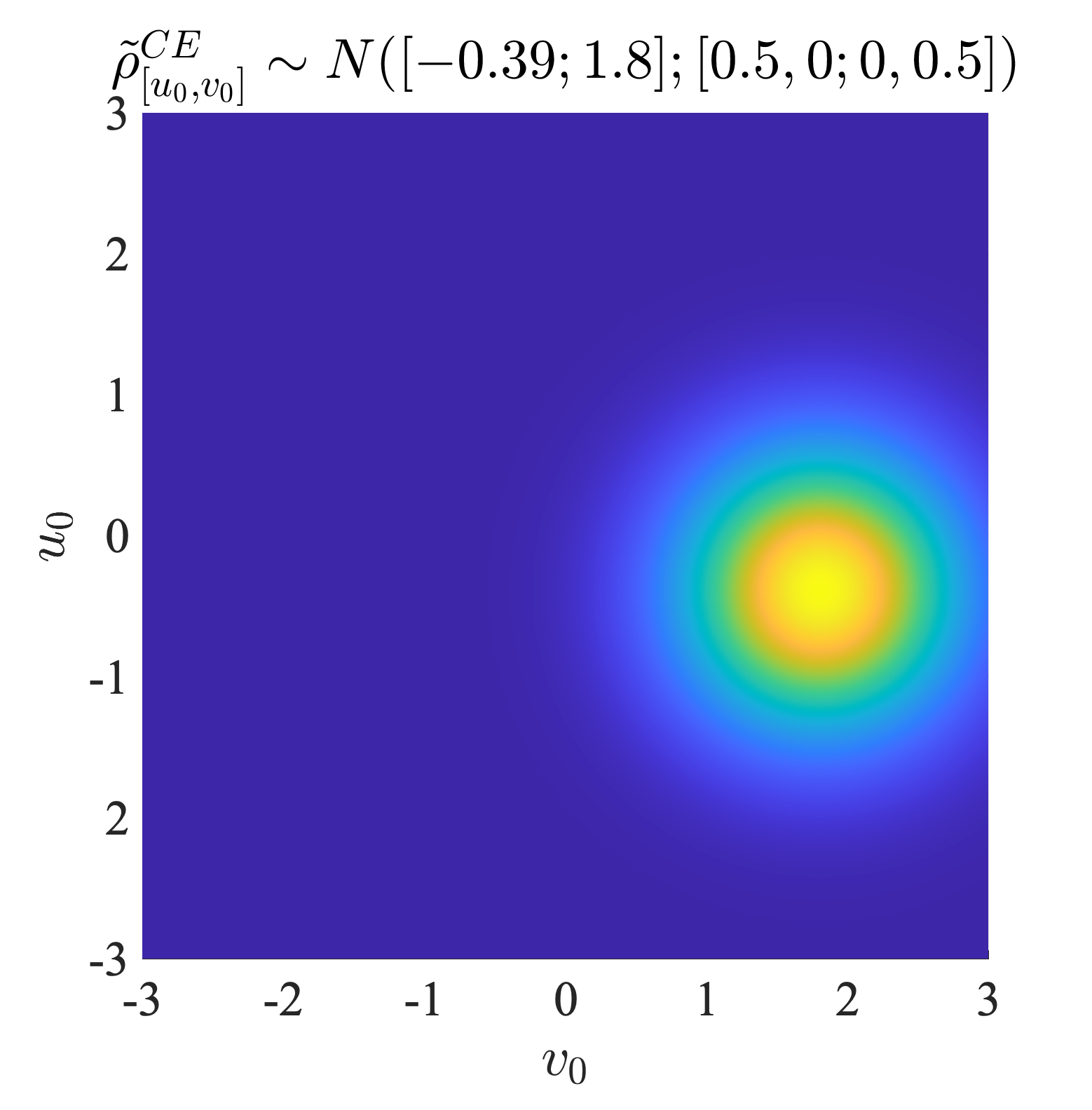}\\
	\hspace*{-0.65cm}
	\includegraphics[height=6.5cm, width=6.5cm]{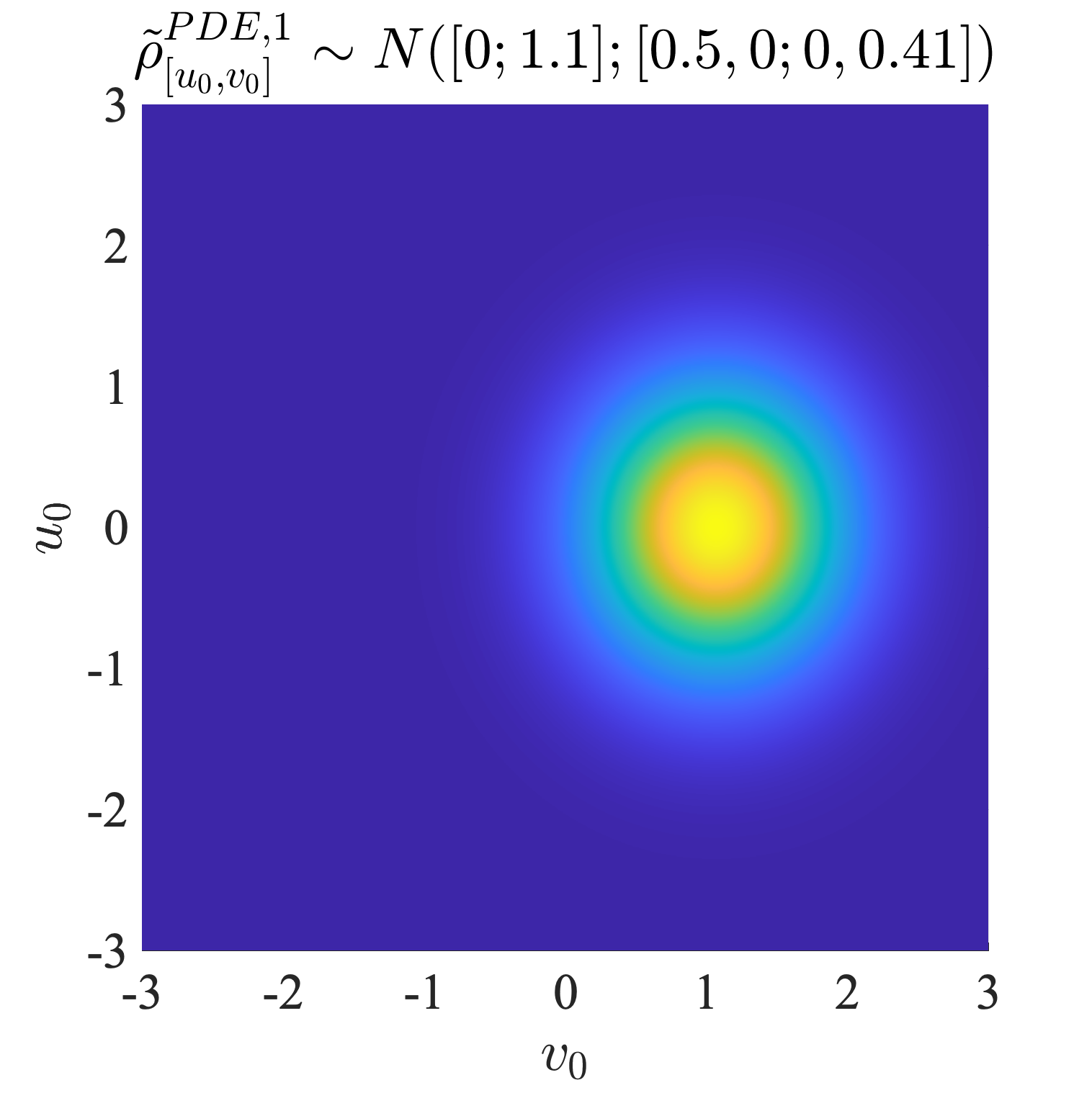}
	\hspace*{-0.6cm}
	\includegraphics[height=6.5cm, width=6.5cm]{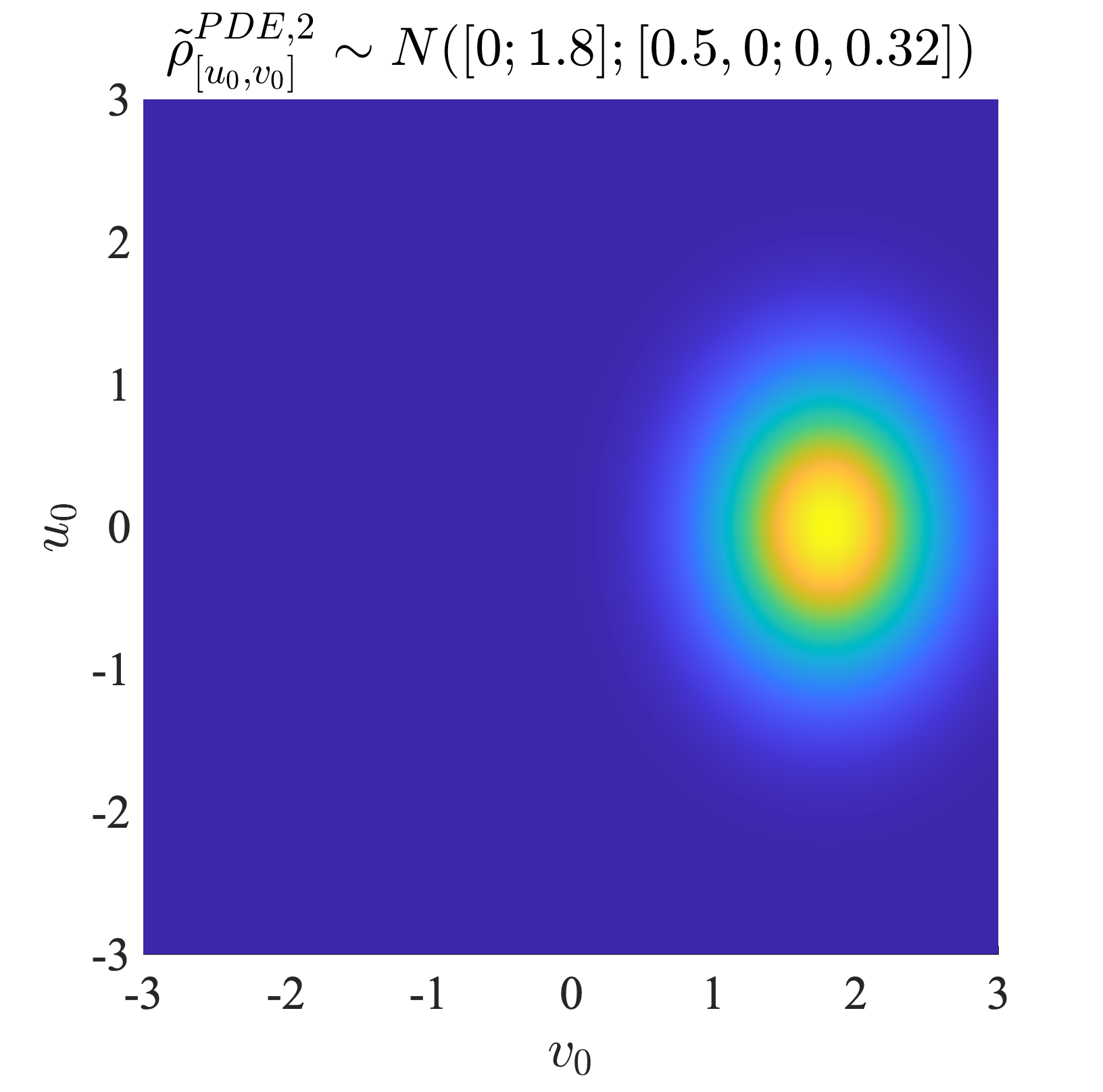}
	
	\caption{\textbf{Langevin dynamics.} A comparison of the original initial denstiy (top left) with the approximately optimal IS initial densities based on the multilevel CE method (top right), the PDE approaches wrt only initial condition (bottom left) and wrt both initial condition and Wiener process (bottom right). }
	\label{fig:Langevinrho0}
\end{figure}

Similarly to the previous DW example, from Table~\ref{table:Langevin}, we can observe that the combination of the PDE-based method wrt the initial distribution and Wiener paths reduces the variance significantly compared to its single approaches and the multilevel CE approach. Figure~\ref{fig:langevinci_corr} shows the consistency of the estimates produced by different methods for the case of correlated initial condition $\begin{bmatrix}
	1 & 0.7 \\
	0.7 & 1
\end{bmatrix}$, where the $95\%$ confidence intervals (CI) overlap and converge as the MC sample size increases, and the IS wrt both $\rho_0$ and $W_t$ provides the highest level of variance reduction when compared to alternative methods. Once again, Figure~\ref{fig:LangevinrobustnessPDE} exhibits a robustness of methods by displaying $95\%$ CI of two independent pilot runs for the parameter setting indicated in Table~\ref{table:Langevin}.

\begin{table}[h!]
	\begin{adjustbox}{width=\columnwidth,center}
		\begin{tabular}{|c|c|c|c|c|l|}
			\hline
			\multirow{2}{*}{$\cK$} & \multicolumn{5}{c|}{\textbf{Estimator with 95\% CI} } \\ \cline{2-6}  
			&  \multicolumn{1}{c|}{\multirow{2}{*}{$\hat{\alpha}^{MC}$}} & \multicolumn{1}{c|}{\multirow{2}{*}{$\hat{\alpha}^{CE,\tilde{\rho}_{0}}$} }& \multicolumn{1}{c|}{\multirow{2}{*}{$\hat{\alpha}^{PDE,\tilde{\rho}_{0}}$}} & \multicolumn{1}{c|}{\multirow{2}{*}{$\hat{\alpha}^{PDE,W_t}$} }&\multicolumn{1}{c|}{\multirow{2}{*}{$\hat{\alpha}^{PDE,\text{both}}$} }\\  [2ex]\hline
			\multirow{2}{*}{$\mathbf{2}$}    & \multicolumn{1}{c|}{$\mathbf{2.53 \times 10^{-2}}$}& \multicolumn{1}{c|}{$\mathbf{2.55 \times 10^{-2}}$}     &  \multicolumn{1}{c|}{$\mathbf{2.53 \times 10^{-2}}$} & \multicolumn{1}{c|}{$\mathbf{2.54 \times 10^{-2}}$} & \multicolumn{1}{c|}{$\mathbf{2.56 \times 10^{-2}}$} \\
			& $[0.025, 0.026]$& $[0.0252, 0.0256]$&$[0.0251, 0.0255]$ &$[0.0252, 0.0256]$&$[0.02559, 0.02569]$\\  \hline
			
			\multirow{2}{*}{$\mathbf{2.5}$}    & \multicolumn{1}{c|}{$\mathbf{4.2 \times 10^{-3}}$}& \multicolumn{1}{c|}{$\mathbf{4.2 \times 10^{-3}}$}     &  \multicolumn{1}{c|}{$\mathbf{4.25 \times 10^{-3}}$} & \multicolumn{1}{c|}{$\mathbf{4.17 \times 10^{-3}}$} & \multicolumn{1}{c|}{$\mathbf{4.27\times 10^{-3}}$} \\
			& $[0.0041, 0.0043]$& $[0.0041, 0.0043]$&$[0.0042, 0.0043]$ &$[0.0041, 0.0042]$&$[0.00426, 0.00428]$\\  \hline
			\multirow{2}{*}{$\mathbf{3}$}    & \multicolumn{1}{c|}{$\mathbf{4.7 \times 10^{-4}}$}& \multicolumn{1}{c|}{$\mathbf{4.6 \times 10^{-4}}$}     &  \multicolumn{1}{c|}{$\mathbf{4.7 \times 10^{-4}}$} & \multicolumn{1}{c|}{$\mathbf{4.8 \times 10^{-4}}$} & \multicolumn{1}{c|}{$\mathbf{4.68 \times 10^{-4}}$} \\
			& $[4.3\text{e-}04, 5.1\text{e-}04]$& $[4.5\text{e-}04, 4.8\text{e-}04]$&$[4.5\text{e-}04, 4.8\text{e-}04]$ &$[4.6\text{e-}04, 4.9\text{e-}04]$&$[4.67\text{e-}04, 4.69\text{e-}04]$\\  \hline
			\multirow{2}{*}{$\mathbf{3.5}$}    & \multicolumn{1}{c|}{$\mathbf{3.8 \times 10^{-5}}$}& \multicolumn{1}{c|}{$\mathbf{3.3 \times 10^{-5}}$}     &  \multicolumn{1}{c|}{$\mathbf{3.4 \times 10^{-5}}$} & \multicolumn{1}{c|}{$\mathbf{3.3 \times 10^{-5}}$} & \multicolumn{1}{c|}{$\mathbf{3.43 \times 10^{-5}}$}\\
			& $[2.6\text{e-}05, 5.0\text{e-}05]$& $[3.0\text{e-}05, 3.6\text{e-}05]$&$[3.1\text{e-}05, 3.7\text{e-}05]$ &$[3.1\text{e-}05, 3.6\text{e-}05]$&$[3.42\text{e-}05, 3.44\text{e-}05]$\\  \hline
		\end{tabular}
	\end{adjustbox}
\newline
\vspace*{0.2 cm}
	\begin{adjustbox}{width=\columnwidth,center}
		\begin{tabular}{|c|ccccl|l|l|l|l|}
			\hline
			\multirow{2}{*}{$\cK$}     & \multicolumn{5}{c|}{\textbf{Relative statistical error}}   &\multicolumn{4}{c|}{\textbf{Variance reduction}}\\\cline{2-10}  
			& \multicolumn{1}{c|}{\multirow{2}{*}{$\epsilon_{st}^{MC}$} }& \multicolumn{1}{c|}{\multirow{2}{*}{$\epsilon_{st}^{CE,\tilde{\rho}_{0}}$}} & \multicolumn{1}{c|}{\multirow{2}{*}{$\epsilon_{st}^{PDE,\tilde{\rho}_{0}}$} }&\multicolumn{1}{c|}{\multirow{2}{*}{$\epsilon_{st}^{PDE,W_t}$} }& \multicolumn{1}{c|}{\multirow{2}{*}{$\epsilon_{st}^{PDE,\text{both}}$} }& \multirow{2}{*}{$\frac{\mathbb{V}_{MC}}{\mathbb{V}_{CE,\tilde{\rho}_{0}}}$}& \multirow{2}{*}{$\frac{\mathbb{V}_{MC}}{\mathbb{V}_{PDE,\tilde{\rho}_{0}}}$} & \multirow{2}{*}{$\frac{\mathbb{V}_{MC}}{\mathbb{V}_{PDE,W_t}}$}&\multirow{2}{*}{$\frac{\mathbb{V}_{MC}}{\mathbb{V}_{PDE, \text{both}}}$}\\  [2ex]\hline
			$\mathbf{2}$          & \multicolumn{1}{c|}{$\mathbf{1.21\%}$}     &  \multicolumn{1}{c|}{$\mathbf{1.18\%}$} & \multicolumn{1}{c|}{$\mathbf{0.8\%}$} & \multicolumn{1}{c|}{$\mathbf{0.7\%}$}&\multicolumn{1}{c|}{$\mathbf{0.2\%}$} &\multicolumn{1}{c|}{$\mathbf{1.1} $} & \multicolumn{1}{c|}{$\mathbf{2.6} $} &\multicolumn{1}{c|}{$\mathbf{3.5} $} &\multicolumn{1}{c|}{$\mathbf{35} $} \\ \hline
			$\mathbf{2.5}$          & \multicolumn{1}{c|}{$\mathbf{3.0\%}$}     &  \multicolumn{1}{c|}{$\mathbf{2.5\%}$} & \multicolumn{1}{c|}{$\mathbf{1.5\%}$} & \multicolumn{1}{c|}{$\mathbf{1.5\%}$}&\multicolumn{1}{c|}{$\mathbf{0.2\%}$} &\multicolumn{1}{c|}{$\mathbf{1.3} $} & \multicolumn{1}{c|}{$\mathbf{3.9} $} &\multicolumn{1}{c|}{$\mathbf{5.8} $} &\multicolumn{1}{c|}{$\mathbf{180} $} \\ \hline
			$\mathbf{3}$          & \multicolumn{1}{c|}{$\mathbf{9.1\%}$}     &  \multicolumn{1}{c|}{$\mathbf{3.8\%}$} & \multicolumn{1}{c|}{$\mathbf{3.5\%}$} & \multicolumn{1}{c|}{$\mathbf{2.9\%}$}&\multicolumn{1}{c|}{$\mathbf{0.2\%}$} &\multicolumn{1}{c|}{$\mathbf{5.7} $} & \multicolumn{1}{c|}{$\mathbf{6.6} $} &\multicolumn{1}{c|}{$\mathbf{9.3} $} &\multicolumn{1}{c|}{$\mathbf{2000} $} \\ \hline
			$\mathbf{3.5}$          & \multicolumn{1}{c|}{$\mathbf{31.8\%}$}     &  \multicolumn{1}{c|}{$\mathbf{8.9\%}$} & \multicolumn{1}{c|}{$\mathbf{9.1\%}$} & \multicolumn{1}{c|}{$\mathbf{7.3\%}$}&\multicolumn{1}{c|}{$\mathbf{0.3\%}$} &\multicolumn{1}{c|}{$\mathbf{15.2} $} & \multicolumn{1}{c|}{$\mathbf{15.4} $} &\multicolumn{1}{c|}{$\mathbf{24.7} $} &\multicolumn{1}{c|}{$\mathbf{12156} $}  \\ \hline
		\end{tabular}
	\end{adjustbox}
	\caption{\textbf{ Langevin dynamics.} The model parameters: $\kappa=2^{-5}\pi^2$, $\cT=1$. Simulation parameters: $T=1$, $\Delta t= 0.01$, $J=10^6$, $[u_0,v_0] \sim N([0\;0], [0.5\; 0; 0\; 0.5])$. PDE-based methods uses $\Delta x_{\mathbf{PDE}}=0.006$, $\Delta t_{\mathbf{PDE}}=\Delta x_{\mathbf{PDE}}/5$.}
	\label{table:Langevin}
\end{table}

\begin{figure}[h!]
	\includegraphics[height=6cm, width=8cm]{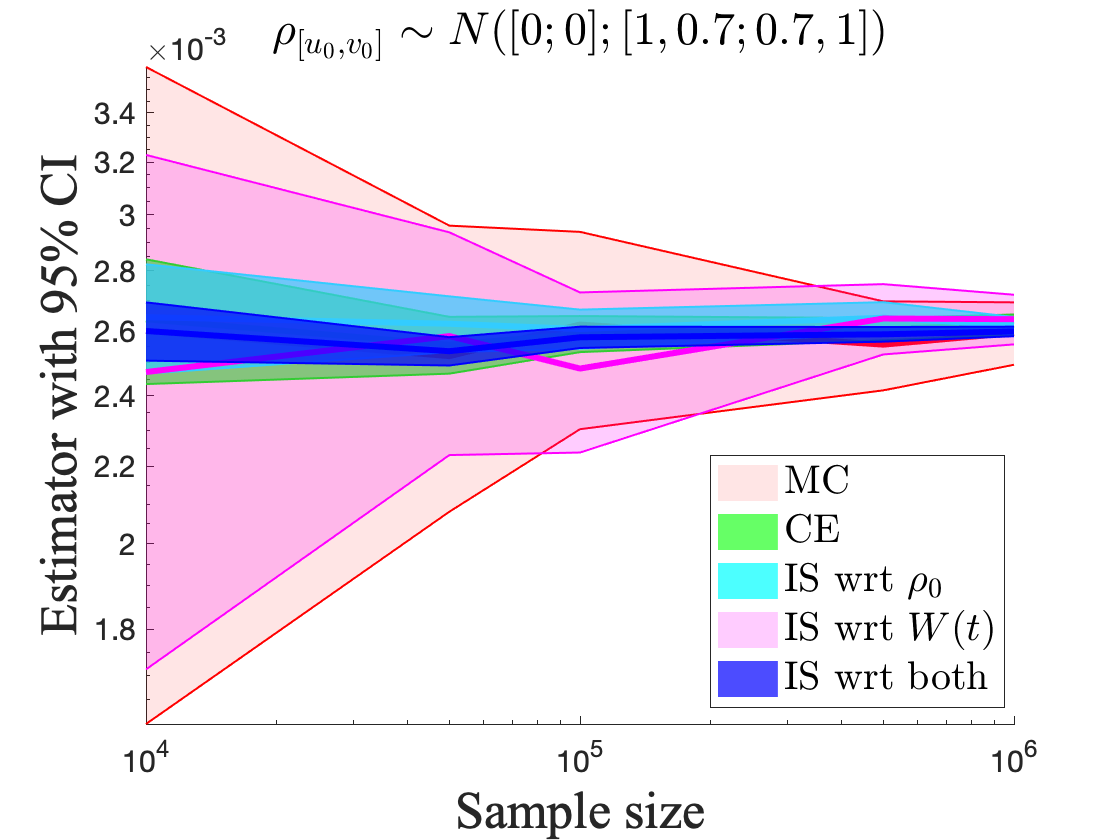}
	\caption{\textbf{ Langevin dynamics.} The model parameters: $\kappa=2^{-5}\pi^2$, $\cT=1$. Simulation parameters: $T=1$, $\Delta t= 0.01$, $\cK=2.5$, $[u_0,v_0] \sim N([0\;0], [1\; 0.7; 0.7\; 1])$. PDE-based methods uses $\Delta x_{\mathbf{PDE}}=0.006$, $\Delta t_{\mathbf{PDE}}=\Delta x_{\mathbf{PDE}}/5$.}
	\label{fig:langevinci_corr}
\end{figure}

\begin{figure}[h!]
	\hspace*{-0.6cm}
	\includegraphics[height=6.3cm, width=6.3cm]{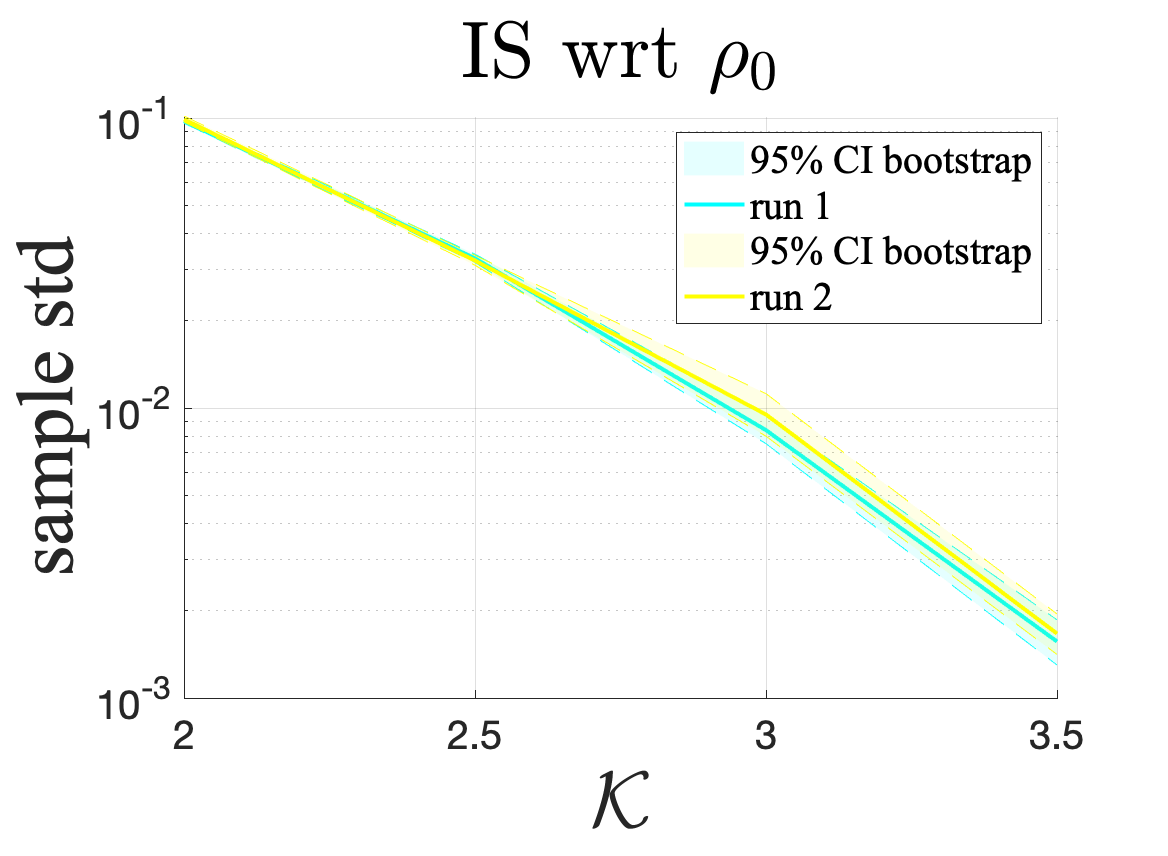}
	\hspace*{-0.1cm}
	\includegraphics[height=6.3cm, width=6.3cm]{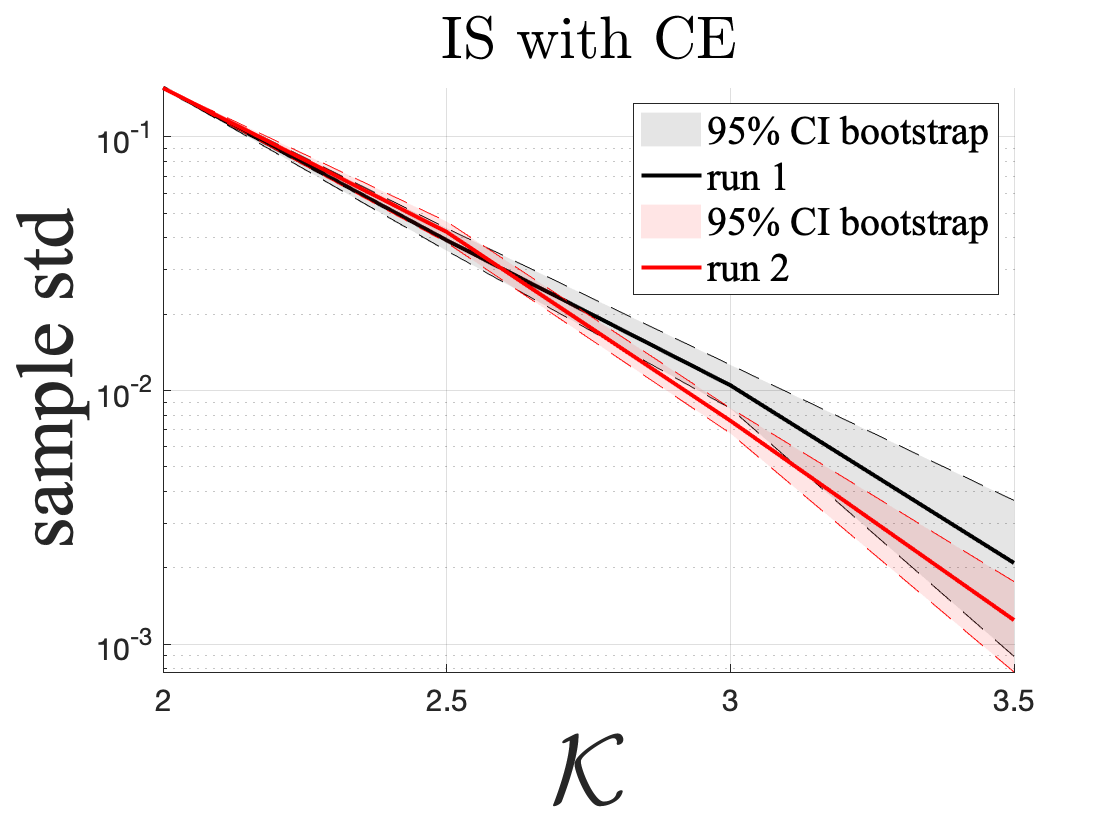}\\
	\hspace*{-0.6cm}
	\includegraphics[height=6.3cm, width=6.3cm]{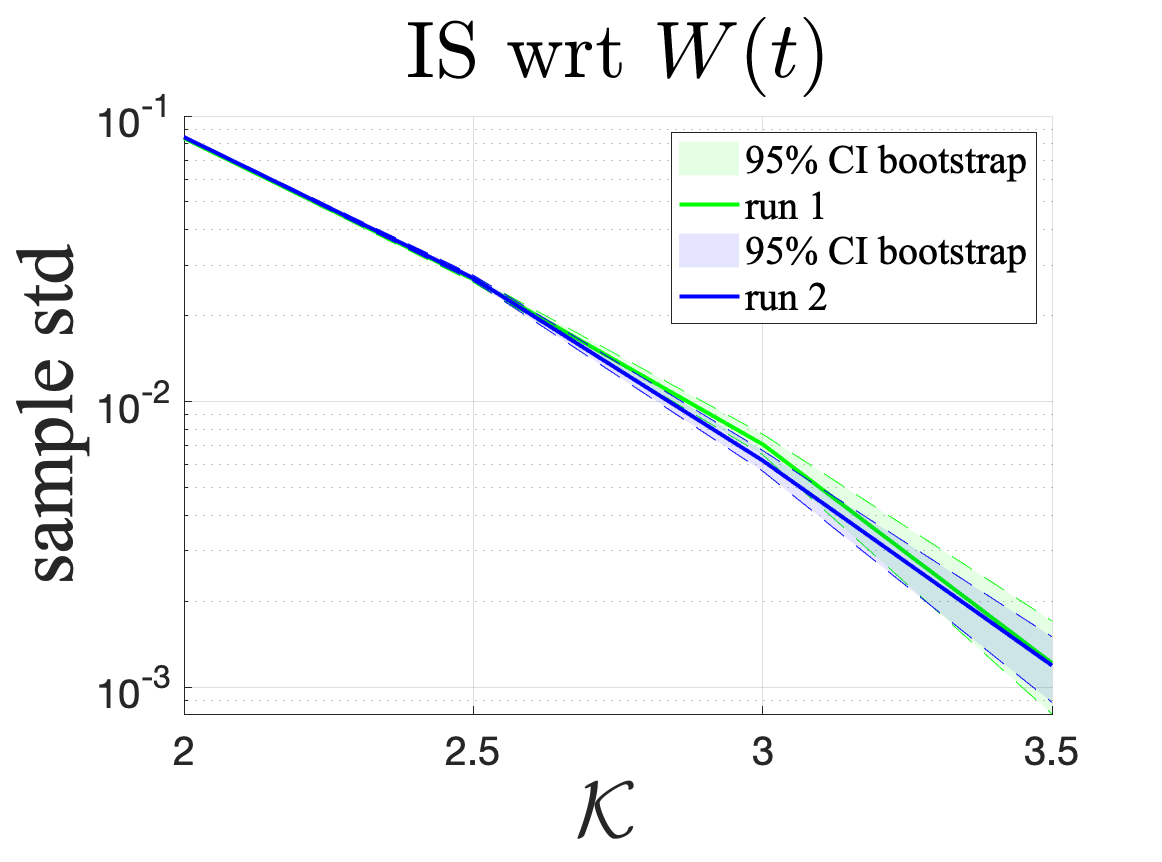}
	\hspace*{-0.1cm}
	\includegraphics[height=6.3cm, width=6.3cm]{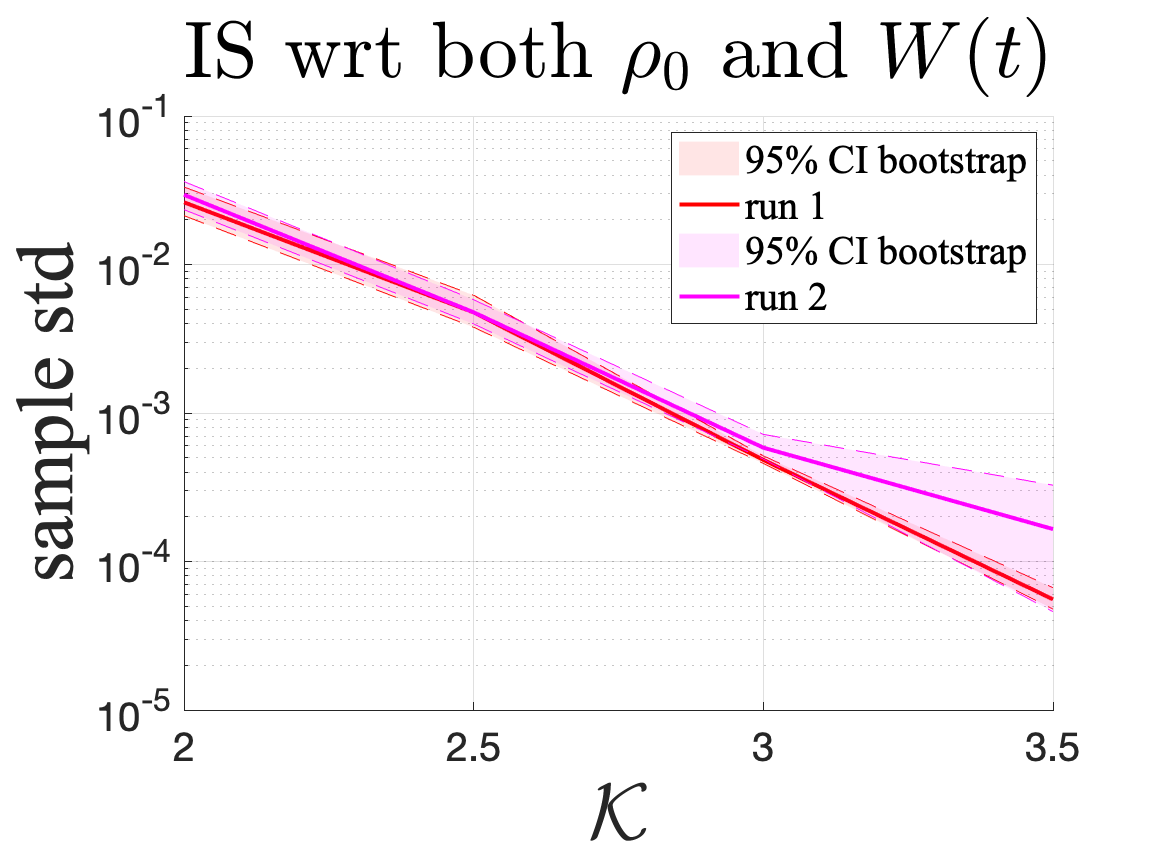}
	\caption{\textbf{Langevin dynamics.} Two independent simulation results of the bootstrapping for a sample standard deviation produced by the CE- and PDE-based methods with parameters of the  Table~\ref{table:Langevin} and with the bootstrap number $B=10^4$.}
	\label{fig:LangevinrobustnessPDE}
\end{figure}

\subsection{A noisy Charney-deVore model}\label{ssec:cdv}
We consider the six-dimensional Charney–deVore (CdV) model polluted with additive Gaussian noise~\cite{grafke2019numerical}
\begin{equation}\label{eq:CharneyDevore}
	\begin{cases} 
		du_1=(\tilde{\gamma}_1u_3-C(u_1-u_1^*))dt+ \sqrt{2b}dW_1, \\
		du_2=(-(\alpha_1u_1-\beta_1)u_3-Cu_2-\delta_1u_4u_6)dt+\sqrt{2b}dW_2,\\
		du_3=((\alpha_1u_1-\beta_1)u_2 - \gamma_1u_1-Cu_3+\delta_1u_4u_5)dt+\sqrt{2b}dW_3,\;\\
		du_4=(\tilde{\gamma}_2u_6-C(u_4-u_4^*)+\eta(u_2u_6-u_3u_5))dt+\sqrt{2b} dW_4,\;\\
		du_5=(-(\alpha_2u_1-\beta_2)u_6-Cu_5-\delta_2u_3u_4)dt+\sqrt{2b}dW_5,\;			\\	
		du_6=((\alpha_2u_1-\beta_2)u_5-\gamma_2u_4-Cu_6+\delta_2u_2u_4)dt+\sqrt{2b}dW_6,\;
	\end{cases}\,
\end{equation}
where for $m\in{1,2}$,
\begin{equation*}
	\begin{split}
		&\alpha_m=\frac{8\sqrt{2}}{\pi}\frac{m^2}{4m^2-1}\frac{q^2+m^2-1}{b^2+m^2}, \quad \beta_m=\frac{\beta q^2}{q^2+m^2},\\
		&\gamma_m=\gamma\frac{ \sqrt{2}q}{\pi}\frac{4m^3}{(4m^2-1)(q^2+m^2)}, \quad \tilde{\gamma}_m=\gamma \frac{\sqrt{2}q}{\pi}\frac{4m}{4m^2-1},\\
		&\delta_m =\frac{64\sqrt{2}}{15\pi}\frac{q^2-m^2+1}{q^2+m^2}, \quad \eta=\frac{16\sqrt{2}}{5\pi}.
	\end{split}
\end{equation*}
The CdV was among the pioneering atmospheric models that demonstrated multiple invariant measures. The model~\eqref{eq:CharneyDevore} features two distinct metastable modes: one is the "zonal" mode and the other is the "blocked" mode. These modes correspond to atmospheric blocking phenomena commonly observed in meteorology~\cite{charney1979multiple, crommelin2004mechanism, majda2012filtering}.

In all numerical tests, we use the parameters provided in~\cite{crommelin2004mechanism} which result in a chaotic behavior of dynamics: zonally symmetric forcing terms $u_1^*=0.95$ and $u_4^*=-0.76095$; the thermal relaxation damping timescale $C=0.1$ corresponding to 10 days; the orographic height $\gamma=0.2$ corresponding to a $200 \mathrm{m}$ amplitude; the Coriolis parameter $\beta=1.25$ defining a central latitude of $45^o$; the channel width–length ratio $q=0.5$ corresponding to a channel of $6300\mathrm{km}\times 1600\mathrm{km}$. For simplicity, we consider the rare event defined by the first component, i.e., the projection $P_1=[1\; 0\; 0\; 0\; 0\; 0]$. We set the initial condition $\boldsymbol{u}_0:=[u_1(0)\; u_2(0)\; u_3(0)\; u_4(0)\; u_5(0)\; u_6(0)]$ such that the first component is sampled from the blocked region
$$\boldsymbol{u}_0\sim N([0.7650, 0.2288, -0.2990, -0.3657, -0.1636, 0.3108], \sigma_0^2I)$$ 
where $I$ is a 6-dimensional identity matrix and the rare event threshold is defined from the zonal region. The metastable feature of the dynamics can be observed in Figure~\ref{fig:CdVtraj}, which illustrates a sample trajectory of the tracked component $u_1$ for a given different values of the diffusion parameter $b=0,\, 0.0001,\, 0.001,\, 0.01$ over time $T=1000$. 

\begin{figure}[h!]
	\includegraphics[height=7.5cm, width=12cm]{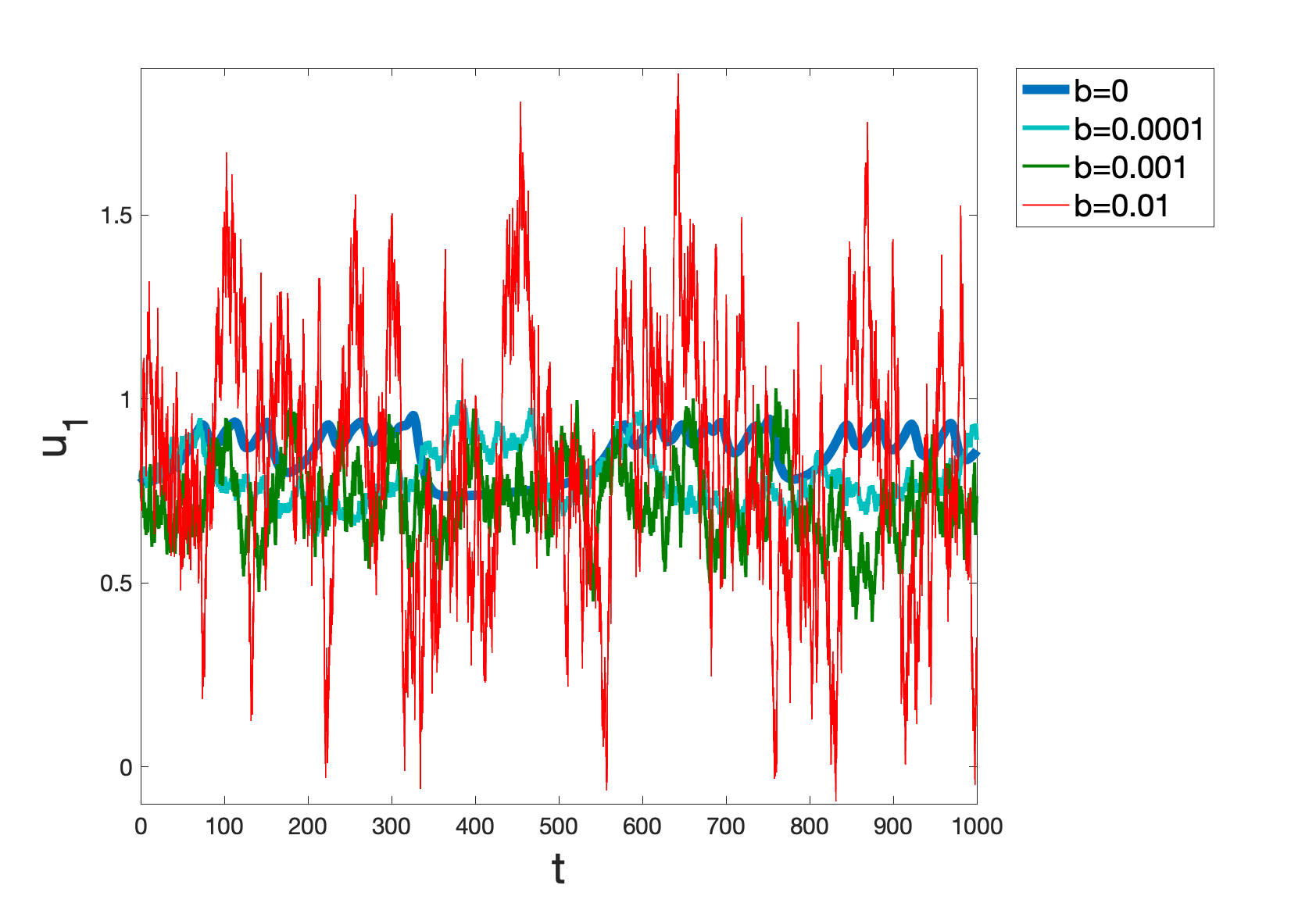}
	\caption{\textbf{Charney-deVore model.} Sample trajectories of the first component corresponding to different values of the diffusion parameter $b$ over $T=1000$ time with $\Delta t=0.01$.}
	\label{fig:CdVtraj}
\end{figure}

	\begin{table}[h!]
	\begin{adjustbox}{width=\columnwidth,center}
		\begin{tabular}{|c|c|c|c|c|l|}
			\hline
			\multirow{2}{*}{$\cK$} & \multicolumn{5}{c|}{\textbf{Estimator with 95\% CI} } \\ \cline{2-6}  
			&  \multicolumn{1}{c|}{\multirow{2}{*}{$\hat{\alpha}^{MC}$}} & \multicolumn{1}{c|}{\multirow{2}{*}{$\hat{\alpha}^{CE,\tilde{\rho}_{0}}$} }& \multicolumn{1}{c|}{\multirow{2}{*}{$\hat{\alpha}^{PDE,\tilde{\rho}_{0}}$}} & \multicolumn{1}{c|}{\multirow{2}{*}{$\hat{\alpha}^{PDE,W_t}$} }&\multicolumn{1}{c|}{\multirow{2}{*}{$\hat{\alpha}^{PDE,\text{both}}$} }\\  [2ex]\hline
			\multirow{2}{*}{$\mathbf{1.1}$}  & \multicolumn{1}{c|}{$\mathbf{3.1 \times 10^{-2}}$}& \multicolumn{1}{c|}{$\mathbf{3.1 \times 10^{-2}}$}     &  \multicolumn{1}{c|}{$\mathbf{3.1 \times 10^{-2}}$} & \multicolumn{1}{c|}{$\mathbf{3.08\times 10^{-2}}$} & \multicolumn{1}{c|}{$\mathbf{3.08 \times 10^{-2}}$} \\
			& $[0.0308, 0.0314]$&$[0.0308, 0.0315]$ &$[0.0308, 0.0314]$ &$[0.0307, 0.0309]$&$[0.0308, 0.0309]$\\  \hline
			
			\multirow{2}{*}{$\mathbf{1.2}$}    & \multicolumn{1}{c|}{$\mathbf{4.2\times 10^{-3}}$}& \multicolumn{1}{c|}{$\mathbf{4.2\times 10^{-3}}$}     &  \multicolumn{1}{c|}{$\mathbf{4.2 \times 10^{-3}}$} & \multicolumn{1}{c|}{$\mathbf{4.1 \times 10^{-3}}$} & \multicolumn{1}{c|}{$\mathbf{4.11 \times 10^{-3}}$} \\
			& $[0.0041, 0.0043]$& $[0.0040, 0.0043]$&$[0.0041, 0.0043]$ &$[0.00409, 0.00412]$&$[0.00411, 0.00412]$\\  \hline
			\multirow{2}{*}{$\mathbf{1.3}$}    & \multicolumn{1}{c|}{$\mathbf{3.3 \times 10^{-4}}$}& \multicolumn{1}{c|}{$\mathbf{3.4 \times 10^{-4}}$}     &  \multicolumn{1}{c|}{$\mathbf{3.44 \times 10^{-4}}$} & \multicolumn{1}{c|}{$\mathbf{3.5 \times 10^{-4}}$} & \multicolumn{1}{c|}{$\mathbf{3.45 \times 10^{-4}}$} \\
			& $[2.9\text{e-}04, 3.7\text{e-}04]$& $[3.0\text{e-}04, 3.8\text{e-}04]$ &$[3.2\text{e-}04, 3.8\text{e-}04]$ &$[3.42\text{e-}04, 3.46\text{e-}04]$&$[3.44\text{e-}04, 3.46\text{e-}04]$\\  \hline
			\multirow{2}{*}{$\mathbf{1.4}$}    & \multicolumn{1}{c|}{$\mathbf{2.1\times 10^{-5}}$}& \multicolumn{1}{c|}{$-$}     &  \multicolumn{1}{c|}{$\mathbf{1.5 \times 10^{-5}}$} & \multicolumn{1}{c|}{$\mathbf{1.8\times 10^{-5}}$} & \multicolumn{1}{c|}{$\mathbf{1.81\times 10^{-5}}$}\\
			& $[1.2\text{e-}05, 3.0\text{e-}05]$& &$[8.9\text{e-}06, 2.2\text{e-}05]$ &$[1.78\text{e-}05, 1.81\text{e-}05]$&$[1.80\text{e-}05, 1.81\text{e-}05]$\\  \hline
		\end{tabular}
	\end{adjustbox}
	\newline
	\vspace*{0.2 cm}
	\begin{adjustbox}{width=\columnwidth,center}
		\begin{tabular}{|c|ccccl|l|l|l|l|}
			\hline
			\multirow{2}{*}{$\cK$}     & \multicolumn{5}{c|}{\textbf{Relative statistical error}}   &\multicolumn{4}{c|}{\textbf{Variance reduction}}\\\cline{2-10}  
			& \multicolumn{1}{c|}{\multirow{2}{*}{$\epsilon_{st}^{MC}$} }& \multicolumn{1}{c|}{\multirow{2}{*}{$\epsilon_{st}^{CE,\tilde{\rho}_{0}}$}} & \multicolumn{1}{c|}{\multirow{2}{*}{$\epsilon_{st}^{PDE,\tilde{\rho}_{0}}$} }&\multicolumn{1}{c|}{\multirow{2}{*}{$\epsilon_{st}^{PDE,W_t}$} }& \multicolumn{1}{c|}{\multirow{2}{*}{$\epsilon_{st}^{PDE,\text{both}}$} }& \multirow{2}{*}{$\frac{\mathbb{V}_{MC}}{\mathbb{V}_{CE,\tilde{\rho}_{0}}}$}& \multirow{2}{*}{$\frac{\mathbb{V}_{MC}}{\mathbb{V}_{PDE,\tilde{\rho}_{0}}}$} & \multirow{2}{*}{$\frac{\mathbb{V}_{MC}}{\mathbb{V}_{PDE,W_t}}$}&\multirow{2}{*}{$\frac{\mathbb{V}_{MC}}{\mathbb{V}_{PDE, \text{both}}}$}\\  [2ex]\hline
			$\mathbf{1.1}$          & \multicolumn{1}{c|}{$\mathbf{1.1\%}$}     &  \multicolumn{1}{c|}{$\mathbf{1.1\%}$} & \multicolumn{1}{c|}{$\mathbf{1.0\%}$} & \multicolumn{1}{c|}{$\mathbf{0.2\%}$}&\multicolumn{1}{c|}{$\mathbf{0.1\%}$} &\multicolumn{1}{c|}{$\mathbf{1.0} $} & \multicolumn{1}{c|}{$\mathbf{1.1} $} &\multicolumn{1}{c|}{$\mathbf{21} $} &\multicolumn{1}{c|}{$\mathbf{78} $} \\ \hline
			$\mathbf{1.2}$          & \multicolumn{1}{c|}{$\mathbf{3.0\%}$}     &  \multicolumn{1}{c|}{$\mathbf{3.1\%}$} & \multicolumn{1}{c|}{$\mathbf{2.6\%}$} & \multicolumn{1}{c|}{$\mathbf{0.4\%}$}&\multicolumn{1}{c|}{$\mathbf{0.2\%}$} &\multicolumn{1}{c|}{$\mathbf{1.0} $} & \multicolumn{1}{c|}{$\mathbf{1.3} $} &\multicolumn{1}{c|}{$\mathbf{75} $} &\multicolumn{1}{c|}{$\mathbf{296} $} \\ \hline
			$\mathbf{1.3}$          & \multicolumn{1}{c|}{$\mathbf{11\%}$}     &  \multicolumn{1}{c|}{$\mathbf{12\%}$} & \multicolumn{1}{c|}{$\mathbf{8.5\%}$} & \multicolumn{1}{c|}{$\mathbf{0.5\%}$}&\multicolumn{1}{c|}{$\mathbf{0.2\%}$} &\multicolumn{1}{c|}{$\mathbf{0.8} $} & \multicolumn{1}{c|}{$\mathbf{1.4} $} &\multicolumn{1}{c|}{$\mathbf{394} $} &\multicolumn{1}{c|}{$\mathbf{1985} $} \\ \hline
			$\mathbf{1.4}$          & \multicolumn{1}{c|}{$\mathbf{43\%}$}     &  \multicolumn{1}{c|}{$\mathbf{-\%}$} & \multicolumn{1}{c|}{$\mathbf{42\%}$} & \multicolumn{1}{c|}{$\mathbf{0.7\%}$}&\multicolumn{1}{c|}{$\mathbf{0.3\%}$} &\multicolumn{1}{c|}{$\mathbf{-} $} & \multicolumn{1}{c|}{$\mathbf{2.0} $} &\multicolumn{1}{c|}{$\mathbf{5395} $} &\multicolumn{1}{c|}{$\mathbf{38467} $}  \\ \hline
		\end{tabular}
	\end{adjustbox}
	\caption{\textbf{A noisy Charney-deVore model with $b=0.01$.} Simulation parameters: $\sigma_0^2=0.0025$, $T=1$, $\Delta t= 0.01$,  $J = 10^6$. Numerical results for the last line of the CE method is missing due to the impracticability of the algorithm for the given parameter setting and a limit of computer capacity to accomplish the simulation.}
	\label{table:CharneydeVore}
\end{table}

In a similar fashion to previous examples, Table~\ref{table:CharneydeVore} provides numerical results of the CdV model with the diffusion parameter $b=0.01$. In this case, we observe a considerable variance reduction in the PDE-based IS technique wrt the Wiener process and a much more significant reduction in the IS approach wrt both the initial condition and the Wiener process. Similar to the Double Well problem in Section~\ref{ssec:dw}, it is expected to get a variance reduction in IS wrt the initial condition if much larger values of $\sigma_0^2$ are considered in the initial condition. Furthermore, we provide the bootstrapping results for the sample standard deviation in Figure~\ref{fig:LangevinrobustnessPDE} for the parameter setting indicated in Table~\ref{table:CharneydeVore} as numerical evidence of the robustness of the considered methods. We have to note that when a diffusion coefficient is very small, the simple numerical PDE solvers with uniform discretization (see Appendix~\ref{appx:PDEsolver}) may not provide accurate solutions, and special adaptive mesh refinement strategies may be needed to solve the PDEs numerically. Instead, the large deviation principle-based methods might be more efficient in a very small noise regime~\cite{grafke2019numerical, vanden2013data}. 

\begin{figure}[h!]
	\hspace*{-0.6cm}
	\includegraphics[height=6.3cm, width=6.3cm]{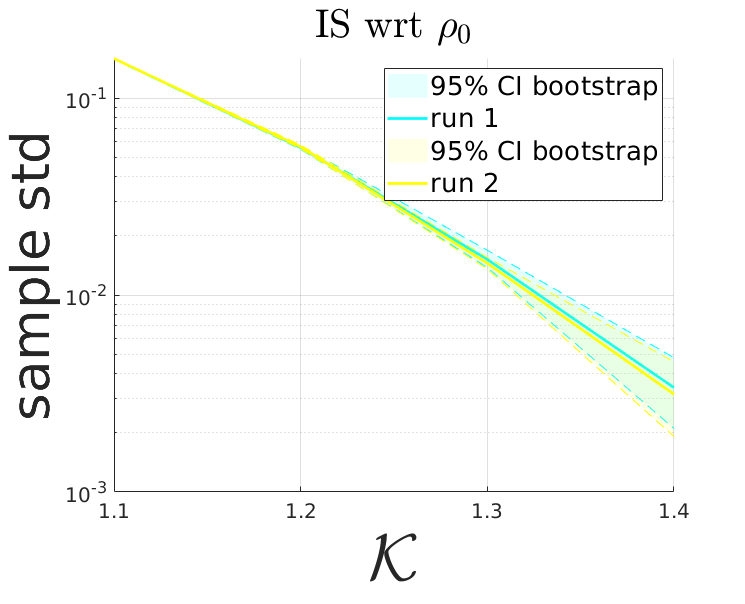}
	\hspace*{-0.1cm}
	\includegraphics[height=6.3cm, width=6.3cm]{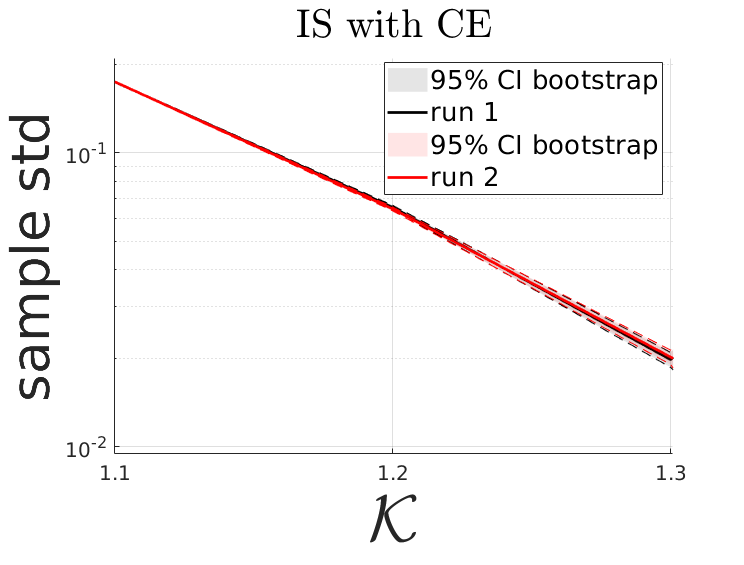}\\
	\hspace*{-0.6cm}
	\includegraphics[height=6.3cm, width=6.3cm]{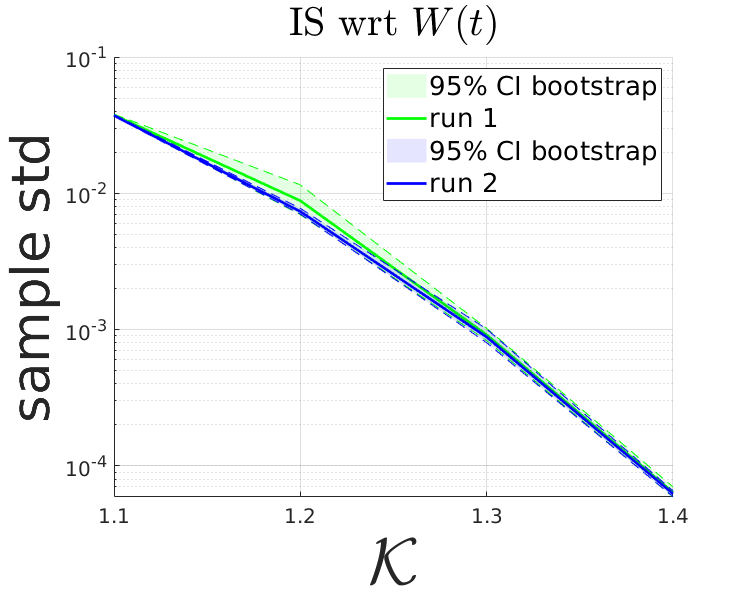}
	\hspace*{-0.1cm}
	\includegraphics[height=6.3cm, width=6.3cm]{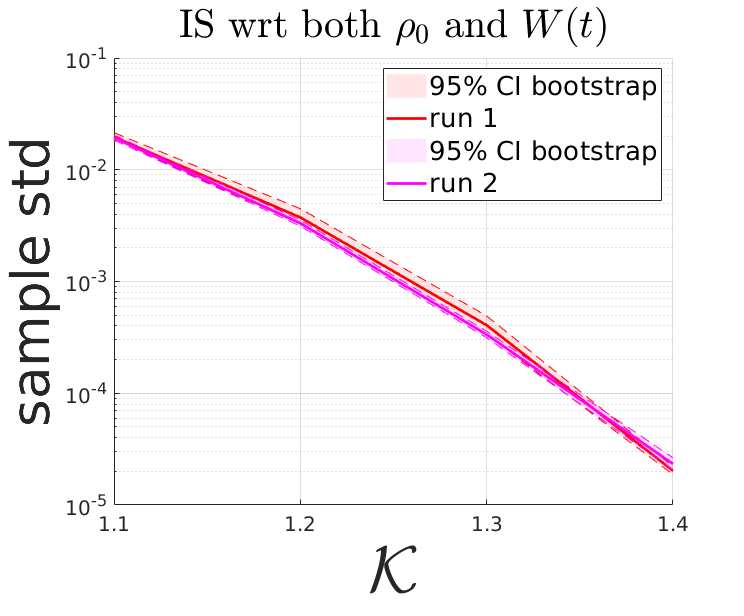}
	\caption{\textbf{Charney-deVore model.} Two independent simulation results of the bootstrapping for a sample standard deviation produced by the CE- and PDE-based methods with parameters of the Table~\ref{table:CharneydeVore} and with the bootstrap number $B=10^4$.}
	\label{fig:CharneydeVoremodel}
\end{figure}

\section{Conclusion}\label{sec:conclusion}
In this study, we made use of IS techniques to monitor a rare event probability exceeding a critical threshold in relation to the running maximum associated with the solution of an SDE. We conducted this investigation within the framework of EnKF. In the interval between two observations in the EnKF method, we put forward three IS strategies: one that changes a measure wrt the SDE's initial condition, another that involves a change of measure wrt the Wiener process using a SOC formulation, and a third that combines these IS approaches for both the initial condition and the Wiener process. All of these IS methods hinge on approximating the solution of the KBE while considering specific boundary conditions. In more than one-dimensional setting, we employed a dimension reduction technique through the Markovian projection. This technique allows us to approximate the KBE solution by solving a simpler one-dimensional PDE. It was demonstrated through numerical examples that the proposed estimation methods significantly reduced variance compared to the standard MC and the multilevel CE techniques.


\appendix

\section{A discrete $L^2$ regression} \label{appx:L2regression}
We want to approximate $\check{a}$ and $\check{b}$, in~\eqref{eq:surrogateDriftAndDiff}, by linear combinations of bivariate polynomials defined as 
\[
\psi_p(t,s) = t^{p_1} s^{p_2}, \qquad p=(p_1, p_2) \in \bN^2.
\]
Restrict $p$ to a tensor product index set $\Lambda(\varpi) \subset \bN^2$ such that $\max\limits_{p\in \Lambda} \max\limits_{k=1,2} p_k\leq\varpi$ and let $\cV_{\Lambda(\varpi)} = \mbox{span} \{\psi_p, p \in\Lambda(\varpi)\}$.
We approximate $\check{a}$ and $\check{b}$ by $\hat{\check{a}}$ and $\hat{\check{b}}$ respectively, obtained by computing
\begin{equation}\label{eq:L2regression}
	\small
	\begin{cases} 
		\hat{\check{a}} =\arg \min\limits_{h_a\in \cV_{\Lambda}} \frac{1}{NM}\sum_{n=0}^{N-1} \sum_{i=1}^{M} |P_1 a (u^{(i)}(t_n))-h_a(t_n, P_1 u^{(i)}(t_n))|^2,\\
		\hat{\check{b}}^2 =\arg \min\limits_{h_b\in \cV_{\Lambda}} \frac{1}{NM}\sum_{n=0}^{N-1} \sum_{i=1}^{M} |(P_1 b b^T P_1^T) (u^{(i)}(t_n))-h_b( t_n, P_1 u^{(i)}(t_n))|^2.
	\end{cases}\,
\end{equation}

First, in order to ensure the stabilty, we consider \textit{the orthonormalised basis expansion} $\bar{\psi}_p$ via orthonormalising $\psi_p$ 
by the modified Gram-Schmidt algorithm wrt a discrete scalar product for $\forall p, q \in \Lambda (\varpi)$:
\[
<\psi_p, \psi_q> = \frac{1}{NM}\sum_{n=0}^{N-1} \sum_{i=1}^{M}  \psi_p (t_n, P_1 u^{(i)}(t_n))  \psi_q (t_n, P_1 u^{(i)}(t_n)).
\]
The modified Gram-Schmidt algorithm produces a QR decomposition for non-orthonormalised matrix $\Psi$ as follows:
\[
\underbrace{\left[
	\begin{array}{c|c|c|c}
		&    &   &  \\
		\psi_p^1  & \psi_p^2   & ...  &\psi_p^{\# \Lambda}	 \\
		&    &   &  \\
	\end{array}
	\right]}_{\Psi}=
\underbrace{\left[
	\begin{array}{c|c|c|c}
		&    &   &  \\
		\bar{\psi}_p^1  & \bar{\psi}_p^2   & ...  &\bar{\psi}_p^{\# \Lambda}	 \\
		&    &   &  \\
	\end{array}
	\right]}_{Q}
\underbrace{\left[
	\begin{array}{c|c|c|c}
		r_{11}& r_{12}   &   & r_{1\# \Lambda} \\
		& r_{22}   & ...  &r_{2\# \Lambda}	 \\
		&    & ...  &  \\
		&    &   &  r_{\# \Lambda\# \Lambda}\\
	\end{array}
	\right]}_{R}
\]
where $Q^TQ=I$ and $R$ is the upper-triangular matrix. Note that since $Q=\Psi R^{-1}$, we can use the entries of $R^{-1}$ to construct the orthonormalised function $\bar{\psi}_p(t,s)$.

Then, solving~\eqref{eq:L2regression} is equivalent to solving the least squares problem for the expansion coefficients $a_p$ and $b_p$:
\begin{equation}\label{h_fun}
	\begin{cases} 
		\hat{h}_a (t,s) = \sum_{p\in\Lambda}	a_p \bar{\psi}_p (t,s),\\
		\hat{h}_b (t,s) = \sum_{p\in\Lambda}	b_p \bar{\psi}_p (t,s).
	\end{cases}\,
\end{equation} 
The expansion coefficients in~\eqref{h_fun} are computed via solving the respective normal equations for $\underbar{a} = \{a_p\}_{p\in\Lambda}$ and  $\underbar{b} = \{b_p\}_{p\in\Lambda}$:
\begin{equation}
	\begin{cases} 
		Q^TQ \underbar{a} = Q^T \underbar{f} \\
		Q^TQ \underbar{b} = Q^T \underbar{g},\\
		(\underbar{f})_j = P_1 a(t_n, u^{(i)}(t_n)), \\
		(\underbar{g})_j  = P_1 b b^T P_1^T (t_n, u^{(i)}(t_n)),\\
		Q_{jp} = \bar{\psi}_p(t_n, P_1 u^{(i)}(t_n)) \mbox{ for } p\in\Lambda,\; j=1,...,NM.
	\end{cases}\,
\end{equation} 

\section{Multilevel Cross-Entropy method}\label{appx:CEMethod}
We mainly follow the description of the method provided in~\cite{de2005tutorial, kroese2013cross}. In the framework of the EnKF, the optimal original initial density is assumed to be Gaussian. Let $\rho_{0}(x_0) =\phi(x_0; \bu)$ belong to a family of Gaussian densities with a parameter vector $\bu$ which is approximated by the EnKF ensembles. Then, similarly, to~\eqref{eq:optimdensityPropto}, the optimal importance sampling density is
\begin{equation}
	\tilde{\rho}_{0}(x_0) \propto \phi(x_0; \bu) \sqrt{\bP(\max_{0
			\leq k\leq K } P_1 \bar{u}_{n,k}\geq \cK |\bar{u}_{n,0}=x_0)}.
\end{equation}

We want to approximate $\tilde{\rho}_0(x_0)$ by a member of the Gaussian family, $\phi(x_0; \bzeta)$, for some parameters $\bzeta$.
We consider the \textit{Kullback-Leibler divergence} defined as
\begin{equation}
	\begin{split}
		\mathcal{D}(\tilde{\rho}_{0}(\cdot),\phi (\cdot; \bzeta))=\int (\log \tilde{\rho}_{0}(x_{0}) - \log \rho_{0} (x_{0}; \bzeta))\tilde{\rho}_{0}(x_0) dx_0.
	\end{split}
	\label{KLdist}
\end{equation}
This is also known as the cross-entropy between the given densities.

Finding a minimizer of~\eqref{KLdist} is equivalent to finding
\begin{equation}
	\begin{split}
		\bzeta^*&=\argmax_{\bzeta}-\mathcal{D}(\tilde{\rho}_{0}(\cdot), \phi(\cdot; \bzeta)) =  \argmax_{\bzeta} \int \tilde{\rho}_{0}(x_0)  \log \phi (x_{0}; \bzeta) dx_0 \\
		& =\argmax_{\bzeta} \int \phi(x_0,\bu)\sqrt{\bP\Big(\max_{0
				\leq k\leq K } P_1 \bar{u}_{n,k}\geq \cK |\bar{u}_{n,0}=x_0\Big)} \log\phi (x_0; \bzeta)dx_{0}\\
		& = \argmax_{\bzeta} \E_{\bw} \Big[ \sqrt{\bP\Big(\max_{0
				\leq k\leq K } P_1 \bar{u}_{n,k}\geq \cK |\bar{u}_{n,0}=X_0\Big)} \log \phi (X_0; \bzeta) L(X_0; \bu, \bw)  \Big]
	\end{split}
	\label{OptimProblem}
\end{equation}
where the likelihood ratio is defined by
\[
L(x_0;\bu,\bw)=\frac{\rho_{0} (x_0; \bu)}{\phi(x_0; \bw)}.
\]
$\E_{\bw}[\cdot]$ means the expectaction is taken wrt the density $\phi(\cdot; \bw)$, which is, in our setting, another family of Gaussian densities with an arbitrary parameter vector $\bw$. The reason of introducing it is to estimate $\bzeta^{*}$ using random samples $\{X_0^{(j)}\}_{j=1}^J$ from $\phi (\cdot; \bw)$. Then, the approximate optimization becomes
\begin{equation}\label{eq:CEoptimproblem}
\max_{\bzeta} \underbrace{\frac{1}{J} \sum_{j=1}^{J}\sqrt{\bP(\max_{0
			\leq k\leq K } P_1 \bar{u}_{n,k}\geq \cK |\bar{u}_{n,0}=X_0^{(j)})} \log\rho_{0} (X_0^{(j)}; \bzeta) L(X_0^{(j)}; \bu, \bw)}_{-\hat{D}}.
\end{equation}
Assuming that $\hat{D}$ is convex and differentiable with respect to $\bzeta$, we solve the system of equations
\[
\frac{1}{J} \sum_{j=1}^{J}\sqrt{\bP(\max_{0
		\leq k\leq K } P_1 \bar{u}_{n,k}\geq \cK |\bar{u}_{n,0}=X_0^{(j)})} \nabla \log\rho_{0} (X_0^{(j)}; \bzeta) L(X_0^{(j)}; \bu, \bw) =\boldsymbol{0}
\]
where the gradient is with respect to $\bzeta$ and $\boldsymbol{0}$ denotes a zero vector in size of $\bu$. We can often compute the analytical solution for the above problem if the random variable distributions belong to a natural exponential family, which is the case in the EnKF. 

It is crucial to note that the optimization problem~\eqref{eq:CEoptimproblem} can be difficult to carry out due to the rareness of the events $\{\max_{0
	\leq k\leq K } P_1 \bar{u}_{n,k}\geq \cK|\bar{u}_{n,0}=X_0^{(j)} \}$. We can employ a so-called \textit{multilevel} algorithm to overcome the difficulty with the rare event probability~\cite{de2005tutorial}. The main idea is to construct a sequence of reference parameters $\{\bzeta_\ell, \ell\geq 0\}$ and a sequence of threshold levels $\{\cK_\ell, \ell\geq 1\}$ (see Algorithm~\ref{alg:mlCE}). 

\begin{algorithm}[H]\label{alg:mlCE}
	\DontPrintSemicolon
	
	
	{  Define $\hat{\bzeta}_0=\bu$. Set $\ell=1$. 
		
		Generate  $\{\bar{u}_{n,0}^{(i)}\}_{i=1}^{J_1} \sim \phi (x_{0}; \hat{\bzeta}_{\ell-1})$ and $\{\bar{u}_{n,1:K}^{(i)}\}_{i=1}^{J_1} \sim \rho_{\bar{u}_{n,1:K}|\bar{u}_{n,0}^{(i)}} (x_{1:K}|x_0)$. Compute the sample $(1-\beta)-$quantile $\hat{\cK}_\ell$ by calculating $\bar{M}_{\Delta t,(i)}^{n}=\max_{0\leq k\leq K} P_1\bar{u}_{n,k}^{(i)}$ and arranging them in ascending order
		\[
		\bar{M}_{\Delta t,(1)}^{n}\leq ...\leq \bar{M}_{\Delta t,(J_1)}^{n}
		\]
		where $\beta$ is an usually not small number, say $\beta\simeq 10^{-2}$. Evaluate $\hat{\cK}_\ell=\bar{M}_{\Delta t,(\lceil (1-\beta)J_1 \rceil )}^{n}$ if $\hat{\cK}_\ell<\cK$, otherwise $\hat{\cK}_\ell=\cK$.
		
		Given $\bzeta_{\ell-1}$, $\hat{\cK}_\ell$ and the \textbf{same} sample $\bar{u}_{n,0:K}^{(1)},...,\bar{u}_{n,0:K}^{(J_1)}$, find
		\[\small
		\hat{\bzeta}_\ell=\argmax_{\bzeta} \frac{1}{J_1} \sum_{i=1}^{J_1} \sqrt{\bP(\max_{0
				\leq k\leq K } P_1 \bar{u}_{n,k}^{(i)}\geq \cK |\bar{u}_{n,0}^{(i)})}L(\bar{u}_{n,0}^{(i)};\bu,\hat{\bzeta}_{\ell-1})\log\phi(\bar{u}_{n,0}^{(i)}; \bzeta) 
		\]
		
		If $\hat{\cK}_\ell<\cK$, set $\ell=\ell+1$ and go to \textbf{Step 2}, otherwise, go to \textbf{Step 5}.
		
		Estimate the rare event probability $\alpha_n$ by
		\[
		\hat{\alpha}_n = \frac{1}{J_2} \sum_{i=1}^{J_2} \bP(\max_{0
			\leq k\leq K } P_1 \bar{u}_{n,k}^{(i)}\geq \cK|\bar{u}_{n,0}^{(i)}) L(\bar{u}_{n,0}^{(i)};\bu,\hat{\bzeta}_{\cL})
		\]
		where $\cL$ denotes the final number of levels/iterations. 	
	}	
	\caption{Multilevel Cross-Entropy method for Rare-Event Probability}
\end{algorithm}
For example, let $\bzeta=[\tilde{\mu}_0, \tilde{\sigma}_0]$ and $\phi(x_0,\bu)=\frac{1}{\sqrt{2\pi}\sigma_0} e^{-(x_0-\mu_0)^2/2\sigma_0^2}$, where $\bu=[\mu_0, \sigma_0]$. Then, in \textbf{Step 3} of the Algorithm~\ref{alg:mlCE}, we have to solve the following system of equations
\begin{equation}
	\begin{cases}
		\frac{1}{J_1} \sum_{i=1}^{J_1}  \sqrt{\bP(\max_{0
				\leq k\leq K } P_1 \bar{u}_{n,k}^{(i)}\geq \cK |\bar{u}_{n,0}^{(i)})} L(\bar{u}_{n,0}^{(i)};\bu,\hat{\bzeta}_{\ell-1})\Big( \frac{\bar{u}_{n,0}^{(i)}-\tilde{\mu}_0}{\tilde{\sigma}_0^2}\Big) =\boldsymbol{0}\\
		\frac{1}{J_1} \sum_{i=1}^{J_1}  \sqrt{\bP(\max_{0
				\leq k\leq K } P_1 \bar{u}_{n,k}^{(i)}\geq \cK |\bar{u}_{n,0}^{(i)})} L(\bar{u}_{n,0}^{(i)};\bu,\hat{\bzeta}_{\ell-1})\Big( \frac{(\bar{u}_{n,0}^{(i)}-\tilde{\mu}_0)^2}{\tilde{\sigma}_0^3}-\frac{1}{\tilde{\sigma}_0}\Big) =\boldsymbol{0}
	\end{cases}
\end{equation}
which have the following analytical solutions
\begin{equation}\label{eq:CErefparameters}
	\begin{cases}
		\tilde{\mu}_{0,\ell}=\frac{\sum_{i=1}^{J_1}  \sqrt{\bP(\max_{0
					\leq k\leq K } P_1 \bar{u}_{n,k}^{(i)}\geq \cK |\bar{u}_{n,0}^{(i)})} L(\bar{u}_{n,0}^{(i)};\bu,\hat{\bzeta}_{\ell-1}) \bar{u}_{n,0}^{(i)}}{\sum_{i=1}^{J_1}  \sqrt{\bP(\max_{0
					\leq k\leq K } P_1 \bar{u}_{n,k}^{(i)}\geq \cK |\bar{u}_{n,0}^{(i)})} L(\bar{u}_{n,0}^{(i)};\bu,\hat{\bzeta}_{\ell-1})},\\
		\\
		\tilde{\sigma}_{0,\ell}^2=\frac{ \sum_{i=1}^{J_1}  \sqrt{\bP(\max_{0
						\leq k\leq K } P_1 \bar{u}_{n,k}^{(i)}\geq \cK |\bar{u}_{n,0}^{(i)})}L(\bar{u}_{n,0}^{(i)};\bu,\hat{\bzeta}_{\ell-1})\Big( \bar{u}_{n,0}^{(i)}-\tilde{\mu}_0\Big)^2}{\sum_{i=1}^{J_1}  \sqrt{\bP(\max_{0
						\leq k\leq K } P_1 \bar{u}_{n,k}^{(i)}\geq \cK |\bar{u}_{n,0}^{(i)})} L(\bar{u}_{n,0}^{(i)};\bu,\hat{\bzeta}_{\ell-1})}.
	\end{cases}
\end{equation}

Note that the CE method with a Gaussian assumption for the initial condition may lead to an infinite variance of the estimator. Therefore, we observe certain non-robustness of the estimator variance produced by the CE method. To resolve this issue, we propose to fix $\tilde{\sigma}_0=\sigma_0$ and shift only $\tilde{\mu}_0$ by~\eqref{eq:CErefparameters}. $\tilde{\rho}_{0}^{CE}$ corresponds to the optimal IS density obtained via the multilevel CE method with the shifted mean and fixed variance.

\section{Numerical PDE solver for the problem with a discontinuity at the boundary}\label{appx:PDEsolver}
Here, we briefly describe the numerical method used to approximate the solution to the auxiliary initial-boundary value problem of the type 
\begin{subequations}\label{eq:BKEpde0T}
\begin{equation}\label{eq:BKEpde0Ta}
		\frac{\partial \gamma}{\partial t}  = -a(x,t)	\frac{\partial \gamma}{\partial x}   -\frac{1}{2} b(x,t)b(x,t)^T	\frac{\partial^2 \gamma}{\partial x}, \quad (x,t)\in (-\infty, \cK) \times (0,T)
\end{equation}
\begin{equation}
		\left\{
	\begin{array}{ll}\label{eq:BKEpde0Tb}
		\gamma(x, T) =0, \quad \quad \;\; \, x<\cK,\\
		\gamma(\cK, t) =1, \quad \qquad t\leq T,\\
		\lim\limits_{x\rightarrow -\infty} \gamma(x,t)=0, \;\, t\leq T.
	\end{array}
	\right.
\end{equation}
\end{subequations}
Compare~\eqref{eq:BKEpde} and~\eqref{eq:BKEpdeSurrogate}.

Since the accuracy of the approximation to~\eqref{eq:BKEpde0T} does not directly affect the accuracy of the estimate of the rare event probability $\alpha_n$, only the efficiency of the proposed IS algorithm, and since even a qualitatively correct approximation can give us a useful IS we may use a simple method here. In our numerical tests, we used the Crank-Nickolson scheme with the uniform discretization of a truncated domain $(x_{\min},\cK)\times (0,T)$ with an artificial boundary condition $\frac{\partial^2\gamma(x,t)}{\partial x^2}=0$ at $x=x_{\min}$ approximated from known neighboring values by a linear extrapolation. However, the jump discontinuity in the corner $(\cK, T)$ requires special care with this very basic approach. We outline our treatment here. 

The constant coefficient version, $a(x,t)\equiv \mathbf{a}$, $b(x,t)\equiv \mathbf{b}$ of~\eqref{eq:BKEpde0T} has the closed-form solution
\begin{equation}\label{eq:analsol}
	\begin{split}
		\gamma_{\textbf{const}}(x,t)=&\frac{1}{2}\mathrm{erfc}\Bigg(\frac{\cK-x-\mathbf{a} (T-t)}{2\sqrt{\frac{\mathbf{b}^2}{2}(T-t)}}\Bigg)\\
		&+\frac{1}{2}e^{\frac{\mathbf{a}(\cK-x)}{\mathbf{b}^2/2}}\mathrm{erfc}\Bigg(\frac{\cK-x+\mathbf{a} (T-t)}{2\sqrt{\frac{\mathbf{b}^2}{2}(T-t)}}\Bigg),
	\end{split}
\end{equation}
where $\mathrm{erfc}(z):=\frac{2}{\sqrt{\pi}}\int_{z}^{\infty} e^{-t^2}dt$ is a complementary error function; see e.g.,~\cite{harleman1963longitudinal}. Let $\gamma_{\textbf{fr}}(x,t)$ be the solution to the constant coefficient problem with the frozen coefficients $\mathbf{a}=a(\cK,T)$ and $\mathbf{b}=b(\cK,T)$, i.e., 
\begin{equation}\label{eq:freezcoeffsol}
	\footnotesize
	\begin{split}
		\gamma_{\textbf{fr}}(x,t)=&\frac{1}{2}\mathrm{erfc}\Bigg(\frac{\cK-x-a(\cK,T) (T-t)}{\abs{b(\cK,T)}\sqrt{2(T-t)}}\Bigg)\\
		&+\frac{1}{2}e^{\frac{a(\cK,T)(\cK-x)}{b^2(\cK,T)/2}}\mathrm{erfc}\Bigg(\frac{\cK-x+a(\cK,T) (T-t)}{\abs{b(\cK, T)}\sqrt{2(T-t)}}\Bigg).
	\end{split}
\end{equation}	

We divide the computational domain into three regions, see Figure~\ref{fig:pdesolveridea}, and for appropriately chosen $(\Delta x_{\textbf{PDE}}, \Delta t_{\textbf{PDE}})$ and $(\Delta x_{\textbf{fr}}, \Delta t_{\textbf{fr}})$ approximate $\gamma$ as follows:

\begin{figure}[h!]
	\includegraphics[width=0.7\linewidth]{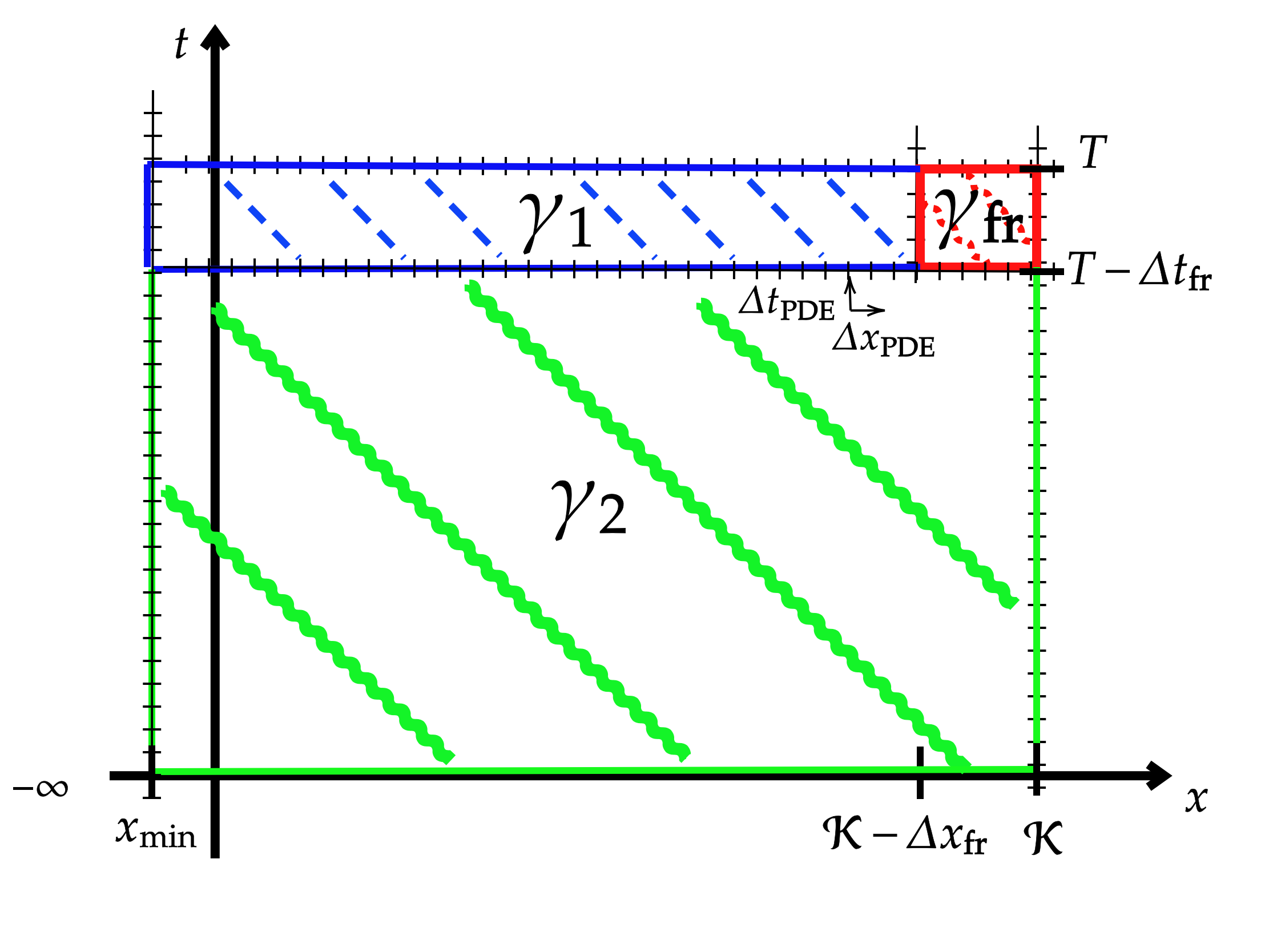}
	\caption{Division of the problem domain for the PDE solver} 
	\label{fig:pdesolveridea}
\end{figure}

\textbf{I.} In the rectangle $[\cK-\Delta x_{\textbf{fr}}, \cK]\times [T-\Delta t_{\textbf{fr}}, T]$, we approximate $\gamma(x,t)$ by $\gamma_{\textbf{fr}}(x,t)$. Observe that $\gamma_{\textbf{fr}}$ decays very fast toward $0$ as $x$ decreases. 

\textbf{II.} In $[x_{\min}, \cK-\Delta x_{\textbf{fr}}) \times [T-\Delta t_{\textbf{fr}}, T]$, we use the finite difference approximation of~\eqref{eq:BKEpde0Ta} with the initial-boundary conditions
\begin{equation*}
	\left\{
	\begin{array}{ll}
		\gamma (x,T)=0, \qquad \qquad \qquad \qquad \qquad \qquad \qquad \qquad \qquad \quad  \quad  x\in (x_{\min}, \cK-\Delta x_{\textbf{fr}}),\\
		\gamma (\cK-\Delta x_{\textbf{fr}}, t)=\gamma_{\textbf{fr}}(\cK-\Delta x_{\textbf{fr}}, t), \qquad \qquad \qquad \qquad \qquad \qquad \; \, t\in [T-\Delta t_{\textbf{fr}}, T], \\
        \frac{\partial^2\gamma(x_{\min},t)}{\partial x^2}=0, \qquad \qquad \qquad \qquad \qquad \qquad \qquad \qquad \qquad \qquad \, \; \; t\in[T-\Delta t_{\textbf{fr}}, T].
	\end{array}
	\right.
\end{equation*}
Denote this approximation $\gamma_1$. 

\textbf{III.} In $[x_{\min},\cK]\times [0,T-\Delta t_{\textbf{fr}}]$, we use the same finite difference approximation of~\eqref{eq:BKEpde0Ta} with the initial-boundary conditions

\begin{equation*}
	\begin{split}
		&\gamma(x,T-\Delta t_{\textbf{fr}}) = \left\{
		\begin{array}{ll}
			\gamma_1(x,T-\Delta t_{\textbf{fr}}), \;\; x\in (x_{\min}, \cK-\Delta x_{\textbf{fr}}),\\
			\gamma_{\textbf{fr}}(x,T-\Delta t_{\textbf{fr}}), \; x\in [\cK-\Delta x_{\textbf{fr}}, \cK),
		\end{array}
		\right.\\
		&\gamma(\cK, t) =1, \qquad  \quad \; \, t\in [0,T-\Delta t_{\textbf{fr}}],\\
	   &\frac{\partial^2\gamma(x_{\min},t)}{\partial x^2}=0, \quad t\in [0,T-\Delta t_{\textbf{fr}}].
	\end{split}
\end{equation*}
Denote this approximation $\gamma_2$. 

We patch together the numerical approximation from the three subdomains and denote the PDE solution obtained via this strategy by $\gamma^{\textbf{PDE}}_{\textbf{fr}}(x,t)$ and the corresponding control by $\xi^{\textbf{PDE}}_{\textbf{fr}}(x,t)$. 

The appropriate choice of the pair $(\Delta x_{\mathbf{fr}}, \Delta t_{\mathbf{fr}})$ depends on the discretization parameters $(\Delta x_{\mathbf{PDE}}, \Delta t_{\mathbf{PDE}})$ of the finite-difference scheme for $\gamma_1$ and $\gamma_2$. One can choose $(\Delta x_{\mathbf{fr}}, \Delta t_{\mathbf{fr}})$ such that the gradient of the closed-form solution at $(\cK-\Delta x_{\mathbf{fr}}, T-\Delta t_{\mathbf{fr}})$ should closely match with $(\Delta x_{\mathbf{PDE}}, \Delta t_{\mathbf{PDE}})$. In the numerical tests in Section~\ref{sec:numerics}, we use $\Delta x_{\mathbf{PDE}}=0.005$, $\Delta t_{\mathbf{PDE}}=\frac{\Delta x_{\mathbf{PDE}}}{2}$ with $x_{\min}=-5$ in the Double Well problem, and $\Delta x_{\mathbf{PDE}}=0.006$, $\Delta t_{\mathbf{PDE}}=\frac{\Delta x_{\mathbf{PDE}}}{5}$ with $x_{\min}=-3$ in the Langevin problem. Correspondingly, the choice $(\Delta x_{\mathbf{fr}}, \Delta t_{\mathbf{fr}})$ used to define the red region in Figure~\ref{fig:pdesolveridea} is provided in Table~\ref{table:redregion} for both problems given the different thresholds $\cK$.
\begin{table}[]
	\begin{tabular}{|c|cc|llll}
		\cline{1-3} \cline{5-7}
		\multirow{2}{*}{$\cK$} & \multicolumn{2}{c|}{\begin{tabular}[c]{@{}c@{}}\textbf{Double Well}\\\end{tabular}} & \multicolumn{1}{l|}{} & \multicolumn{1}{l|}{\multirow{2}{*}{$\cK$}} & \multicolumn{2}{l|}{\begin{tabular}[c]{@{}l@{}}\textbf{Langevin}\\ \end{tabular}} \\ \cline{2-3} \cline{6-7} 
		& \multicolumn{1}{c|}{$\Delta x_{\mathbf{fr}}$}                                 & $\Delta t_{\mathbf{fr}}$                              & \multicolumn{1}{l|}{} & \multicolumn{1}{l|}{}                   & \multicolumn{1}{l|}{$\Delta x_{\mathbf{fr}}$}                        & \multicolumn{1}{l|}{$\Delta t_{\mathbf{fr}}$}                        \\ \cline{1-3} \cline{5-7} 
		\textbf{0}                  & \multicolumn{1}{c|}{0.185}                                   & 0.0225                                 & \multicolumn{1}{l|}{} & \multicolumn{1}{l|}{\textbf{2}}                  & \multicolumn{1}{l|}{0.170}                        & \multicolumn{1}{l|}{0.00625}                      \\ \cline{1-3} \cline{5-7} 
		\textbf{0.5}                & \multicolumn{1}{c|}{0.195}                                   & 0.0250                                 & \multicolumn{1}{l|}{} & \multicolumn{1}{l|}{\textbf{2.5}}                & \multicolumn{1}{l|}{0.121}                        & \multicolumn{1}{l|}{0.0025}                       \\ \cline{1-3} \cline{5-7} 
		\textbf{1 }                 & \multicolumn{1}{c|}{0.185}                                   & 0.0225                                 & \multicolumn{1}{l|}{} & \multicolumn{1}{l|}{\textbf{3}}                  & \multicolumn{1}{l|}{0.120}                        & \multicolumn{1}{l|}{0.0025}                       \\ \cline{1-3} \cline{5-7} 
		\textbf{1.2 }               & \multicolumn{1}{c|}{0.185}                                   & 0.0225                                 & \multicolumn{1}{l|}{} & \multicolumn{1}{l|}{\textbf{3.5}}                & \multicolumn{1}{l|}{0.1365}                       & \multicolumn{1}{l|}{0.00375}                      \\ \cline{1-3} \cline{5-7} 
		\textbf{1.5 }               & \multicolumn{1}{c|}{0.185}                                   & 0.0225                                 &                       & \multicolumn{3}{l}{\multirow{4}{*}{}}                                                                                                           \\ \cline{1-3}
		\textbf{2  }                & \multicolumn{1}{c|}{0.185}                                   & 0.0225                                 &                       & \multicolumn{3}{l}{}                                                                                                                            \\ \cline{1-3}
		\textbf{2.5 }             & \multicolumn{1}{c|}{0.200}                                   & 0.0275                                 &                       & \multicolumn{3}{l}{}                                                                                                                            \\ \cline{1-3}
		\textbf{3}                  & \multicolumn{1}{c|}{0.195}                                   & 0.0250                                 &                       & \multicolumn{3}{l}{}                                                                                                                            \\ \cline{1-3}
	\end{tabular}
\caption{The choice for $(\Delta x_{\mathbf{fr}}, \Delta t_{\mathbf{fr}})$ used in the numerical PDE solver for the problems considered in Section~\ref{sec:numerics}.}
\label{table:redregion}
\end{table}

\section{Adaptive time-stepping scheme for IS simulation}\label{appx:AdaptiveScheme}

The optimal control given by~\eqref{eq:optcontrol1d} blows up as $t\rightarrow T$. For this reason, it is not appropriate to use uniform time discretizations in the Euler-Maruyama approximation of the controlled SDE. Our numerical experiments used a simple adaptive time-stepping scheme based on the solution $\gamma_{\textbf{fr}}$ to the frozen coefficient problem. The control corresponding to $\gamma_{\textbf{fr}}$ is 
\begin{equation}\label{eq:optcontfrozencoeff}
	\footnotesize
	\begin{split}
		\xi_{\textbf{fr}}(x,t)=& b(\cK,T) \frac{\partial (\log(\gamma_{\textbf{fr}}(x,t)))}{\partial x}= \frac{b(\cK,T) }{\gamma_{\textbf{fr}}(x,t)}\frac{\partial (\gamma_{\textbf{fr}}(x,t))}{\partial x}\\
		&=\frac{b(\cK,T) }{\gamma_{\textbf{fr}}(x,t)} \Bigg[\frac{1}{\abs{b(\cK,T) }\sqrt{2\pi (T-t)}}\Bigg(e^{-(R^{-}(x))^2}+e^{-(R^{+}(x))^2+\frac{2a(\cK,T) (\cK-x)}{b(\cK,T) ^2}}\Bigg)\\
		&-\frac{a(\cK,T) }{b(\cK,T) ^2}e^{\frac{2a(\cK,T) (\cK-x)}{b(\cK,T) ^2}}\mathbf{erfc}(R^{+}(x))\Bigg],
	\end{split}
\end{equation}	
where $R^{\pm}(x):=\frac{\cK-x\pm a(\cK,T) (T-t)}{\abs{b(\cK,T) }\sqrt{2(T-t)}}$. 

Asymptotically, as $t\rightarrow T$, for fixed $x$ 
\[
\xi_{\textbf{fr}}(x,t)\cong \xi_{\textbf{fr}}^{\infty}(x,t):=\frac{1}{b(\cK, T)} \Big(\frac{\cK-x}{T-t}\Big)-a(\cK, T).
\]

We extend this asymptotic formula for the control to non-constant coefficients as 
\[
\xi^{\infty}(x,t)=\frac{1}{b(x, t)} \Big(\frac{K-x}{T-t}\Big)-a(x, t).
\]
and use the control $\xi^{\infty}(x,t)$ for $t\in (T-\delta t, T]$, while we use $\xi^{\textbf{PDE}}_{\textbf{fr}}(x,t)$, based on the numerical approximation $\gamma^{\textbf{PDE}}_{\textbf{fr}}(x,t)$, for $t\in [0, T-\delta t]$.

In Algorithm~\ref{alg:adaptivetimestep}, below, we provide a sketch of the adaptive time-stepping scheme with the Brownian bridge technique used in Monte Carlo simulations of all methods (Crude MC, IS wrt $\rho_0$, IS wrt $W(t)$, IS wrt both $\rho_0$ and $W(t)$ and IS with CE) considered in this work. 

\begin{algorithm}[H]
	\DontPrintSemicolon
	
	\KwInput{ \small $T$, $K$, $\cK$, model parameters, $\delta t$ and $\varepsilon$.}
	Define the uniform Euler-Maruyama timestep: $\Delta t=\frac{T}{K}$; 
	
	Set the optimal control such that
	\begin{equation}
		\small
		\xi(x,t)= 	\left\{
		\begin{array}{ll}
			\xi_{\textbf{fr}}^{\textbf{PDE}}(x,t),   \mbox{   if  } t \in [0, T-\delta t],\\
			\xi^{\infty}(x,t), \quad \mbox{  if  }  t \in (T-\delta t, T];
		\end{array}
		\right.
	\end{equation}

	$k=0$, $t_{n,k}=0$, $\mathrm{h} = \Delta t $;
	
	Set the initial condition $\bar{u}_{n,0}$ according to each method;
	
	\While{$(t_{n,k}\leq T)$, $(P_1\bar{u}_{n,k}<\cK)$ and $(\mathrm{h}>\varepsilon)$}
	{   \textit{Adaptive rule:} $\mathrm{h}= \min (\Delta t, \frac{T-t_{n,k}}{2})$;
		
		$\Delta W^{\bP}_{n,k}\sim N(0,\mathrm{h}\mathbb{I}_{d_W})$;
		
		Simulate the trajectory 
		$\bar{u}_{n,k+1}={\bar{u}_{n,k}+a(\bar{u}_{n,k})\mathrm{h}+b(\bar{u}_{n,k})\xi(\bar{u}_{n,k},t_{n,k})\mathrm{h}+b(\bar{u}_{n,k})\Delta W^{\bP}_{n,k}}^{\red{\ast}};$
		
		Compute the corresponding likelihood according to each method;
		
		\textit{Hitting test via Brownian bridge technique}:
		
		\nonumber{compute the exiting probability $q=\exp\prt{-\frac{2\max(\cK-P_1\bar{u}_{n,k},0)\max(\cK-P_1 \bar{u}_{n,k+1},0)}{(P_1b(\bar{u}_{n,k})b(\bar{u}_{n,k})^TP_1^T)\mathrm{h}}}$;}
		
		sample a uniform random variable $r \sim U(0,1)$;
		
	    \eIf {$r<q$} {set  $P_1 \bar{u}_{n,k+1}=\cK$ and the stopping time
		$\tau=t_{n,k}+0.5\mathrm{h}$;}{$t_{n,k}=t_{n,k}+\mathrm{h}$;}
		$k=k+1$;			
	}
	\caption{Adaptive time-stepping scheme with Brownian bridge technique}
	\label{alg:adaptivetimestep}
\end{algorithm}
$^{\red{\ast}}$ $\xi(x,t)=0$ in the crude MC, IS wrt $\rho_0$ and CE simulations.

In both numerical tests in Section~\ref{sec:numerics}, we set $\delta t=10 \Delta t$ and $\varepsilon=10^{-6}$ as this choice was sufficient to ensure a high level of variance reduction for each method. 
\end{justify}
\medskip

{\bf Acknowledgments } This work was supported by the KAUST Office of
Sponsored Research (OSR) under Award No. URF/1/2584-01-01 and the
Alexander von Humboldt Foundation. E.~von Schwerin, G.~Shaimerdenova and R.~Tempone are
members of the KAUST SRI Center for Uncertainty Quantification in
Computational Science and Engineering.

For the purpose of
open access, the authors have applied a Creative Commons Attribution (CC BY) licence to any
Author Accepted Manuscript version arising from this submission.
\bibliography{references}
\bibliographystyle{plain}

\end{document}